\documentclass[12pt, a4paper, final, oneside, onecolumn, notitlepage]{report}

\usepackage{graphicx}
\usepackage{amsmath}
\usepackage{amsfonts}
\usepackage{amssymb}
\usepackage{amsthm}
\usepackage[utf8]{inputenc}
\usepackage[usenames]{color}
\usepackage{enumerate}
\usepackage{stackengine}
\usepackage{accents}
\usepackage{relsize}
\usepackage{mathtools}
\usepackage{textgreek}
\usepackage{upgreek}
\usepackage{ulem}
\usepackage[T1]{fontenc}
\usepackage{dsfont}
\usepackage{bbold}
\usepackage{mathrsfs}
\usepackage{yfonts}
\usepackage{stmaryrd}
\usepackage{array}
\usepackage{stackrel}
\stackMath

\newtheorem{thm}{\small\bf Theorem}[section]

\newtheorem{lem}[thm]{\small\bf Lemma}
\newtheorem{cor}[thm]{\small\bf Corollary}
\newtheorem{re}[thm]{\small\bf Remark}
\newtheorem{df}[thm]{\small\bf Definition}

\numberwithin{equation}{section}

\newcommand{\be}{\begin{eqnarray}}
\newcommand{\ee}{\end{eqnarray}}
\newcommand{\bd}{\begin{displaymath}}
\newcommand{\ed}{\end{displaymath}}
\newcommand{\bee}{\begin{eqnarray*}}
\newcommand{\eee}{\end{eqnarray*}}

\newcommand{\R}{\mathbb{R}}
\newcommand{\N}{\mathbb{N}}
\newcommand{\Z}{\mathbb{Z}}

\newcommand{\T}{\mathcal{T}}
\newcommand{\TT}{\mathbb{T}}

\newcommand{\XX}{\mathbb{X}}
\newcommand{\x}{\mathbf{x}}

\newcommand{\ttt}{\mathbf{t}}

\newcommand{\rrr}{\mathbf{r}}
\newcommand{\sss}{\mathbf{s}}

\newcommand{\aae}{\textsl{\ae}}
\newcommand\tsup[2][2]{%
 \def\useanchorwidth{T}%
  \ifnum#1>1%
    \stackon[-1.3ex]{\tsup[\numexpr#1-1\relax]{#2}}{\mathchar"307E}%
  \else%
    \stackon[-1ex]{#2}{\mathchar"307E}%
  \fi%
}

\newenvironment{namelist}[1]{%
\begin{list}{}
     {
      
      \settowidth{\labelwidth}{#1}
      \setlength{\leftmargin}{1.1\labelwidth}
               }
      }{%
\end{list}}
\def\rightsq{\urcorner\!\!\!\lrcorner}
\def\leftsq{\ulcorner\!\!\!\llcorner}

\newlength{\dhatheight}
\newcommand{\doublehat}[1]{%
	\settoheight{\dhatheight}{\ensuremath{\hat{#1}}}%
	\addtolength{\dhatheight}{-0.35ex}%
	\hat{\vphantom{\rule{1pt}{\dhatheight}}%
		\smash{\hat{#1}}}}

\newcommand{\bs}[1]{\boldsymbol{#1}}
\newcommand{\ssucc}{\mathbin{\succ\kern-.3em\succ}}
\newcommand{\pprec}{\mathbin{\prec\kern-.3em\prec}}
\newcommand{\clo}{\operatorname{clos}}
\newcommand{\sm}{\setminus}
\newcommand{\eps}{\varepsilon}
\newcommand{\bbt}{\mathbb{T}}

\usepackage{fancyhdr}
\begin{titlepage}
\title{The Convex Peano Curve Does Exist}
\author{Adam Paszkiewicz}
\date{}
\end{titlepage}

\renewcommand{\maketitle}{

\bibliographystyle{plain}

\newcommand{\linia}{\rule{\linewidth}{0.4mm}}
    \vspace*{1cm}
    \begin{center}\small
		\textsc{ }\\
		\vspace{1.5cm}
    \end{center}
    \vspace{2cm}
    \noindent\linia
    \begin{center}
     \LARGE \textsc{The Convex Peano Curve Does Exist} \\

         \end{center}
     \linia
    \vspace{0.5cm}
    	\begin{center}
					\large \textsc{Adam Paszkiewicz}
					\\
			\end{center}
		\vspace{1cm}
		\vspace{2cm}

    \vspace*{\stretch{6}}
    \begin{center}
    \textsc{Łódź}, 2024.06.12.
    \end{center}
}

\begin{document}

\maketitle
\thispagestyle{empty}

\newpage
\thispagestyle{empty}
\tableofcontents

\newpage
\underline{Abstract.} We refer here to the surprising construction made by Giu\-sep\-pe Peano in 1890. He gave an
example of a continuous function (called now the Peano curve) from the unit interval to the whole unit square.
We show here the existence of a more general space-filling curve with additional convexity properties. More
precisely: by $\mathbb{T} \subset \mathbb{R}^2$ we denote a convex closed and bounded set and we show that there
exists a continuous surjection $f \colon [0,1] \mapsto \mathbb{T}$ for which the image of any interval is a convex set.


\underline{Classification:} 52A10, 26E40, 46G12

\newpage
\chapter*{1. Introduction}
\addcontentsline{toc}{chapter}{Introduction}
\stepcounter{chapter}

\section{Main result}

The Peano curve is the first example of a space-filling curve which was discovered by Giuseppe Peano
in 1890 \cite{Peano}.  Peano's curve is a surjective, continuous function from the unit interval onto the unit square.
Peano was motivated by the result of Georg Cantor stating that these two sets have the same cardinality. Some authors
use the name space-filling curve for the Peano curve and some of its generalizations \cite{Gugenheimer}.

\begin{df}
\label{df1.1.1}
The space-filling curve $f \colon [0,1] \mapsto \mathbb{T} \subset \mathbb{R}^2$ which is a continuous surjection
such that the image of any subinterval of $[0,1]$ (closed or not) is a convex set we call the \underline{convex Peano curve}.
\end{df}

\noindent
In this paper we present the following main result: \\

\begingroup
\setcounter{thm}{0}
\renewcommand\thethm{\Alph{thm}}
\begin{thm}
\label{thm0}
For every convex, closed, bounded non-empty set $\mathbb{T} \subset \mathbb{R}^2$ we can find a ''convex Peano
curve $f$ filing $\mathbb{T}$''. If the set $\mathbb{T}$ has non-empty interior, then the $f$ curve can additionally
transform non-degenerate to a point interval into a convex set with non-empty interior.
\end{thm}
\endgroup
\setcounter{thm}{0}

\section{Relationships with the results known so far}

The problem of the existence of a convex Peano curve has been known since the 80s of the last century. J. Pach and
C.A. Rogers in the paper \cite{PachRogers} stated that this problem is known to them as a question posed by
M. Michalik and A. Wieczorek from the Polish Academy of Sciences.

The partial result obtained so far is as follows: for some set $\mathbb{T}$ with non-empty interior there
is a continuous surjection $f \colon [0,1] \mapsto \mathbb{T}$, for which the images $f([0 , \alpha])$ are
convex for every $\alpha \in [0,1]$. Independently of each other, such constructions were obtained by J. Pach and C.A.
Rogers in paper \cite{PachRogers} when $\mathbb{T}$ is a square and A. Vince and D.C. Wilson in paper \cite{VinceWilson}
when $\mathbb{T}$ is a ball.

As a corollary, we obtain that there is a continuous surjection $f_1 \colon [0,1]$ $\mapsto \mathbb{T}$ for
which the images $f_1([0,\alpha])$ and also $f_1([\beta, 1])$ are convex for all $\alpha, \beta \in [0,1]$. Just define
$$
f_1(u) = \left\{ \begin{array}{lcl}
      f(2u) & \hbox{for} & u \in [0,\frac{1}{2}], \\
      f(2-2u) & \hbox{for} & u \in (\frac{1}{2},1]. \end{array} \right.
$$
This has been noted in the work \cite{PachRogers}.

The literature on Peano curves is very rich. We refer here to monograph \cite{Sagan}. More recent results
on the richness of such curves can be found in \cite{Albuquerque+all}, \cite{Bartoszewicz+all} and in the papers cited there. We will note here
only a few known results related to the role of Theorem \ref{thm0} and the natural directions of further research.

The natural and in some sense the strongest generalization of G. Peano's construction is the classical result known
as the Hahn-Ma\-zur\-kie\-wicz theorem: A topological Hausdorff space $Y$ is an image for a continuous surjection defined
on $[0,1]$ if and only if $Y$ is compact, connected, locally connected and it satisfies the second axiom of countability
(possess a countable topological base) . Such $Y$ spaces are called Peano spaces (cf. \cite{Willard} Theorem 31.5).

Call any subset $\mathbb{C}$ of a topological linear space {\em a convex Peano space} if there is a continuous
surjection $f \colon [0,1] \mapsto \mathbb{C}$ such that the image of any interval is a convex set. Describing convex
Peano spaces (starting with $\mathbb{R}^n$ subsets) is a natural direction for further research.

Recently in paper \cite{BochiMilet}, J. Bochi and P.H. Milet constructed a continuous surjection $f \colon [0,\infty) \mapsto
\mathbb{R}^2$ that satisfies $f(0) = \mathbf{0}$ and the following conditions: For every $\alpha > 0$, \\
--- the image of $f([0,\alpha])$ is a convex set in $\mathbb{R}^2$, \\
--- its boundary $\mathbb{C}_{\alpha}$ is the image $\varphi_{\alpha} ([0,2 \pi))$ for some function  $\varphi_{\alpha} \colon
[0, 2\pi)$ $\mapsto \mathbb{C}_{\alpha} \subset \mathbb{R}^2$, \\
--- the function $\varphi_{\alpha}$ is a bijection of the $C_{\infty}$ class, \\
--- the functions $\varphi_{\alpha}$ and  $\psi_{\alpha}(\theta) := (\alpha \cos \theta, \alpha \sin \theta)$, $\theta \in [0, 2\pi)$,
we can make as closed in the norm $C_{\infty }$ as we  want.

This result reinforces the results obtained in \cite{PachRogers}, \cite{VinceWilson}. \\

Additionally, Bochi and Milet obtained a surprising property: $f(\alpha) \in \mathbb{C}_{\alpha}$ for $\alpha \in [0,1]$.
They also showed that \\
--- there is a continuous mapping $\mathbb{R}^2 \setminus \{0\} \ni \mathbf{x} \mapsto j(\mathbf{x})$ into the
   set of lines in $\mathbb{R}^2$ such that for $\mathbf{x} \in \mathbb{C}_{\alpha}$ the line $j(\mathbf{x})$ is tangent to
   $\mathbb{C}_{\alpha}$ at the point $\mathbf{x}$; \\
--- the map $\mathbf{x} \mapsto j(\mathbf{x})$ is nowhere differentiable. \\

This result suggests the possibility of examining the maximum regularity of the boundaries of the sets $f([\alpha, \beta])$,
$0 \leqslant \alpha \leqslant \beta \leqslant 1$, for the convex Peano curve appearing in Theorem A. It remains an open
problem the existence of a convexity-preserving surjection $f \colon [0,\infty) \mapsto \mathbb{R}^2$ (or equivalently
$f \colon [0,1) \mapsto \mathbb{R}^2$ ) filling the entire $\mathbb{R}^2$
space. We believe that such a mapping exists.

\medskip

Another classical (and rather simple) result is to find the maximal H\"{o}lder exponent of the Peano curve
$f \colon [0,1] \mapsto [0,1] \times [0,1]$. Already the classical construction of Peano shows that $f$
can have the H\"{o}lder exponent $\frac{1}{2}$ (see e.g. \cite{Kennedy}). There is no analogous $f$ function
that has an exponent greater than $\frac{1}{2}$ (closer to the $1$-exponent of the Lipschitz mapping), cf.
\cite{Falkoner}, Proposition 2.3. We believe that there are convex Peano curves with any H\"{o}lder exponent
less than $\frac{1}{4}$ and no such curves exist for any exponent greater than $\frac{1}{4}$. However, there are further
technical difficulties here.

\medskip

The existence of the convex Peano curve also has the following, we believe, interesting consequences in the
theory of linear topological spaces. Let $X$ and $Y$ be linear spaces, and let $D \subset X$ be a convex set. $f \colon
D \mapsto Y$ is called an \textit{affine mapping} if for any $\mathbf{x}_1, \mathbf{x}_2 \in D$ and $\lambda \in [0,1]$ we
have $ f(\lambda \mathbf{x}_1 + (1-\lambda) \mathbf{x}_2) = \lambda f(\mathbf{x}_1) + (1-\lambda) f(\mathbf{x} _2)$. Whereas
$f$ is a \textit{convex mapping} if only the images $\{f(\lambda \mathbf{x}_1 + (1-\lambda) \mathbf{x}_2) \colon \lambda \in [0,1]\}$
are convex sets for any $\mathbf{x}_1, \mathbf{x}_2 \in D$. A straight line in a linear space is any one-dimensional affine subspace
$\{ \lambda f(\mathbf{y}_1) + (1-\lambda) f(\mathbf{y}_2) \colon \lambda \in \mathbb{R} \}$, for $\mathbf{y}_1, \mathbf{y}_2
\in X$, $\mathbf{y}_1 \neq \mathbf{y}_2$. Sch\"{o}pf's theorem (see \cite{Schopf} states that for linear topological
spaces $X$ and $Y$ the convex and G\^{a}teaux-differentiable mapping $f \colon X \mapsto Y$ such that the image $f(X)$ is
not contained in any straight line is already an affine map.

\medskip

Studying possible strengthening of Sch\"{o}pf's theorem (e.g. by weakening the differentiability condition of $f$) is included,
among others, in paper \cite{Horbaczewska}. The following interesting tool for studying linear topological spaces $X$ and $Y $
is shown there. For the linearly topological spaces $X$ and $Y$ and the map $f \colon X \rightarrow Y$ whose image is not
contained in any straight line, the following implication holds: if $f$ preserves the collinearity of the points, then $f$ is an
affine map (see \cite{Horbaczewska}, Theorem 2).

\medskip

For the linearly topological space $X$, the locally convex space $Y$, and the continuous and convex function $f \colon X
\rightarrow Y$, consider the functions:
$$
f_{\Lambda_1 \Lambda_2} (\mathbf{x}) := (\Lambda_1 (f(\mathbf{x})) , \Lambda_2 (f(\mathbf{x}))),
$$
for arbitrary linear continuous functionals $\Lambda_1, \Lambda_2 \in Y^{\ast}$. If there was no convex Peano curve, then all the
$f_{\Lambda_1, \Lambda_2}$ functions (and consequently the $f$ function itself) would have to maintain collinearity of the
points; in particular every continuous convex function $f \colon X \mapsto Y$ would automatically be affine. This is not the
case here, our Theorem \ref{thm0} is the counterexample. It can also be seen that each reinforcement of Theorem \ref{thm0} defines a whole collection
of impossible reinforcements of Sch\"{o}pf's Theorem.

For example, for any $m,n \in \mathbb{N}$, $n \geqslant 2$, and a convex continuous function $f \colon \mathbb{R}^m \rightarrow
\mathbb{R}^n$ for which $f(\mathbb{R}^m)$ is not contained in any line does not have to be affine. Based on Theorem A, we will
prove this at the end of Chapter 3 (Corollary 3.4.1).

\medskip

Another classic problem that we will mention here is the characterization of all Borel measures on topological spaces,
which are images of the Lebesgue measure on $[0,1]$ under continuous mappings. So the problem is to describe all Peano
measures defined as $\nu(B) = \lambda(f^{-1}(B))$, where $f\colon [0,1] \mapsto X$ is a continuous function, $X$ is a
topological space, $B \in \mathcal{B}(X)$ a Borel set, and $\lambda$ is the Lebesgue measure on $\mathbb{R}$. Research
on such measures was initiated by H. Steinhaus in \cite{Steinhaus} in 1936 and discovered again by A.M. Garsia in \cite{Garsia} in
1976 for $X = \mathbb{R}^2$. The final answer when $X$ is a Hausdorff space is given by the following theorem by A.V. Kolesnikov (see
\cite{Kolesnikov} and also \cite{Bogachev} vol. 2, sect. 9.7 and the bibliographic notes in \cite{Bogachev}): \\
{\it $\nu$ is a Peano measure if and only if it is a probability Radon measure having a support that is connected, locally connected,
compact and metrizable $($i.e., the induced topology is metrizable$)$.}

Similarly, when $f \colon [0,1] \mapsto \mathbb{R}^2$ is a convex Peano curve, then $\nu$ we call a convex Peano measure.
We will deal with these issues in the next publication.

\section{Some more important notations.}

In the paper, we use many typical symbols: $\mathbb{N}, \mathbb{Z}$ and $\mathbb{R}$ respectively denote sets of:
natural numbers (without zero), integers and real numbers. By $\| \mathbf{x}\|$ we denote the Euclidean
norm of the vector $\mathbf{x} \in \mathbb{R}^2$, $\Gamma$ is the Cantor set, $2^S$
denotes the class of all subsets of the set $S$. By $\hbox{clos}(S)$, $\hbox{conv}(S)$,
$\partial S$, $\hbox{int}(S)$, $d(S)$ we denote respectively closure, convex hull, boundary, interior and diameter of the set
$S$. For simplicity $\hbox{conv}(\{\mathbf{x}, \mathbf{y}\}) \equiv [\mathbf{x}, \mathbf{y}]$ is the interval connecting points
$\mathbf{x}$ and $\mathbf{y}$, and $[\x, \x] = \{ \x \}$. For algebraic operations, we adopt notations
of the following type:
$$
\begin{array}{lcl}
a \cdot S & = & \left\{ as \colon s \in S \right\}, \\
a+ S & = & \left\{ a+s \colon s \in S \right\},
\end{array} \qquad \hbox{ for }
a \in \mathbb{R}, \, S \subset \mathbb{R}.
$$
Due to the nature of our the considerations, from now on we will also use some special notation:
\begin{namelist}{lll}
\item{$\bullet$} $\mathcal{T}$ denotes the family of compact,  convex, non-empty interior subsets of $\mathbb{R}^2$.

\item{$\bullet$} $\mathbb{T} \in \mathcal{T}$.

\item{$\bullet$} In this work, most of the constructions is based on the appropriate use of sequences of the partition
     of the set $\mathbb{T} \subset \mathbb{R}^2$.  To emphasize that the order of
     elements in a given partition is important to us, we will use the notation in the form of a vector -- $\bigl( t(j); 1\leqslant j \leqslant N\bigr)$. We assume that $\bigcup_{j=1}^N t(j) = \mathbb{T}$ but we do not require that the elements of a given
     partition are disjoint, in fact, as we shall see later, they cannot be disjoint.

\item{$\bullet$} Given a fixed sequence $( t(k) \colon 1 \leqslant k \leqslant N)$ of subsets of the set $\mathbb{T}$
     (elements of the set $2^{\mathbb{T}}$), then
     $$
     t \left(\leftsq k, j\rightsq \right) := \left\{ \begin{array}{lcl}
         \emptyset & \hbox{ for } & k >j; \\
         \bigcup_{j=k}^{j} t(j) & \hbox{ for } & k \leqslant j.
         \end{array} \right.
     $$

\item{$\bullet$} The partition $( t(k) \colon 1 \leqslant k \leqslant N)$ of subsets of the set $\mathbb{T}$ we
      call {\it sets population} if $t(\leftsq k,j\rightsq) \in \mathcal{T}$ for every $1 \leqslant k \leqslant j \leqslant N$.

\end{namelist}

\medskip

\noindent
We will use the simplified notation for subsets of $\mathbb{R}^2$, which is best explained with examples:
$$
\bigl( |x| \leqslant a, |y| \in [b(x), c(x)) \bigr) := \bigl\{ (x,y) \in \mathbb{R}^2 \colon |x| \leqslant a, |y| \in [b(x), c(x)) \bigr\},
$$
$$
\bigl( |x| \leqslant a \bigr) := \bigl\{ (x,y) \in \mathbb{R}^2 \colon |x| \leqslant a \bigr\}
$$
\medskip
if its clear that we use $\bigl( |x| \leqslant a \bigr)$ for a subset of $\R^2$.

\medskip

\section{Scheme of this paper}

The basic scheme of work is as follows:

For given compact and convex set $\TT \subset \R^2$ we reduced problem of construction of convex Peano curve filling $\TT$ to construction of some sequence of partitions of set $\TT$ (Chapter 2.).

Then we show that such sequence of partitions, called Fundamental Constructions, exists (Chapter 3.).

\medskip

The aim of Chapter 2. will be carried out by analysing
the Hausdorff characterization proof of compact sets (reduced to subsets of $\mathbb{R}^2$):

\medskip

\noindent
\begin{thm} {\bf (Hausdorff).}
\label{thm1}
The set $\mathbb{T} \subset \mathbb{R}^2$ is compact if and only if it is the image of the Cantor set
for some continuous function.
\end{thm}

The only difficulty in the proof is the construction of a continuous surjection $f \colon \Gamma \rightarrow
\mathbb{T}$ filling the given compact set $\mathbb{T}$. A complete proof
can be found in \cite{Sagan}, Section 6.6. We will use the following lemma based on the fact that any compact set in
$\mathbb{R}^2$ can be covered with a finite number of balls of arbitrarily small diameter.

\medskip

\noindent
\begin{lem}
\label{lem1.4.2}
Let $\mathbb{T} \subset \mathbb{R}^2$ be a compact set. There exist $N(j), n(j) \in \mathbb{N}$ for $j \in
\mathbb{N}$ satisfying $N(j+1) = N(j) n(j)$ and sequences of closed sets $\bigl(t(j; K), 1\leqslant K \leqslant N(j)\bigr)$ such that

{\bf 1)} $\bigcup_{1 \leq K \leq N(1)} t(1;K) = \TT$;

{\bf 2)} $t(j;K) = \bigcup_{(K-1)n(j) +1 \leq k \leq Kn(j)} t(j+1; k)$ for any $1 \leq K \leq N(j)$;

{\bf 3)} $\max_{1 \leqslant K \leqslant N(j)} d \bigl( t(j; K) \bigr) \rightarrow 0$ if $j \rightarrow \infty$.
\end{lem}

{\noindent}
{\bf Sketch of the proof of Hausdorff Theorem} we start ({\bf Step 1}) by selecting for the numbers $N(j), n(j)$, $j \in \mathbb{N}$, obtained in Lemma \ref{lem1.4.2} such a partition of the topological space $\Gamma$ into {open-closed sets }
$(I(j; K); 1 \leqslant K \leqslant N(j))$ for which the conditions are satisfied:
\begin{namelist}{ll}
\item{1')} $\bigcup_{1 \leq K \leq N(1)} I(1;K) = \Gamma$;
\item{2')} for any $1 \leq K \leq N(j)$ we have
 \begin{equation*}
 I(j;K) = \bigcup_{(K-1)n(j)+1 \leq k \leq Kn(j)} I(j+1;k);
 \end{equation*}
\item{3')} $\max_{1 \leq K \leq N(j)} d \bigl( I(j; K) \bigr) \rightarrow 0$  for 
$j \rightarrow \infty$.
\end{namelist}

\noindent
{\bf Step 2.} For any point $u \in \Gamma$ we find the sequence $(K(j), j \in \mathbb{N})$ such that
  $\{u\} = \bigcap_{j} I(j; K(j))$. Let
  $$
  \{ \mathbf{x}\} = \bigcap_{j\in \mathbb{N}} t\bigl(j; K(j) \bigr).
  $$
Finally we show that the mapping $\Gamma \ni u \mapsto f(u) = \mathbf{x} \in \mathbb{T}$ is uniquely determined
for $u \in \Gamma$ and it is a continuous surjection onto $\mathbb{T}$. \qed

\vspace{3mm}

A modification of this reasoning will lead us to the construction of a convex Peano curve on $[0,1]$
filling the set $\TT \in \T$, with $\operatorname{int}(f([a,b])) \neq \emptyset$ for $0 \leq a < b \leq 1$. Partitions $(I(j; K), 1 \leqslant K \leqslant N(j))$
of the Cantor set $\Gamma$ will be replaced with equal partitions of the interval $[0,1]$:
$$
I(j; K) := \Bigl[ \textstyle{\frac{K-1}{N(j)}, \frac{K}{N(j)}} \Bigr] \quad \hbox{for} \quad 1 \leqslant K \leqslant N(j).
$$
The partitions $(t(j; K); 1 \leqslant K \leqslant N(j))$ of $\TT$ will be constructed in the Chapter 3
(we will call it the Fundamental Construction). We will make sure that in addition to the conditions
$1)\div 3)$ in the appropriately changed equivalent of Lemma 3, the following condition is also satisfied: \\

\noindent
{\bf 4)} the sets $\bigcup_{K \leqslant J \leqslant L} t(j; J)$ are convex and with non-empty interior for all $j \in
 \mathbb{N}$, $1 \leqslant K \leqslant L \leqslant N(j)$. \\

Satisfying the condition 4) is the greatest difficulty in this construction as described in Chapter $3$.

\medskip

Vince and Wilson's hypothesis in \cite{VinceWilson} (Conjecture 4.3) comes to the statement that such a
construction is not possible. Therefore, the Fundamental Construction becomes one of the main results of this paper.

\newpage
\chapter*{2. Existence of the convex Peano curve - conditional proof}
\addcontentsline{toc}{chapter}{2. Existence of the convex Peano curve - conditional proof}
\stepcounter{chapter}

In this section, we will provide the {\it conditional} proof of the main result of the work. {\it Conditional},
because assuming that the Fundamental Construction is already done; in fact it will be done only in the next chapter.

\medskip

\begin{namelist}{ll}

\item{$\bullet$} We say that the function $g \colon [0,1] \mapsto \mathbb{R}^2$ is {\it non--resting} if $g((a,b))$ has non-empty interior for all $0\leqslant a < b \leqslant 1$.
\end{namelist}

\medskip
\noindent
The following theorem is not yet the Theorem \ref{thm0}. since it contains the Fundamental Construction as an assumption.

\medskip

\noindent
\begin{thm}
\label{thm2.0.1}
Let $\mathbb{T} \in \mathcal{T}$, $\bigl( t(j; K) \colon$ $1 \leqslant K \leqslant N(j) \bigr) \in \T^{N(j)}$ for $j, N(j) \in
\mathbb{N}$, be a sequences of set populations such that: \\

\begin{namelist}{lll}

\item[$1)$] $\bigcup_{1 \leq K \leq N(1)} t(1;K) = \TT$; \\

\item[$2)$] for every $j\in \mathbb{N}$ there exists $n(j) \in \N$ satisfying $N(j +1) = N(j) n(j)$ and $t(j; K) = t(j +1; [(K-1)n(j) +1, Kn(j)])$ for $1 \leq K \leq N(j)$; \\

\item[$3)$] $\max \Bigl\{d\bigl(t(j;K) \colon 1 \leq K \leqslant N(j) \Bigr\} \rightarrow 0$ for $j \rightarrow \infty$; \\
\end{namelist}

\noindent
Then there exists a non-resting convex Peano curve $f \colon [0,1] \mapsto \mathbb{T}$ filling the set $\mathbb{T}$.
\end{thm}

\medskip

\noindent
{\bf Proof of Theorem \ref{thm2.0.1}.} {\bf Step I.} For the sequence of sets populations $(t(j;K); 1 \leq K \leq N(j))$ there exists a continuous function $f \colon [ 0,1] \mapsto \mathbb{T}$ satisfying for all $j \in \mathbb{N}$
\begin{equation}
\label{eq2.0.1}
f \left( \left[ \textstyle{\frac{K-1}{N(j)}, \frac{K}{N(j)}} \right] \right) = t(j; K) \qquad
1 \leqslant K \leqslant N(j).
\end{equation}

{\bf Proof of step I.} Let $I(j; K) = [\frac{K-1}{N(j)}, \frac{K}{N(j)}]$ for $1 \leq K \leq N(j), j \in \mathbb{N}$. Condition 2) yields the following relations between the coverings $(I(j;K); 1 \leq K \leq N(1))$ and $(t(j;K); 1 \leq K \leq N(j))$. For $j_0 \in \mathbb{N}, 1 \leq K_0 \leq N(j_0)$ and for any $j > j_0$, we have
\begin{equation*}
t(j_0; K_0) = \bigcup_{1 \leq K \leq N(j), I(j;K) \subset I(j_0; K_0)} t(j; K).
\end{equation*}

So according to condition 3), if a set $S \subset \mathbb{T}$ satisfies the implication,
\begin{equation}
\label{eq*}
I(j;K) \subset I(j_0;K_0) \Longrightarrow S \cap t(j; K) \neq \emptyset,
\end{equation}
then the set $S$ is dense in $t(j_0;K_0)$.

For $u \in [0,1]$, let's arbitrarily choose numbers $K(j)$ for $j \in \mathbb{N}$ satisfying \\
$u \in I(j; K(j)), I(j +1; K(j +1)) \subset I(j; K(j))$ for $j \in \mathbb{N}$. Then there exists a point $f(u) \in \mathbb{T}$ defined by the condition $\{ f(u) \} = \bigcap_{j \in \mathbb{N}} t(j; K(j))$, which means $f(u) \in t(j; K(j))$ for $j \in \mathbb{N}$. If numbers $\bar{K}(l)$ for $j \in \mathbb{N}$ are analogously defined for $\bar{u} \in [0,1]$, then for $j_0 \in \mathbb{N}$ the condition $|u-\bar{u}| < \frac{1}{N(j_0)}$ implies $I(j_0; K(j_0)) \cap I(j_0; \bar{K}(j_0)) \neq \emptyset$, which leads to $|K(j_0) - \bar{K}(j_0)| \leq 1$, hence $t(j_0; K(j_0)) \cup t(j_0; \bar{K}(j_0)) \in \T$, implying $t(j_0; K(j_0)) \cap t(j_0; \bar{K}(j_0)) \neq \emptyset$, and finally $||f(u) - f(\bar{u})|| \leq 2 \max \{ d(t(j_0;K): 1 \leq K \leq N(j_0)) \}$. According to condition 3), we have proven the continuity of the obtained mapping $f: [0,1] \rightarrow \mathbb{T}$.

If $u \in \text{ int}(I(j_1; K_1))$, then the corresponding sequence $(K(j); j \in \mathbb{N})$ satisfies $K(j_1) = K_1$. This proves the inclusion $f(\text{ int}(I(j_1;K_1))) \subset t(j_1; K_1)$ for $1 \leq K_1 \leq N(j_1), j_1 \in \mathbb{N}$. The set $S = f(\text{int}(I(j_0; K_0)))$ thus satisfies $S \subset t(j_0; K_0)$ and condition (\ref{eq*}). As observed, $S$ is then dense in $I(j_0;K_0)$, proving (\ref{eq2.0.1}) for the continuous function $f$.

\medskip
{\bf Step II.} The function $f \colon [0,1] \mapsto \mathbb{T}$ obtained in Step I is a non-resting
convex Peano curve filling $\TT$.

\medskip
\noindent
{\bf Proof of Step II.} By property $(2.0.1)$ for $j = 1$ we see that $f([0,1]) = \mathbb{T}$. By property $(2.0.1)$
and since $t(j; K) \in \mathcal{T}$, where $\mathcal{T}$ is a family of sets with non-empty interiors we get that $f$ is a
non-resting function. It remains to prove that for every $0 \leq a < b \leq 1$ the set $f([a,b])$ is convex.

We choose sequences $(K_{j}), (\bar{K}_{j})$ satisfying $a \in I(j; K_{j}), b \in I(j; \bar{K}_{j})$ for $j \in \N$. Then we have 
$$
 [a,b] \subset \left[ \textstyle{ \frac{K_{j} -1}{N(j)}, \frac{\bar{K}_{j}}{N(j)}} \right] =
\bigcup_{J= K_{j}}^{\bar{K}_{j}} I(j; J).
$$
Consequently
$$
f\bigl( [a,b] \bigr) \subset f \Bigl( \left[ \textstyle{ \frac{K_{j} -1}{N(j)}, \frac{\bar{K}_{j}}{N(j)}} \right] \Bigr) = \bigcup_{J= K_{j}}^{\bar{K}_{j}} f \Bigl( I(j; J) \Bigr) =
\bigcup_{J= K_{j}}^{\bar{K}_{j}} t(j;J)
,
$$
by (2.0.1), and the last set is convex by the definition of sets population. Thus we have that
$$
s := \bigcap_{j \in \mathbb{N}} f \Bigl( \left[ \textstyle{ \frac{K_{j} -1}{N(j)}, \frac{\bar{K}_{j}}{N(j)}} \right] \Bigr)
$$
is also a convex set as an intersection of the family of convex sets.
On the other hand $s = \bigcap_{j \in \N} (f([a,b]) \cup f(\Delta_{j}))$ for
\begin{equation*}
\Delta_{j} := \left[ \textstyle{ \frac{K_{j} -1}{N(j)}, \frac{\bar{K}_{j}}{N(j)}} \right] \setminus [a,b]
\end{equation*}
and by continuity of function $f$, $\bigcap_{j \in \N} f(\Delta_{j}) \subset \{ f(a), f(b) \}$. The obtained equality $f([a,b]) = s$ completes the proof of Step II and the Theorem \ref{thm2.0.1}. \qed

\newpage
\chapter*{3. Fundamental Construction}
\addcontentsline{toc}{chapter}{3. Fundamental Construction}
\stepcounter{chapter}

\section{The idea of souls populations}

Recall that $\T$ is a family of compact, convex, non-empty interior subsets of $\R^2$.

\begin{df}
\label{df3.1.1}
The pair of sets $(t,t_1)$ is a \underline{soul} if $t_1 \subset t; t,t \setminus t_1 \in \T$. Sets $t$ and $t_1$ are called respectively \underline{base} and \underline{disturbance} of the soul.
\end{df}

\begin{namelist}{lll}
\item{$\bullet$} $\tau_0$ will be used for a space of all souls.
\end{namelist}

For sequences of sets $\ttt = (t(k); 1 \leq k \leq m)$, $\ttt_1 = (t_1(k); 1 \leq k \leq m)$, $\sss = (s(k); 1 \leq k \leq m)$ and any set $\mathfrak{m} \subset \Z$, we use the notations

\begin{namelist}{lll}
\item{$\bullet$} $\Delta_1 \ttt = \sss$ if $s(k) = t(k) \setminus t(k-1)$ for $ 2 \leq k \leq m$, $s(1) = \emptyset$;
\item{$\bullet$} $\Delta_{-1} \ttt = \sss$ if $s(k) = t(k) \setminus t(k+1)$ for $ 1 \leq k \leq m-1$, $s(m) = \emptyset$;
\item{$\bullet$} $\mathbf{J} \ttt = t(1), \mathbf{S} \ttt = t(m)$ (the first and the least entry of a sequence $\ttt$);
\item{$\bullet$} $\ttt_1 |_{\mathfrak{m}} = \sss$ if $s(k) = \left\{ \begin{array}{lcl}
      t_1(k) & \hbox{for} & k \in \mathfrak{m} \\
      \emptyset & \hbox{for} & k \notin \mathfrak{m} \end{array} \right.$ for $k \in \{1,\dots, m \}$;
\end{namelist}
useful and frequenly used will be shorthand
\begin{namelist}{lll}
\item{$\bullet$} $\ttt_1 |_1 = \ttt_1 |_{2 \N -1}, \ttt_1 |_{-1} = \ttt_1 |_{2 \N}$;   
\item{$\bullet$} $<\ttt,\ttt_1> = ((t(k), t_1(k)); 1 \leq k \leq m)$ if $(t(k), t_1(k))$ are souls for $1 \leq k \leq m$.
\end{namelist}

The sequences $\ttt, \mathbf{t}_1$ and the number $m$ will be respectively referred to as the \textit{base}, the \textit{disturbances} and the \textit{length} for the sequence of souls $<\ttt, \ttt_1>= ((t(k), t_1(k)); 1 \leq k \leq m)$.

The notation $((t, t_1)(k); 1 \leq k \leq m)$ will always denote an abbreviation for $((t(k), t_1(k)); 1 \leq k \leq m)$.

\begin{df}
\label{df3.1.2}
The sequence of souls $<\ttt, \ttt_1>$ of length $M$ is \underline{regular} if $M \in 2 \N$ and $\ttt_1 |_1 = \Delta_1 \ttt, \ttt_1 |_{-1} = \Delta_{-1} \ttt$.
\end{df}

Then immediately
\begin{lem}
\label{lem3.1.3}
If $((t,t_1)(K); 1 \leq K \leq M)$ is a regular sequence, then
\begin{enumerate}[a)]
\item $t(2J-1) = t(2J)$ for $1 \leq J \leq \frac{M}{2}$ and $t(2J) \setminus t_1(2J) = t(2J+1) \setminus t_1(2J+1)$ for $1 \leq J \leq \frac{M}{2} - 1$;
\item $t(K) \cap t(K+1) \in \mathcal{T}$ for $1 \leq K \leq M-1$;
\item if sequences $<\ttt, \ttt_1>, <\bar{\ttt}, \bar{\ttt_1}>$ are regular and $\mathbf{S} \ttt = \mathbf{J} \bar{\ttt}$, then the concatenation $<\ttt, \ttt_1> \circ <\bar{\ttt}, \bar{\ttt}_1> = <\ttt \circ \bar{\ttt}, \ttt_1 \circ \bar{\ttt}_1>$ is also a regular sequence.
\end{enumerate}
\end{lem}

\begin{namelist}{lll}
\item{$\bullet$} We say that sets $r, s \in \mathbb{R}^2$ are \textit{separated} if inf$\{ || \mathbf{x} - \mathbf{y} ||: \mathbf{x} \in s, \mathbf{y} \in r \} > 0$ or $r = \emptyset$ or $s = \emptyset$.

\item{$\bullet$} For sequences of sets $\mathbf{q} = (q(k); 1 \leq k\leq m), \mathbf{r} = (r(k); 1 \leq k \leq m), \mathbf{s} = (s(k); 1 \leq k \leq m)$, we define the relation $\mathbf{q} \prec_{\mathbf{s}} \mathbf{r}$ if for any indices $k,l \in \{1,\dots,m \}, k <l$, the following alternatives hold: $q(k) = r(l)$ or $q(k), r(l)$ are separated or $q(k) \in s([k+1,l-1])$ or $r(l) \in s([k+1, l-1])$.
\end{namelist}

\begin{df}
\label{df3.1.4}
A sequence of souls $<\ttt, \ttt_1>$ is \underline{consistent} if $\ttt_1 |_{-1} \prec_{\ttt} \ttt_1 |_1$.
\end{df}

\begin{df}
\label{df3.1.5}
We will call a sequence $<\ttt, \ttt_1>$ a \underline{population of souls} if it is both regular and consistent.
\end{df}

The main result of this section is:

\begin{thm}
\label{thm3.1.6}
If $<\ttt, \ttt_1>$ is a population of souls, then the base $\ttt = (t(K); 1 \leq K \leq M)$ forms a population of sets (meaning $t([K,L]) \in \mathcal{T}$ for $1 \leq K \leq L \leq M$).
\end{thm}

\begin{lem}
\label{lem3.1.7}
If $<\ttt, \ttt_1> = ((t,t_1)(K); 1 \leq K \leq M)$ is a population of souls, then:
\begin{enumerate}[1)]
\item $t(K) \cap t(K+1) \neq \emptyset$ for $1 \leq K \leq M-1$;
\item $t(K) \setminus t(K+1), t(K+1) \setminus t(K)$ are separated for $1 \leq K \leq M-1$;
\item for $K,L \in \{1, \dots, M\}, L-K \geq 2, p := t(K) \setminus t([K+1, L-1]), q := t(L) \setminus t([K+1, L-1])$, we have $p=q$ or $p,q$ are separated.
\end{enumerate}
\end{lem}

\textbf{Proof.} Condition 1) follows from Lemma \ref{lem3.1.3},b).

Let $(s(K); 1 \leq K \leq M) := \Delta_1 \ttt, (r(K); 1 \leq K \leq M) := \Delta_{-1} \ttt$. Since $r(K) \cap s(K+1) = \emptyset$, the alternative in the definition of the relation $\prec_{\ttt}$ implies that $r(K),s(K+1)$ are separated, hence we have 2).

Let's prove 3). For $1 \leq K < L \leq M$, denote $J=L-K$. Inductively, for $J = 2,3, \dots, M-1$, we see easily that $p \in \{ r(K), \emptyset \}, q \in \{ s(L), \emptyset \}$. Condition 3) again results from the alternative in the definition of $\prec_{\ttt}$. \qed

\medskip

The consequences of conditions 1), 2), 3) are obtained based on the following interesting fact.

\begin{lem}
\label{lem3.1.8}
If for sets $t,\bar{t} \in \mathcal{T}$, the differences $t \setminus \bar{t}, \bar{t} \setminus t$ are separated and $t \cap \bar{t} \neq \emptyset$, then $t \cup \bar{t} \in \mathcal{T}$.
\end{lem}

{\bf Proof.} We need to prove that $[a,b] \subset t\cup \overline{t}$ for every $a \in t$, $b \in \overline{t}, a \neq b$. To simplify, let's assume that $a = (0,0), b = (0,1)$. If there exists a point $(0,y)
\in t \cap \overline{t}$, for some $y \in (0,1)$ then the result is immediate. In the opposite case we have
$$
\forall_{y \in (0,1)}: (0,y) \not\in t \cap \bar{t} \quad \hbox{ and } \quad \exists_{x_1 \neq 0} \exists_{y_1 \in \mathbb{R}}: (x_1,y_1) \in t \cap \overline{t}.
$$
It is enough to consider the case $x_1 =1$. Continuing the proof by contradiction, suppose that there exists
$y_0  \in (0,1)$ such that $(0,y_0) \not\in t \cup \overline{t}$. We define two functions, for $x \in [0,1]$,
\begin{eqnarray*}
g(x) & = & \max \left\{ y \colon (x,y) \in t \right\}, \\
\overline{g}(x) & = & \min \left\{ y \colon (x,y) \in \overline{t} \right\}.
\end{eqnarray*}
These functions are continuous because $t, \overline{t} \in \mathcal{T}$, and they have the
following properties
$$
g(0) < y_0 < \overline{g}(0), \qquad g(1) \geqslant y_1 \geqslant \overline{g}(1).
$$
Thus there exists $x_1 \in (0,1]$ the intersection point of $g$ and $\overline{g}$ and we have
$$
y_2 := g(x_1) = \overline{g}(x_1), \qquad g(x) < \overline{g}(x) \hbox{ for all } x \in [0, x_1).
$$
Then $(x_1, y_2) \in t \cap \overline{t}$ and $(x_1,y_2)$ is the limit of points $(x, g(x)) \in t \setminus
\overline{t}$ and points $(x, \overline{g}(x)) \in \overline{t} \setminus t$. This gives a contradiction
because the sets $t \setminus \overline{t}$ and $\overline{t} \setminus t$ are separated by assumption. \qed

\begin{lem}
\label{lem3.1.9}
A sequence of sets $(t(K); 1 \leq K \leq M)$ satisfying conditions 1), 2), 3) in Lemma \ref{lem3.1.7} forms a population of sets.
\end{lem}

{\bf Proof.} We will show that for the sequence $\bigl( (t(K) \colon 1 \leqslant
K \leqslant N\bigr)$ the set $t(\leftsq K,L \rightsq )$ is convex for any $1 \leqslant K \leqslant L \leqslant N$.
We will use mathematical induction for $j = L-K  = 0, 1, \dots, N-1$. \\

--- For $j=0$ we have that $t([K,L]) = t(K)$, which is assumed convex. \\

--- For $j=1$, we use the Lemma \ref{lem3.1.8} for $t = t(L)$ and $\overline{t} = t(L-1)= t(K)$. According to conditions 1) and 2) the sets
$t \setminus \bar{t}$ and $\overline{t} \setminus t$ are separated, $t \cap \overline{t} \neq \emptyset$ so the set $t(K) \cup t(L)$
is convex. \\

--- Assume that $j \geqslant 2$ and $t(\leftsq K,L \rightsq )$ is a convex set for all $L-K \leqslant j-1$. For fixed
$K,L$ satisfying the condition $L-K = j$ we have that the sets $t(\leftsq K+1, L\rightsq ), t(\leftsq K, L-1 \rightsq )$
are convex. If $t(K) \setminus t(\leftsq K+1, L-1 \rightsq ), t(L) \setminus t(\leftsq K+1, L-1 \rightsq )$ are separated
then using Lemma \ref{lem3.1.8} for $t = t(\leftsq K, L-1 \rightsq )$ and $\overline{t} = t(\leftsq K+1, L
\rightsq )$ we get convexity $t(\leftsq K, L \rightsq )$. Otherwise, in view of alternatives in condition $3)$ we have $t(K) \setminus t(\leftsq K+1, L-1 \rightsq ) = t(L) \setminus t(\leftsq K+1, L-1 \rightsq )$. But then $t( \leftsq K,L\rightsq ) = t( \leftsq K,L-1 \rightsq )$, which by inductive assumption
is also convex.\qed

\hfill \break
\textbf{Proof of Theorem \ref{thm3.1.6}.} We apply Lemma \ref{lem3.1.7} and Lemma \ref{lem3.1.9}.

\section{Offspring of souls - axiomatic approach}

Let us fix families of souls $\tau,\bar{\tau}\subset\tau_0$.
A somewhat complex definition of the offspring of souls being a function
$\tau\ni(t,t_1)\mapsto \left\langle\bs{t}'(t,t_1),\bs{t}'_1(t,t_1)\right\rangle\in\bar{\tau}^{m'}$
will be preceded with an introduction of some additional notation. For sequences 
$\mathbf{r}=\left(r(i); 1\leq i\leq n\right)$ and
$\mathbf{s}=\left(s(i); 1\leq i\leq n\right)$ of subsets of $\mathbb{R}^2$, and a set $\bar{r}\subset\mathbb{R}^2$ we shall write
\begin{namelist}{lll}
	\item{$\bullet$} $\bs{r}\perp\bs{s}$, if for all $1\leq i\leq n$ we have $r(i)\cap s(i)=\emptyset$; we will say that \(\bs{r}\) and \(\bs{s}\) are \textit{orthogonal};
	\item{$\bullet$} $\bs{r}=\bar{r}\otimes\bs{s}$, if $r(i)=\bar{r}\cap s(i)$ for all $1\leq i\leq n$;
	\item{$\bullet$} $\bigcup\bs{s}:=\bigcup_{1\leq i\leq n}s(i)$;
	\item{$\bullet$} $\bs{s}\ssucc\bs{p}$, for a sequence \(\bs{p}=\left(p(j);1\leq j\leq k\right)\), if for some indices $1\leq i(1)<\ldots<i(k)\leq n$ we have $s(i(j))=p(j)$, for $j\in\{1,\ldots,k\}$, and $s(i)=\emptyset$, for $i\notin\{i(1),\ldots,i(k)\}$;
	\item{$\bullet$} $\bs{r}\equiv\bs{s}$, if sequences $\bs{r}$ and $\bs{s}$ satisfy \(\bs{r}\ssucc\bs{p}\), \(\bs{s}\ssucc\bs{p}\), for some \(\bs{p}\);
	\item{$\bullet$} $\bs{r}\ll \bs{s}$, if for all indices $k,l\in\{1,\ldots,n\}$, $k<l$, the sets $r(k)$ and $s(l)$ are separated.
\end{namelist}

The notation \(\bs{J}\bs{s}\), \(\bs{S}\bs{s}\), \(\bs{s}|_{\mathfrak{m}}\), \(\bs{s}|_{1}\), and \(\bs{s}|_{-1}\) will be used in accordance with Section 3.1.

Let us fix a number $m'\in2\mathbb{N}$ and thus the index set $\mathfrak{m}'=\{1,\ldots,m'\}$. With $\mathfrak{M}'$ we shall denote the space of all sequences $\bs{s}=\left(s(k);k\in\mathfrak{m}'\right)$ of subsets of $\R^2$.

The first step to define the offspring of souls consists in the following:
\begin{df}
\label{df3.2.1}
	Function $\tau\ni(t,t_1)\mapsto\bs{s}(t,t_1)\in\mathfrak{M}'$ will be called a \emph{dust} if for all $(t,t_1),(\bar{t},\bar{t}_1)\in\tau$ we have
	\begin{itemize}
	\item[1)] $\left(t_1\otimes\bs{s}(t,t_1)\right) \perp \left((t\setminus t_1)\otimes\bs{s}(t,t_1)\right)$;
	\item[2)] $t_1\otimes \bs{s}(t,t_1) \ll (t\setminus t_1)\otimes\bs{s}(t,t_1)$;
	\item[3)] ${t_1=\bar{t}_1} \implies {t_1\otimes\bs{s}(t,t_1)\equiv{ \bar{t}_1\otimes\bs{s}(\bar{t},\bar{t}_1)}}$;
	\item[4)] ${t\setminus t_1 = \bar{t}\setminus \bar{t}_1}
	\implies
	{(t\setminus t_1)\otimes\bs{s}(t,t_1)} \equiv {(\bar{t}\setminus\bar{t}_1)\otimes \bs{s}(\bar{t},\bar{t})}$.
	\end{itemize}
	Function  $\tau\ni(t,t_1)\mapsto\bs{r}(t,t_1)\in\mathfrak{M}'$ will be called an \emph{anti-dust} if for all $(t,t_1),(\bar{t},\bar{t}_1)\in\tau$ we have
	\begin{itemize}
		\item[5)] $t=\bar{t} \implies \bs{r}(t,t_1)\equiv\bs{r}(\bar{t},\bar{t}_1)$.
	\end{itemize}
	Function  $\tau\ni(t,t_1)\mapsto\bs{q}(t,t_1)\in\mathfrak{M}'$ will be called a \emph{filling} if for all $(t,t_1),(\bar{t},\bar{t}_1)\in\tau$ we have
	\begin{itemize}
		\item[6)] \(t=\bigcup\bs{q}(t,t_1)\),
		\item[7)] ${t\setminus t_1 = \bar{t}\setminus \bar{t}_1}
		\implies \bs{J}\bs{q}(t,t_1)=\bs{J}\bs{q}(\bar{t},\bar{t}_1)$,
		\item[8)] $t=\bar{t} \implies \bs{S}\bs{q}(t,t_1)=\bs{Sq}(\bar{t},\bar{t}_1)$.
	\end{itemize}
\end{df}

Now we can define the offspring of souls.
\begin{df}
\label{df3.2.2}
	Function $\tau\ni(t,t_1)\mapsto\left\langle \bs{t}'(t,t_1),\bs{t}'_1(t,t_1)\right\rangle\in\bar{\tau}^{m'}$ will be called \emph{offspring of souls} (or \emph{offspring of length $m'$ on $\tau$ in $\bar{\tau}$}) if
	\begin{itemize}
		\item[1)] for $(t,t_1)\in\tau$ the function value $\left\langle \bs{t}'(t,t_1),\bs{t}'_1(t,t_1)\right\rangle$ is a population of souls,
		\item[2)] the mapping \(\tau\ni(t,t_1)\mapsto\bs{t}'_1(t,t_1)|_1\) is a dust,
		\item[3)] the mapping \(\tau\ni(t,t_1)\mapsto\bs{t}'_1(t,t_1)|_{-1}\) is an anti-dust,
		\item[4)] the mapping \(\tau\ni(t,t_1)\mapsto\bs{t}'(t,t_1)\) is a filling.
	\end{itemize}
\end{df}

The value of the function $\langle \mathbf{t}'(t,t_1), \mathbf{t}_1'(t,t_1) \rangle$ will be called the \textit{offspring of the parent} $(t,t_1)$.

For any sequence $\mathbf{a}=(a(k); 1\leqslant k\leqslant m')$, we define the inversion operation $I\mathbf{a}=(a(m'+1-k); 1\leqslant k \leqslant m')$ and recursively $I^0\mathbf{a}=\mathbf{a}$, $I^K\mathbf{a}=I(I^{K-1}\mathbf{a})$ for $K\in \N$. For any sequence $(\mathbf{a}(K); 1\leqslant K\leqslant M)$ composed of sequences $\mathbf{a}(K)=(a(K,k); k\in\mathfrak{m}')$, we call the \textit{anti-order} the concatenation:

\begin{equation}
\label{eq3.2.1}
\mathcal{A}(\mathbf{a}(K); 1\leqslant K\leqslant M)=I^0\mathbf{a}(1)\circ I^1\mathbf{a}(2)\circ\cdots\circ I^{M-1}(\mathbf{a}(M)).
\end{equation}

For a sequence of souls $\langle \mathbf{t}, \mathbf{t}_1 \rangle \in (\tau_0)^{m'}$, we have $I\langle \mathbf{t}, \mathbf{t}_1 \rangle=\langle I\mathbf{t}, I\mathbf{t}_1 \rangle$, and for a sequence composed of sequences of souls $(\langle \mathbf{t}(K), \mathbf{t}_1(K) \rangle; 1\leqslant K\leqslant)$ in the space $(\tau_0)^{m'}$, we have for its anti-order $$\mathcal{A}(\langle \mathbf{t}(K), \mathbf{t}_1(K) \rangle; 1\leqslant K\leqslant M)= \langle \mathcal{A}(\mathbf{t}(K); 1\leqslant K\leqslant M); \mathcal{A} (\mathbf{t}_1(K); 1\leqslant K\leqslant M) \rangle.$$

Now, we formulate the main result of the section.

\begin{thm}
\label{thm3.2.3}
Let $\tau\ni(t,t_1)\mapsto\langle \mathbf{t}'(t,t_1), \mathbf{t}'_1(t,t_1) \rangle \in \bar{\tau}^{m'}$ be the offspring of souls, and let $((t(K), t_1(K)); 1\leqslant K\leqslant M) \in \tau^M$ be the population of souls. Then the anti-order

\begin{eqnarray}
\label{eq3.2.2}
&& \langle \mathbf{t}^A, \mathbf{t}^A_1 \rangle =\\
\nonumber &&=\langle \mathcal{A}(\mathbf{t}'(t(K), t_1(K)); 1\leqslant K\leqslant M), \mathcal{A}(\mathbf{t}_1'(t(K), t_1(K)); 1\leqslant K\leqslant M) \rangle
\end{eqnarray}
is a new population of souls, and for $(t^A(\aae); 1\leqslant \aae \leqslant Mm')=\mathbf{t}^A$, we have

\begin{equation}
\label{eq3.2.3}
t(K)=t^A([(K-1)m'+1, Km']) \ \text{for} \ 1\leqslant K\leqslant M.
\end{equation}
\end{thm}

Easily, we obtain

\begin{lem}
\label{lem3.2.4}
The sequence of souls defined by the anti-order (\ref{eq3.2.2}) is regular and satisfies (\ref{eq3.2.3}).
\end{lem}

\medskip
\noindent
\textbf{Proof.}
Notice that for any regular sequence of souls $\langle \mathbf{t}, \mathbf{t}_1 \rangle$, inversion yields a new regular sequence of souls of the form $\langle I\mathbf{t}, I\mathbf{t}_1 \rangle$.

To achieve the regularity of sequence (\ref{eq3.2.2}), according to Corollary 3.1.3.c), it suffices to show for $(t'(K,k); k\in \mathfrak{m}')=\mathbf{t}'(t(K), t_1(K))$, $(\bar{t}'(K,k); k\in\mathfrak{m}')=I\mathbf{t}'(t(K), t_1(K))$, and any $1\leqslant J\leqslant M/2$, that $t'(2J-1, m')=\bar{t}'(2J, 1)$ and $\bar{t}'(2J, m')=t'(2J+1, 1)$ (when $J\leqslant M/2-1$). Thus, we need to show $t'(2J-1, m')=t'(2J, m')$ and $t'(2J, 1)=t'(2J+1, 1)$ (when $J\leqslant M/2-1$). Since the function $\tau \ni (t,t_1)\to \mathbf{t}'(t,t_1)$ is a filling, according to conditions 7) and 8) in Definition 3.2.1, it suffices to use the equality $t(2J-1)=t(2J)$ and $t(2J)\smallsetminus t_1(2J)=t(2J+1)\smallsetminus t_1(2J+1)$ (when $J\leqslant M/2-1$) (compare with lemma 3.1.3.a)).

Equality (\ref{eq3.2.3}) follows from condition 6) in Definition 3.2.1. \qed

\medskip

The proof of the consistency of the anti-order (\ref{eq3.2.2}) is cumbersome. The population of souls $(t(K), t_1(K); 1\leqslant K\leqslant M)$ will remain fixed. We denote the offspring of the parent $(t(K), t_1(K))$ shortly as $\langle \mathbf{t}'(K), \mathbf{t}'_1(K) \rangle = \langle \mathbf{t}'(t(K),t_1(K)), \mathbf{t}'_1(t(K), t_1(K)) \rangle$, and also denote $(t^A(\aae); 1\leqslant\aae\leqslant Mm')=\mathbf{t}^A$, $(t^A_1(\aae); 1\leqslant\aae\leqslant Mm')=\mathbf{t}^A_1$ for the anti-order (\ref{eq3.2.2}).

\begin{namelist}{lll}
\item{$\bullet$} We set the index sets $\mathfrak{m}(K)=(K-1)m'+\mathfrak{m}'$ for $1\leqslant K\leqslant M$ (thus $\mathfrak{m}(1)=\mathfrak{m}'$); $\mathfrak{m}^A=\mathfrak{m}(1)\cup \cdots \cup \mathfrak{m}(M)=\{1, \cdots, Mm'\}$. Then we denote the spaces of sequences: $\mathbf{s}\in\mathfrak{M}^A$ if $\mathbf{s}=(s(\aae); \aae\in\mathfrak{m}^A)$, $s(\aae)\subset \mathbb{R}^2$ for $\aae\in\mathfrak{m}^A$; $\mathfrak{M}(K)=\{\mathbf{s}\in\mathfrak{M}^A : \mathbf{s}|_{\mathfrak{m}(K)}=\mathbf{s} \}$ for $1\leqslant K\leqslant M$.
\end{namelist}

We will frequently use the following operations on sequences of sets. For $\mathbf{p}=(p(k); k\in\mathfrak{m}')$, $\mathbf{s}=(s(\aae);\aae\in\mathfrak{m}^A)$, and $K\in\{ 1,\cdots,M \}$, we define
\begin{namelist}{lll}
\item{$\bullet$} $\mathbf{p}\!\!\nnearrow_1^K=\mathbf{s} \text{ if } \mathbf{s}\in\mathfrak{M}(K) \text{ and } s(\aae)=p(\aae-(K-1)m') \text{ for } \aae\in\mathfrak{m}(K)$,
\item{$\bullet$} $\mathbf{p}\!\!\nnearrow_{-1}^K=(I\mathbf{p})\!\!\nnearrow_1^K$.
\end{namelist}

For $r\subset\mathbb{R}^2$ and $\varepsilon\in\{-1,1\}$, we have $(r\otimes\mathbf{p})\!\!\nnearrow_{\varepsilon}^K=r\otimes(\mathbf{p}\!\!\nnearrow_{\varepsilon}^K)$, which justifies the notation $r\otimes\mathbf{p}\!\!\nnearrow_{\varepsilon}^K$.

For $\mathbf{r}=(r(\aae);\aae\in\mathfrak{m}^A)$, $\mathbf{s}=(s(\aae);\aae\in\mathfrak{m}^A)$, where $\mathbf{r}, \mathbf{s}\in\mathfrak{M}^A$, we denote 
\begin{namelist}{lll}
\item{$\bullet$} $\mathbf{r}\oplus\mathbf{s}=(r(\aae)\cup s(\aae); \aae\in\mathfrak{m}^A)$. As usual, $\bigoplus\limits_{1\leqslant i\leqslant k}\mathbf{s}_k=\mathbf{s}_1\oplus\cdots\oplus\mathbf{s}_k$ for $\mathbf{s}_1,\cdots,\mathbf{s}_k\in\mathfrak{M}^A$.

\item{$\bullet$} $\mathbf{r}\stackrel{-1}{\equiv}\mathbf{s}$ if, for some decreasing bijection $\{\aae\in\mathfrak{m}^A:r(\aae)\neq\emptyset\}\leftrightarrow\{\lambda\in\mathfrak{m}^A:s(\lambda)\neq\emptyset\}$, we have the implication $\aae\leftrightarrow\lambda\Longrightarrow r(\aae)=s(\lambda)$.
\end{namelist}

Notice that in the adopted notation we have $\mathbf{t}^A, \mathbf{t}_1^A\in\mathfrak{M}^A$, and

\begin{equation}
\label{eq3.2.6}
\mathbf{t}^A=\bigoplus\limits_{1\leqslant J\leqslant M/2}(\mathbf{t}'(2J-1)\!\!\nnearrow_1^{2J-1}\oplus\mathbf{t}'(2J)\!\!\nnearrow_{-1}^{2J}),
\end{equation}

\begin{equation}\label{eq:piatka}
\mathbf{t}_1^A=\bigoplus\limits_{1\leqslant J\leqslant M/2}(\mathbf{t}_1'(2J-1)\!\nnearrow_1^{2J-1}\oplus\mathbf{t}_1'(2J)\!\nnearrow_{-1}^{2J}).
\end{equation}

The relation $\stackrel{-1}{\equiv}$ is symmetric and we have the following conditions (of course $\mathbf{s}|_{\eta}\!\!\nnearrow_{\varepsilon}^K=(\mathbf{s}|_{\eta})\!\!\nnearrow_{\varepsilon}^K$ for $\varepsilon, \eta\in\{1,-1\}$). For sequences $\rrr = (r(\aae); \aae \in \mathfrak{m}^A), \sss = (s(\lambda); \lambda \in \mathfrak{m}^A)$ we denote
\begin{namelist}{lll}
\item{$\bullet$} $\rrr < \sss \text{ if for } \aae, \lambda \in \mathfrak{m}^A, r(\aae) \neq \emptyset, s(\lambda) \neq \emptyset \text{ we have } \aae < \lambda$
\end{namelist}

\begin{cor}
\label{cor3.2.5}
For sequences $\mathbf{p},\mathbf{q}\in\mathfrak{M}'$, numbers $K,L\in\{1,\cdots,M/2\}$ and set $s\subset\mathbb{R}^2$ we have
\begin{enumerate}
\item[A.] If $\mathbf{p}\equiv\mathbf{q}$, then $\mathbf{p}\!\!\nnearrow_1^K\stackrel{-1}{\equiv}\mathbf{q}\!\!\nnearrow_{-1}^L$.
\item[B.] If $\mathbf{p}<<\mathbf{q}$, then $\mathbf{p}\!\!\nnearrow_1^K<<\mathbf{q}\!\nnearrow_1^K$ and $\mathbf{q}\!\nnearrow_{-1}^K<<\mathbf{p}\!\!\nnearrow_{-1}^K$.
\item[C.] For $\varepsilon, \eta\in\{1,-1\}$ we have $(s\otimes\mathbf{p})|_{\eta}\!\!\nnearrow_{\varepsilon}^K=s\otimes(\mathbf{p}|_{\eta}\!\!\nnearrow_{\varepsilon}^K)$ (hence we will write $s\otimes\mathbf{p}|_{\eta}\!\!\nnearrow_{\varepsilon}^K$).
\item[D.] For $\varepsilon, \eta\in\{1,-1\}$ we have $\mathbf{p}|_{\eta}\!\!\nnearrow_{\varepsilon}^K=\mathbf{p}\!\!\nnearrow_{\varepsilon}^K\!\!|_{\eta\varepsilon}$.
\item[E.] We have equivalences $\mathbf{p}<\mathbf{q}\Longleftrightarrow\mathbf{p}\!\nnearrow_1^K<\mathbf{q}\!\nnearrow_1^K\Longleftrightarrow\mathbf{q}\!\nnearrow_{-1}^K<\mathbf{p}\!\nnearrow_{-1}^K$ and, for $\varepsilon,\eta\in\{1,-1\}$, the implication $K<L\Longrightarrow\mathbf{p}\!\nnearrow_{\varepsilon}^K<\mathbf{q}\!\nnearrow_{\eta}^L$.
\end{enumerate}
\end{cor}

\medskip
\noindent
\textbf{Proof.}
For $K=L=1$, for any sequence $\mathbf{p}\in\mathfrak{M}'$ and for $(\bar{s}(\aae); \aae\in\mathfrak{m}^A)=\mathbf{p}\!\nnearrow_1^1$, $(\bar{\bar{s}}(\lambda)$; $\lambda\in\mathfrak{m}^A)=\mathbf{p}\!\nnearrow_{-1}^1$ we have $(\bar{s}(k); k\in\mathfrak{m}')=\mathbf{p}$, $(\bar{\bar{s}}(k); k\in\mathfrak{m}')=I\mathbf{p}$ and $\bar{s}(\aae)=\bar{\bar{s}}(\aae)=\emptyset$ for $\aae\in\mathfrak{m}^A\diagdown\mathfrak{m}'$. Therefore, points A. - C. are straightforward properties of the inversion $I$. Since $m'\in2\mathbb{N}$, we have $k\in2\mathbb{N}\Longleftrightarrow m'+1-k\in2\mathbb{N}-1$, which means $(I\mathbf{p})|_{\varepsilon}=I(\mathbf{p}|_{-\varepsilon})$ for $\varepsilon\in\{-1,1\}$, hence D. also holds. For $K,L\in\{1,\cdots,M\}$, points A. - D. follow immediately from the special case $K=L=1$. Condition E. is also immediate. \qed

\medskip

Since the consistency of the soul sequence $\langle\mathbf{t}^A, \mathbf{t}^A_1\rangle$ means the relation $\mathbf{t}^A_1|_{-1} \prec_{\mathbf{t}^A} \mathbf{t}^A_1|_1$, the following properties of the relation $\prec_{\mathbf{t}^A}$ will be useful.

\begin{lem} 
\label{lem3.2.6}
For sequences $\mathbf{r}, \mathbf{s}\in\mathfrak{M}^A$ we have:
\begin{enumerate}
\item[F.] (additivity rule) if $\mathbf{r}=\mathbf{r}_1\oplus\mathbf{r}_2$, $\mathbf{s}=\mathbf{s}_1\oplus\mathbf{s}_2$, $\mathbf{r}_1\perp\mathbf{r}_2$, $\mathbf{s}_1\perp\mathbf{s}_2$, then we have equivalences $\mathbf{r}\prec_{\mathbf{t}^A}\mathbf{s}\Longleftrightarrow\mathbf{r}\prec_{\mathbf{t}^A}\mathbf{s}_i$ for $i\in\{1,2\}\Longleftrightarrow\mathbf{r}_i\prec_{\mathbf{t}^A}\mathbf{s}$ for $i\in\{1,2\}$;
\item[G.] (ordering rule) if $\mathbf{s}<\mathbf{r}$, then $\mathbf{r}\prec_{\mathbf{t}^A}\mathbf{s}$;
\item[H.] (separability rule) if the sets $\bigcup\mathbf{r}$, $\bigcup\mathbf{s}$ are separated, then $\mathbf{r}\prec_{\mathbf{t}^A}\mathbf{s}$;
\item[I.] (containing rule) if $\mathbf{r}\in\mathfrak{M}(K)$, $\mathbf{s}\in\mathfrak{M}(L)$, $L-K\geqslant 2$ and we have the alternative $(\bigcup\mathbf{r}\subset t([K+1, L-1])\vee\bigcup\mathbf{s}\subset t([K+1,L-1]))$, then $\mathbf{r}\prec_{\mathbf{t}^A}\mathbf{s}$;
\item[J.] (shifting rule) for $1\leqslant J\leqslant \frac{M}{2}$; $\mathbf{p}$, $\mathbf{q}\in\mathfrak{M}'$ we have implications  $\mathbf{p}\prec_{\mathbf{t}'(2J-1))}\mathbf{q}\Longrightarrow\mathbf{p}\!\nnearrow_1^{2J-1} \prec_{\mathbf{t}^A} \mathbf{q}\!\nnearrow_1^{2J-1}$ and $\mathbf{p}\prec_{\mathbf{t}'(2J)} \mathbf{q}\Longrightarrow \mathbf{q}\!\nnearrow_{-1}^{2J} \prec_{\mathbf{t}^A} \mathbf{p}\!\nnearrow_{-1}^{2J}$.
\end{enumerate}
\end{lem}

\medskip
\noindent
\textbf{Proof.}
For $\aae, \lambda \in \mathfrak{m}^A$, $r(\aae)\neq\emptyset$, $s(\lambda)\neq\emptyset$, assumptions in points G., H., I. imply (respectively) contradiction; $r(\aae)$, $s(\lambda)$ - separated; $(r(\aae)\subset t^A([\aae+1, \lambda-1])\vee s(\lambda)\subset t^A([\aae+1,\lambda-1]))$. Thus, points G., H., I. are immediate consequences of defining the relation $\prec_{\mathbf{t}^A}$. Similarly, an immediate consequence of this definition is point F.

To show point J., we notice that for any sequences $\mathbf{n}, \mathbf{p}, \mathbf{q} \in \mathfrak{M}'$ and a number $K\in\{1,\cdots,M\}$, we have equivalences $\mathbf{p}\prec_{\mathbf{n}}\mathbf{q}\Longleftrightarrow\mathbf{p}\!\!\!\nnearrow_1^K \ \prec_{\mathbf{n}\nnearrow_1^K} \mathbf{q}\!\!\nnearrow_1^K\Longleftrightarrow\mathbf{q}\!\nnearrow_{-1}^K \ \prec_{\mathbf{n}\nnearrow_{-1}^K} \mathbf{p}\!\!\nnearrow_{-1}^K$. Thus, point J. follows from definition (\ref{eq3.2.6}) of the sequence $\mathbf{t}^A$. \qed

\medskip

We will prove the consistency of the soul sequence $\langle\mathbf{t}^A, \mathbf{t}_1^A\rangle$ using the following notations. For $K\in\{1,\cdots,M\}$, we define
\begin{namelist}{lll}
\item{$\bullet$} $p(K)=\left\{\begin{array}{rcl}
1&\text{if}&K\in2\mathbb{N}-1,\\
-1&\text{if}&K\in2\mathbb{N},
\end{array} \right.$
\end{namelist}
and then
\begin{equation}
\label{eq3.2.9}
\begin{array}{rcl}
\Theta_{-1}(K) & = & \mathbf{t}_1^{'}(K)|_{-1} \!\nnearrow_{p(K)}^K,\\
\Theta_1(K) & = & \mathbf{t}_1^{'}(K)|_1\!\nnearrow_{p(K)}^K,
\end{array}
\end{equation}
and (compare with condition C.)
\begin{equation}
\label{eq3.2.10}
\begin{array}{rcl}
\Theta_1^{0}(K) & = & (t(K)\setminus t_1(K))\otimes\Theta_1(K),\\
\Theta_1^{1}(K) & = & t_1(K)\otimes\Theta_1(K).
\end{array}
\end{equation}
Then $\Theta_{-1}(K)=\Theta_{-1}(K)|_{-p(K)}$, $\Theta_1(K)=\Theta_1^{0}(K)\oplus\Theta_1^{1}(K)$, $\Theta_1^0(K)\perp\Theta_1^1(K)$, and $\Theta_1^{\eta}(K)=\Theta_1^{\eta}(K)|_{p(K)}$ for $\eta\in\{0, 1\}$ (compare with condition D.). A crucial role will be played by the representation arising from \eqref{eq:piatka}, \eqref{eq3.2.9}:
\begin{equation*}
\begin{array}{rcl}
\mathbf{t}_1^A |_1 & = & \bigoplus\limits_{1\leqslant J\leqslant M/2} (\Theta_1(2J-1)\oplus\Theta_{-1}(2J)),\\
\mathbf{t}_1^A |_{-1} & = & \bigoplus\limits_{1\leqslant I\leqslant M/2} (\Theta_{-1}(2I-1)\oplus\Theta_1(2I)).
\end{array}
\end{equation*}
According to the additivity rule F., we have
	
\begin{cor}
\label{cor3.2.7}
The sequence $\langle \mathbf{t}^A, \mathbf{t}_1^A \rangle$ is consistent if and only if for any $I, J \in \{1, \cdots, M/2\}$, we have
\begin{equation}
\label{eq3.2.12}
\Theta_{-1}(2I-1)\oplus\Theta_1(2I)\prec_{\mathbf{t}^A}\Theta_1(2J-1)\oplus\Theta_{-1}(2J).
\end{equation}
\end{cor}

We will establish statement (\ref{eq3.2.12}) based on the following two sets of conditions (Lemmas 3.2.8 and 3.2.9). The first set is a natural consequence of conditions F. - J., but proving the second one requires new significant ideas.

\begin{lem}
\label{lem3.2.8}
For $I, J\in \{1, \cdots, \frac{M}{2}\}$, the following conditions hold:
\begin{enumerate}
\item[I)] $I<J$ implies $\Theta_{-1}(2I-1)\prec_{\mathbf{t}^A}\Theta_1(2J-1)\oplus\Theta_{-1}(2J)$, $\Theta_{-1}(2I-1)\oplus\Theta_1(2I)\prec_{\mathbf{t}^A}\Theta_{-1}(2J)$.
\item[II)] $J-I\geqslant2$ implies $\Theta_1^0 (2I) \prec_{\mathbf{t}^A} \Theta_1(2J-1)\oplus\Theta_{-1}(2J)$, $\Theta_{-1}(2I-1)\oplus\Theta_1(2I)\prec_{\mathbf{t}^A}\Theta_1^0 (2J-1)$.
\item[III)] $J-I\geqslant2$, $t_1(2I)\neq t_1(2J-1)$ implies $\Theta_1^1 (2I) \prec_{\mathbf{t}^A} \Theta_1^1 (2J-1)$.
\item[IV)] $J=I+1$ implies $\Theta_1^1 (2I) \prec_{\mathbf{t}^A} \Theta_1^1 (2J-1)$.
\item[V)] $I=J$ implies $\Theta_{-1} (2I-1) \prec_{\mathbf{t}^A} \Theta_1 (2J-1)$, $\Theta_1 (2I) \prec_{\mathbf{t}^A} \Theta_{-1} (2J)$.
\item[VI)] $I=J$ implies $\Theta_1 (2I) \prec_{\mathbf{t}^A} \Theta_1 (2J-1)$.
\item[VII)] $I>J$ implies $\Theta_{-1}(2I-1)\oplus\Theta_1(2I)\prec_{\mathbf{t}^A}\Theta_1(2J-1)\oplus\Theta_{-1}(2J)$.
\end{enumerate}
\end{lem}

\begin{lem}
\label{lem3.2.9}
For $I, J\in \{1, \cdots, \frac{M}{2}\}$, the following conditions hold:
\begin{enumerate}
\item[$\alpha$)] $J-I \geqslant 2$, $t_1 (2I) = t_1 (2J-1)$ implies $\Theta_1^1 (2I) \prec_{\mathbf{t}^A} \Theta_1^1 (2J-1)$.
\item[$\beta$)] $I=J$ implies $\Theta_{-1} (2I-1) \prec_{\mathbf{t}^A} \Theta_{-1} (2J)$.
\item[$\gamma$)] $J=I+1$ implies $\Theta_1^0 (2I) \prec_{\mathbf{t}^A} \Theta_1 (2J-1)$.
\item[$\delta$)] $J=I+1$ implies $\Theta_1 (2I) \prec_{\mathbf{t}^A} \Theta_1^0 (2J-1)$.
\end{enumerate}
\end{lem}

We will begin by demonstrating, based on Lemmas 3.2.8 and 3.2.9, that the sequence $\langle \mathbf{t}^A, \mathbf{t}_1^A \rangle$ is consistent.

\medskip
\noindent
\textbf{Proof of Theorem \ref{thm3.2.3}.}
According to Lemma \ref{lem3.2.4}, it suffices to show that the sequence $\langle \mathbf{t}^A, \mathbf{t}_1^A \rangle$ is consistent, which, by Corollary \ref{cor3.2.7}, means proving relation (\ref{eq3.2.12}) for $I, J\in \{1, \cdots, \frac{M}{2}\}$.
	
For $J-I\geqslant 2$, relation (\ref{eq3.2.12}) is implied by conditions I), II), III), and $\alpha$) (due to additivity rule F.). For $J=I+1$, relation (\ref{eq3.2.12}) similarly follows from I), IV), $\gamma$), and $\delta$). For $I=J$, relation (\ref{eq3.2.12}) follows from V), VI), and $\beta$). Therefore, in accordance with VII), the proof is concluded. \qed

\medskip
\noindent
\textbf{Proof of Lemma \ref{lem3.2.8}.}

\textbf{Proof of I).} Let $I<J$. Due to the regularity of the sequence $((t(K), t_1(K));1\leqslant K\leqslant M)$, we have $t(2I-1)=t(2I)$ and according to definition (\ref{eq3.2.9}) we have $\bigcup\Theta_{-1}(2I-1)\subset t(2I-1)$ (thus $\bigcup\Theta_{-1}(2I-1)\subset t(2I)$). Moreover, $\Theta_{-1}(2I-1)\in\mathfrak{M}(2I-1)$ while $\Theta_1(2J-1)\in\mathfrak{M}(2J-1)$ and $\Theta_{-1}(2J)\in\mathfrak{M}(2J)$. According to containing rule I. (due to additivity rule F.), we have $\Theta_{-1}(2I-1)\prec_{\mathbf{t}^A}\Theta_1(2J-1)\oplus\Theta_1(2J)$. The second relation in I) is proven similarly.

\textbf{Proof of II).}
Let $J-I\geqslant2$. This time we notice that $t(2I)\smallsetminus t_1(2I)=t(2I)\cap t(2I+1)\subset t(2I+1)$, and $\bigcup\Theta_1^0(2I)\subset t(2I)\smallsetminus t_1(2I)$ (according to (\ref{eq3.2.10})). Containing rule I. thus implies $\Theta_1^0 (2I) \prec_{\mathbf{t}^A} \Theta_1 (2J-1)\oplus\Theta_{-1}(2J)$. The second relation in II) is proven analogously.
	
\textbf{Proof of III).}
Let $J-I\geqslant2$. Since the sequence $((t(K), t_1(K)); 1\leqslant K\leqslant M)$ is consistent, the condition $t_1(2I)\neq t_1(2J-1)$ implies an alternative of conditions $1^{\circ} t_1(2I)$, $t_1(2J-1)$ - separated or $2^{\circ}$ we have an alternative $(t_1(2I)\subset t([2I+1, 2J-2])\vee t_1(2J-1)\subset t([2I+1, 2J-2]))$. According to definition (\ref{eq3.2.10}), we have $\bigcup\Theta_1^1(2I)\subset t_1(2I)$, $\bigcup\Theta_1^1(2J-1)\subset t_1(2J-1)$. In the case of $1^{\circ}$, the relation $\Theta_1^1(2I)\prec_{\mathbf{t}^A} \Theta^1 (2J-1)$ follows from separability rule H. In the case of $2^{\circ}$, it follows from containing rule I.

\textbf{Proof of IV).}
For $J=I+1$, the sets $t_1(2I)$, $t_1(2J-1)=t_1(2I+1)$ must be separated (for a consistent sequence $((t(K), t_1(K)); 1\leqslant K\leqslant M)$) and we use separability rule H.
	
Conditions VI and VII follow from ordering rule G. For example, the equality $I=J$ implies $2I>2J-1$ and we have $\Theta_{-1}(2I)\in\mathfrak{M}(2I)$, $\Theta_{-1}(2J-1)\in\mathfrak{M}(2J-1)$.
	
\textbf{Proof of V).}
Since the offspring $\langle \mathbf{t}'(2I-1), \mathbf{t}_1'(2I-1) \rangle$ of the parent $((t(2I-1), t_1(2I-1))$ is a population of souls, the sequence $\langle \mathbf{t}'(2I-1), \mathbf{t}_1'(2I-1) \rangle$ is consistent. We obtain $\mathbf{t}_1'(2I-1)|_{-1}\prec_{\mathbf{t}'(2I-1)}\mathbf{t}_1'(2I-1)|_1$ and, according to shiffing rule J., $\mathbf{t}_1'(2I-1)|_{-1}\!\nnearrow_1^{2I-1}\prec_{\mathbf{t}^A}\mathbf{t}_1'(2I-1)|_1\!\nnearrow_1^{2I-1}$. Thus, we have proved the condition $\Theta_{-1}(2I-1)\prec_{\mathbf{t}^A} \Theta_1(2I-1)$ (according to definition (\ref{eq3.2.9})). The condition  $\Theta_1(2I)\prec_{\mathbf{t}^A} \Theta_{-1}(2I)$ can be proved analogously. \qed

\medskip

We will conduct the proof of Lemma \ref{lem3.2.9}, which covers conditions $\alpha$) - $\delta$), after describing some general methods regarding the relation $\prec_{\mathbf{q}}$ for certain sequences $\mathbf{q}$ (Lemmas \ref{lem3.2.11} - \ref{lem3.2.13}). We will also use the following properties of sequences $\Theta_{-1}(K)$, $\Theta_1^0(K)$, $\Theta_1^1(K)$ (definitions (\ref{eq3.2.9}), (\ref{eq3.2.10})), resulting from the definition of offspring of souls and conditions A. and B. in Corollary \ref{cor3.2.5}.

\begin{lem}
\label{lem3.2.10}
For any integers $I, J\in\{1,\cdots,M/2\}$, we have the following relations:
\begin{enumerate}
\item[1)] $t_1(2I)=t_1(2J-1)\Longrightarrow\Theta_1^1(2I)\stackrel{-1}{\equiv} \Theta_1^1(2J-1)$;
\item[2)] $\Theta_{-1}(2J-1)\stackrel{-1}{\equiv}\Theta_{-1}(2J)$;
\item[3)] $\Theta_1^0(2I)\stackrel{-1}{\equiv}\Theta_1^0(2I+1)$ when $I\leqslant M/2-1$;
\item[4)] $\Theta_1^1(2I+1)<<\Theta_1^0(2I+1)$ when $I\leqslant M/2-1$;
\item[5)] $Theta_1^0(2I)<<\Theta_1^1(2I)$.
\end{enumerate}
\end{lem}

\medskip
\noindent
\textbf{Proof.}
According to Definition \ref{df3.2.1} of dust and antidust (and Definition \ref{df3.2.2}. of offspring of souls), for $I, J\in\{1,\cdots,\frac{M}{2}\}$, we have the conditions:
\begin{enumerate}
\item[1')] $t_1(2I)=t_1(2J-1)\Longrightarrow t_1(2I)\otimes\mathbf{t}_1'(2I)|_1\equiv t_1(2J-1)\otimes\mathbf{t}_1'(2J-1)|_1$;
\item[2')] $\mathbf{t}_1'(2J-1)|_{-1}\equiv\mathbf{t}_1'(2J)|_{-1}$ (since $t(2J-1)=t(2J)$, according to the regularity of the sequence $((t(K), t_1(K)); 1\leqslant K \leqslant M)$);
\item[3')] $(t(2I)\smallsetminus t_1(2I))\otimes\mathbf{t}_1'(2I)|_1\equiv(t(2I+1)\smallsetminus t_1(2I+1))\otimes\mathbf{t}_1(2I+1)|_1$ when $I\leqslant M/2-1$ (since $t(2I)\smallsetminus t_1(2I)=t(2I+1)\smallsetminus t_1(2I+1)$, according to the regularity of the sequence $((t(K),t_1(K));1\leqslant K\leqslant M))$;
\item[4')] $t_1(2I+1)\otimes\mathbf{t}_1'(2I+1)|_1<<(t(2I+1)\smallsetminus t_1(2I+1))\otimes\mathbf{t}_1'(2I+1)|_1$;
\item[5')] $t_1(2I)\otimes\mathbf{t}_1'(2I)|_1<<(t(2I)\smallsetminus t_1(2I))\otimes\mathbf{t}_1'(2I)|_1$.
\end{enumerate}
According to definitions (\ref{eq3.2.9}), (\ref{eq3.2.10}) of sequences $\Theta_{-1}(K)$, $\Theta_1^0(K)$, $\Theta_1^1(K)$ and implications A. and B. in Corollary \ref{cor3.2.5}, we obtain conditions 1) - 5) in the lemma being proved. \qed

\medskip

Lemma \ref{lem3.2.10} allows us to demonstrate conditions $\alpha$) - $\delta$) occurring in Lemma \ref{lem3.2.9} in a strengthened form.

For sequences $\mathbf{r}=(r(\aae); \aae\in\mathfrak{m}^A)$, $\mathbf{s}=(s(\aae); \aae \in \mathfrak{m}^A)$ belonging to $\mathfrak{M}^A$, let us denote the relation $\mathbf{r}\prec_{\oplus}\mathbf{s}$ if $\mathbf{r}\prec_{\mathbf{r}\oplus\mathbf{s}}\mathbf{s}$. Thus, for sequences $\mathbf{r}$, $\mathbf{s}$ satisfying $r(\aae)\subset t_1'(\aae)$, $s(\lambda)\subset t_1'(\lambda)$ for $\aae, \lambda\in\mathfrak{m}^A$, the relation $\mathbf{r}\prec_{\oplus}\mathbf{s}$ implies $\mathbf{r}\prec_{\mathbf{t}_1^A}\mathbf{s}$ and (a fortiori) $\mathbf{r}\prec_{\mathbf{t}^A}\mathbf{s}$. If $\tilde{\mathbf{r}},\tilde{\tilde{\mathbf{r}}}, \tilde{\mathbf{s}},\tilde{\tilde{\mathbf{s}}}\in\mathfrak{M}^A$ and $\tilde{\mathbf{r}}\perp\tilde{\tilde{\mathbf{r}}}, \tilde{\mathbf{s}}\perp\tilde{\tilde{\mathbf{s}}}$, then for $\mathbf{r}=\tilde{\mathbf{r}}\oplus\tilde{\tilde{\mathbf{r}}}$, $\mathbf{s}=\tilde{\mathbf{s}}\oplus\tilde{\tilde{\mathbf{s}}}$, $\mathbf{n}=\mathbf{r}\oplus\mathbf{s}$, we have implications

\begin{equation}
\label{eq3.2.13}
\mathbf{r}\prec_{\mathbf{n}}\tilde{\mathbf{s}}, \ \mathbf{r}\prec_{\mathbf{n}}\tilde{\tilde{\mathbf{s}}} \Longrightarrow \mathbf{r}\prec_{\oplus}\mathbf{s},
\end{equation}
	
\begin{equation}
\label{eq3.2.14}
\tilde{\mathbf{r}}\prec_{\mathbf{n}}\mathbf{s}, \ \tilde{\tilde{\mathbf{r}}}\prec_{\mathbf{n}}\mathbf{s} \Longrightarrow \mathbf{r}\prec_{\oplus}\mathbf{s}.
\end{equation}

The next three lemmas will allow us to infer conditions $\alpha$) - $\delta$) from Lemma 3.2.10, thus concluding the proof of consistency of the soul sequences $\langle \mathbf{t}^A,\mathbf{t}_1^A \rangle$.

\begin{lem}
\label{lem3.2.11}
If $\mathbf{r}, \mathbf{s}\in\mathfrak{M}^A$, $\mathbf{r}<\mathbf{s}$, and $\mathbf{r}\stackrel{-1}{\equiv}\mathbf{s}$, then $\mathbf{r}\prec_{\oplus}\mathbf{s}$.
\end{lem}
	
\medskip
\noindent
\textbf{Proof.}
For $\mathbf{r}=(r(\aae);\aae\in\mathfrak{m}^A)$, $\mathbf{s}=(s(\lambda); \lambda\in\mathfrak{m}^A)$, let the decreasing bijection $\{ \aae\in\mathfrak{m}^A : r(\aae) \neq \emptyset\}\leftrightarrow \{ \lambda\in\mathfrak{m}^A : s(\lambda)\neq\emptyset \}$ satisfy $\aae\leftrightarrow\lambda\Longrightarrow r(\aae)=s(\lambda)$. If $\aae_0$, $\lambda_0\in\mathfrak{m}^A$ satisfy $r(\aae_0) \neq \emptyset$, $s(\lambda_0)\neq\emptyset$, then we have the alternative $\aae_0\leftrightarrow\lambda_0$ or $(\exists\lambda_1\in\mathfrak{m}^A : \lambda_1<\lambda_0, \aae_0\leftrightarrow\lambda_1)$ or $(\exists\aae_1\in\mathfrak{m}^A : \aae_1>\aae_0, \aae_1\leftrightarrow\lambda_0)$. Hence, we have the disjunction $r(\aae_0)=s(\lambda_0)$ or $r(\aae_0)\subset s(\aae_0+1)\cup\cdots\cup s(\lambda_0-1)$ or $s(\lambda_0)\subset r(\aae_0+1)\cup\cdots\cup r(\lambda_0-1)$. Thus, the relation $\mathbf{r}\prec_{\oplus}\mathbf{s}$ is established.

\begin{lem}
\label{lem3.2.12}
If for sequences $\mathbf{r}$, $\tilde{\mathbf{s}}$, $\tilde{\tilde{\mathbf{s}}}\in\mathfrak{M}^A$ the relations $\mathbf{r}<\tilde{\mathbf{s}}\oplus\tilde{\tilde{\mathbf{s}}}$, $\tilde{\mathbf{s}}\perp\tilde{\tilde{\mathbf{s}}}$, $\mathbf{r}\stackrel{-1}{\equiv}\tilde{\mathbf{s}}$, $\tilde{\tilde{\mathbf{s}}}<<\tilde{\mathbf{s}}$ hold, then $\mathbf{r}\prec_{\oplus}\tilde{\mathbf{s}}\oplus\tilde{\tilde{\mathbf{s}}}$.
\end{lem}
	
\medskip
\noindent
\textbf{Proof.}
According to Lemma \ref{lem3.2.11}, we have the relation $\mathbf{r}\prec_{\oplus}\tilde{\mathbf{s}}$. Let $\mathbf{r}=(r(\aae);\aae\in\mathfrak{m}^A)$, $\tilde{\mathbf{s}}=(\tilde{s}(\lambda); \lambda\in\mathfrak{m}^A)$, $\tilde{\tilde{\mathbf{s}}}=(\tilde{\tilde{s}}(\lambda); \lambda\in\mathfrak{m}^A)$, and let the decreasing bijection $\{ \aae\in\mathfrak{m}^A : r(\aae) \neq \emptyset\}\leftrightarrow \{ \lambda\in\mathfrak{m}^A : \tilde{s}(\lambda)\neq\emptyset \}$ satisfy $\aae\leftrightarrow\lambda\Longrightarrow r(\aae)=\tilde{s}(\lambda)$. Fix $\aae_0, \lambda_0\in\mathfrak{m}^A$ such that $r(\aae_0)\neq\emptyset$, $\tilde{\tilde{s}}(\lambda_0)\neq\emptyset$. Since $\tilde{s}\perp\tilde{\tilde{s}}$, for $\lambda_1\in\{\lambda\in\mathfrak{m}^A:\tilde{s}(\lambda)\neq\emptyset\}$ satisfying $\aae_0\leftrightarrow\lambda_1$, we have $\lambda_1\neq\lambda_0$.

If $\lambda_1>\lambda_0$, then the sets $\tilde{\tilde{s}}(\lambda_0)$ and $\tilde{s}(\lambda_1)$ are separated (due to $\tilde{\tilde{s}}<<\tilde{s}$), hence $\tilde{\tilde{s}}(\lambda_0)$ and $r(\aae_0)$ are separated. If $\lambda_1<\lambda_0$, then $r(\aae_0)=\tilde{s}(\lambda_1)\subset\tilde{s}(\aae_0+1)\cup\cdots\cup\tilde{s}(\lambda_0-1)$. Thus, we have shown that $\mathbf{r}\prec_{\tilde{\mathbf{s}}}\tilde{\tilde{\mathbf{s}}}$ and (a fortiori) $\mathbf{r}\prec_{\mathbf{n}}\tilde{\tilde{\mathbf{s}}}$ for $\mathbf{n}=\mathbf{r}\oplus\tilde{\mathbf{s}}\oplus\tilde{\tilde{\mathbf{s}}}$. Considering the observation $\mathbf{r}\prec_{\oplus}\tilde{\mathbf{s}}$, we have $\mathbf{r}\prec_{\mathbf{n}}\tilde{\mathbf{s}}$, and thus, due to the relation $\tilde{\mathbf{s}}\perp\tilde{\tilde{\mathbf{s}}}$, $\mathbf{r}\prec_{\oplus}\tilde{\mathbf{s}}\oplus\tilde{\tilde{\mathbf{s}}}$, in accordance with (\ref{eq3.2.13}).

\begin{lem}
\label{lem3.2.13}
If for sequences $\tilde{\mathbf{r}}$, $\tilde{\tilde{\mathbf{r}}}$, $\mathbf{s}\in\mathfrak{M}^A$ the relations $\tilde{\mathbf{r}}\perp\tilde{\tilde{\mathbf{r}}}$, $\tilde{\mathbf{r}}\oplus\tilde{\tilde{\mathbf{r}}}<\mathbf{s}$, $\tilde{\mathbf{r}}\stackrel{-1}{\equiv}\mathbf{s}$, $\tilde{\mathbf{r}}<<\tilde{\tilde{\mathbf{r}}}$ hold, then $\tilde{\mathbf{r}}\oplus\tilde{\tilde{\mathbf{r}}}\prec_{\oplus}\mathbf{s}$.
\end{lem}

\medskip
\noindent
\textbf{Proof.}
Again, according to Lemma \ref{lem3.2.11}, we have $\tilde{\mathbf{r}}\prec_{\oplus}\mathbf{s}$. Let $\tilde{\mathbf{r}}=(\tilde{r}(\aae); \aae\in\mathfrak{m}^A)$, $\tilde{\tilde{\mathbf{r}}}=(\tilde{\tilde{r}}(\aae); \aae\in\mathfrak{m}^A)$, $\mathbf{s}=(s(\lambda); \lambda\in\mathfrak{m}^A)$, and let the decreasing bijection $\{\aae\in\mathfrak{m}^A:\tilde{r}(\aae)\neq\emptyset\}\leftrightarrow\{\lambda\in\mathfrak{m}^A:s(\lambda)\neq\emptyset\}$ satisfy $\aae\leftrightarrow\lambda\Longrightarrow\tilde{r}(\aae)=s(\lambda)$. This time, we consider the case when there exist $\aae_0, \lambda_0\in\mathfrak{m}^A$ such that $\tilde{\tilde{r}}(\aae_0)$, $s(\lambda_0)\neq\emptyset$. We find $\aae_1\in\{\aae\in\mathfrak{m}^A:\tilde{r}(\aae)\neq\emptyset\}$ satisfying $\aae_1\leftrightarrow\lambda_0$ (and $\aae_1\neq\aae_0$). If $\aae_1<\aae_0$, then the sets $\tilde{r}(\aae_1)$, $\tilde{\tilde{r}}(\aae_0)$ are separated (due to $\tilde{\mathbf{r}}<<\tilde{\tilde{\mathbf{r}}}$), thus $\tilde{\tilde{r}}(\aae_0)$, $s(\lambda_0)$ are separated, as $\tilde{r}(\aae_1)=s(\lambda_0)$. If $\aae_1>\aae_0$, then $s(\lambda_0)=\tilde{r}(\aae_1)\subset\tilde{r}(\aae_0+1)\cup\cdots\cup\tilde{r}(\lambda_0-1)$. Hence, we have shown that $\tilde{\tilde{\mathbf{r}}}\prec_{\tilde{\mathbf{r}}}\mathbf{s}$. For $\mathbf{n}=\tilde{\mathbf{r}}\oplus\tilde{\tilde{\mathbf{r}}}\oplus\mathbf{s}$, we have $\tilde{\tilde{\mathbf{r}}}\prec_{\mathbf{n}}\mathbf{s}$ and $\tilde{\mathbf{r}}\prec_{\mathbf{n}}\mathbf{s}$ (due to $\tilde{\mathbf{r}}\prec_{\oplus}\mathbf{s}$), hence $\tilde{\mathbf{r}}\oplus\tilde{\tilde{\mathbf{r}}}\prec_{\oplus}\mathbf{s}$, according to (\ref{eq3.2.14}).

\medskip
\noindent
\textbf{Proof of Lemma \ref{lem3.2.9}.}
We will use points 1) - 5) in Lemma \ref{lem3.2.10}.

Lemma \ref{lem3.2.11}, applied for $\mathbf{r}=\Theta_{-1}(2I-1)$, $\mathbf{s}=\Theta_{-1}(2I)$, according to condition 2) in Lemma \ref{lem3.2.10} and the condition $\mathbf{r}<\mathbf{s}$ (following from $\mathbf{r}\in\mathfrak{M}(2I-1)$, $\mathbf{s}\in\mathfrak{M}(2I)$), yields $\Theta_{-1}(2I-1)\prec_{\oplus}\Theta_{-1}(2I)$, hence we have $\beta$).

Similarly, Lemma \ref{lem3.2.11} applied for $\mathbf{r}=\Theta_1^1(2I)$, $\mathbf{s}=\Theta_1^1(2J-1)$, $J-I\geqslant2$, according to condition 1) in Lemma \ref{lem3.2.10}, gives $\alpha$).

Lemma \ref{lem3.2.12}, applied for $\mathbf{r}=\Theta_1^0(2I)$, $\tilde{\mathbf{s}}=\Theta_1^0(2I+1)$, $\tilde{\tilde{\mathbf{s}}}=\Theta_1^1(2I+1)$ (when $I\leqslant M/2-1$), with conditions 3), 4) in Lemma \ref{lem3.2.10} and the relation $\mathbf{r}<\tilde{\mathbf{s}}\oplus\tilde{\tilde{\mathbf{s}}}$ (resulting from $\mathbf{r}\in\mathfrak{M}(2I)$; $\tilde{\mathbf{s}}, \tilde{\tilde{\mathbf{s}}}\in\mathfrak{M}(2I+1)$), gives condition $\gamma$).

Lemma \ref{lem3.2.13}, applied for $\tilde{\mathbf{r}}=\Theta_1^0(2I)$, $\tilde{\tilde{\mathbf{r}}}=\Theta_1^1(2I)$, $\mathbf{s}=\Theta_1^0(2I+1)$ (when $I\leqslant M/2-1$), with conditions 3) and 5) in Lemma \ref{lem3.2.10} and the condition $\tilde{\mathbf{r}}\oplus\tilde{\tilde{\mathbf{r}}}<\mathbf{s}$ (resulting from $\tilde{\mathbf{r}}, \tilde{\tilde{\mathbf{r}}}\in\mathfrak{M}(2I)$, $\mathbf{s}\in\mathfrak{M}(2I+1)$), gives condition $\delta$). \qed

\medskip

Thus, we have completed the proof of Theorem \ref{thm3.2.3}. We have shown that if the sequence of parents constitutes the soul population, then the anti-order of the sequence of offspring of these parents is also a soul population.

\section{Offspring of strictly convex souls}

We will introduce families of souls $\tau,\bar{\tau}\subset\tau_0$, for which we will further construct the offspring of souls 
$\tau\ni(t,t_1)\mapsto\left\langle \bs{t}'(t,t_1),\bs{t}(t,t_1)\right\rangle\in\bar{\tau}^{m'}$, described in Definition \ref{df3.2.2}.

We shall fix the following notation. As usually, for $\epsilon>0$, $B(\bs{x},\varepsilon)$ is the closed disc in $\R^2$ centered at $\bs{x}$, of radius, and
\begin{namelist}{lll}
\item{$\bullet$} $\mathcal{B}_\varepsilon = \left\{B(\bs{x},\varepsilon)\colon\bs{x}\in\R^2\right\}$.
\item{$\bullet$} For $t\in\mathcal{T}$, $\varrho\geq 1$, we shall say that $t$ is \textit{$\varrho$-convex at the point $\bs{x}$} if $\bs{x}\in\partial t$ and
\[
\exists B\in\mathcal{B}_\varrho,\varepsilon>0\colon {\bs{x}\in\partial B\text{, }B(\bs{x},\varepsilon)\cap t \subset B} \text{.}
\]
\end{namelist}
For the fixed set \(\bbt\), we assume that \(\bbt\subset[-1,1]\times[-1,1]\), \(d(\bbt)\leq 1\), and define
\begin{namelist}{lll}
\item{$\bullet$} $
\mathcal{F}_\varrho=\left\{
t\in\mathcal{T}\colon t\subset \bbt\text{ and $t$ is $\varrho$-convex at every }\bs{x}\in(\partial t)\cap\operatorname{int}(\bbt)
\right\}\text.$
\end{namelist}
Then it is clear that $\mathcal{F}_\varrho\subset \mathcal{F}_{\bar{\varrho}}$ if $\varrho<\bar{\varrho}$.
For any set $t\in\mathcal{T}$ and any sequence of sets $\bs{t}=\left(t(k); 1\leq k\leq m\right)\in\mathcal{T}^m$ we denote
\begin{equation*}
	d_x(t)= \sup{\left\{ \lvert x-\bar{x} \rvert\colon (x,y),(\bar{x},\bar{y})\in t \right\}}\text,
\end{equation*}
\begin{equation}\label{eq:3.3.1}
	d_x(\bs{t})=\max{\left\{ d_x(t(k))\colon 1\leq k\leq m\right\}}\text.
\end{equation}
For $\varrho\geq 1$ and $\sigma>0$, we define the families of souls
\begin{namelist}{lll}
\item{$\bullet$} $\tau_\varrho = \left\{ (t,t_1)\in\tau_0\colon t,t\setminus t_1\in\mathcal{F}_\varrho \right\} \text,
	\tau_{\varrho\sigma} = \left\{ (t,t_1)\in\tau_\varrho\colon d(t_1)\leq \sigma \right\} \text.$
\end{namelist}

Throughout this section the numbers $\beta\geq 1$, $\gamma,\gamma'\in (0,1/4)$, and $\beta'=\max{\left\{2\beta,\gamma+1/\gamma \right\}}$ will be fixed. Now we proceed to the proof of the following theorem, fundamental in our construction.

\begin{thm}
\label{thm3.3.1}
\label{thm:fund}
	For the fixed numbers \(\beta\), \(\gamma\), \(\beta'\), \(\gamma'\), there exists a number $m'\in2\mathbb{N}$, independent of \(\bbt\) (which satisfies \(\mathcal{T}\ni\bbt\subset[-1,1]\times[-1,1]\), \(d(\bbt)\leq 1\)), such that one can find offspring of souls
	\[
	\tau_{\beta\gamma}\ni(t,t_1)\mapsto\left\langle \bs{t}'(t,t_1),\bs{t}'_1(t,t_1)\right\rangle
	\in\left(\tau_{\beta'\gamma'}\right)^{m'}
	\]
	satisfying
	\begin{equation}\label{eqn:dx11gam-1}
		d_x(\bs{t}'(t,t_1)) \leq 11\gamma \text.
	\end{equation}
\end{thm}

Let us first note that an equivalent formulation of the above Theorem \ref{thm:fund} can be given. For any set $t\in\mathcal{T}$ and any sequence of sets $\bs{t}=\left(t(k); 1\leq k\leq m\right)\in\mathcal{T}^m$ let us denote

\begin{equation*}
I_{xy}t = \left\{(y,x)\colon (x,y)\in t\right\}\text,
\end{equation*}
\begin{equation*}
d_y(t) = d_x(I_{xy}t)\text,
\end{equation*}
\begin{equation}\label{eq:3.3.3}
d_y(\bs{t}) = \max{\left\{d_y(t(k)\colon 1\leq k\leq m)\right\} } \text.
\end{equation}

\begin{cor}
\label{cor3.3.2}
	Theorem \ref{thm:fund} remains true if condition \eqref{eqn:dx11gam-1} is replaced with
	\[
	d_y(\bs{t}'(t,t_1)) \leq 11\gamma\text,\quad\text{for }(t,t_1)\in\tau_{\beta\gamma}\text.
	\]
\end{cor}

\medskip
\noindent
\textbf{Proof.}
	Let us denote
	\[\mathcal{F}_\varrho^{xy}=\left\{t\in\mathcal{T}\colon t\subset I_{xy}\bbt\text{ and }t\text{ is \(\varrho\)-convex at every }\bs{x}\in(\partial t)\cap\operatorname{int}{(I_{xy}\bbt)}\right\}\text,\]
	and
	\(\tau_\varrho^{xy} = \left\{(t,t_1)\in\tau_0\colon t,t\sm t_1\in\mathcal{F}_\varrho^{xy}\right\}\),
	\(\tau_{\varrho\sigma}^{xy}=\left\{(t,t_1)\in\tau_\varrho^{xy}\colon d(t_1)\leq\sigma\right\}\),
	for \(\varrho\geq 1\), \(\sigma>0\).
	
	Together with the proof of Theorem \ref{thm:fund} we have a proof of the following fact.
	There exists an offspring of souls
	\[\tau_{\beta\gamma}^{xy}\ni(\bar{t},\bar{t}_1)\mapsto \left\langle\bs{t}^{xy}(\bar{t},\bar{t}_1),\bs{t}^{xy}_1(\bar{t},\bar{t}_1)\right\rangle \in \left(\tau_{\beta'\gamma'}^{xy}\right)^{m'}\]
	satisfying \(d_x(\bs{t}^{xy}(\bar{t},\bar{t}_1))\leq 11\gamma\), for \((\bar{t},\bar{t}_1)\in\tau_{\beta\gamma}^{xy}\).
	
For sequence of sets $\bs{t}=\left(t(k); 1\leq k\leq m\right)\in\mathcal{T}^m$ let us denote \(I_{xy}\bs{t} = \left(I_{xy}t(k); 1\leq k\leq m\right)\), then for any sequence of souls \( \langle\bs{t},\bs{t}_1\rangle \), we have the sequence of souls \(\left\langle I_{xy}\bs{t},I_{xy}\bs{t}_1\right\rangle \).

Then, the function
\[
\tau_{\beta\gamma}\ni(t,t_1)\mapsto\left\langle I_{xy}\bs{t}^{xy}(I_{xy}t,I_{xy}t_1),I_{xy}\bs{t}^{xy}_1(I_{xy}t,I_{xy}t_1)\right\rangle
\in\left(\tau_{\beta'\gamma'}\right)^{m'}
\]
satisfies the requirements of Definition \ref{df3.2.2}, thus it is offspring of souls.
Additionally, for \( (t,t_1)\in\tau_{\beta\gamma} \), we have
\[
	d_y(I_{xy}\bs{t}^{xy}(I_{xy}t,I_{xy}t_1)) \leq 11\gamma
\text,
\]which completes the proof. \qed

\subsubsection*{A. Construction of offspring of souls with the use of stations}

We will reduce the proof of Theorem \ref{thm:fund} to a construction of some increasing sequences of sets. In this section we use the following notation in which the numbers $\beta\geq 1$, $\gamma,\gamma'\in(0,1/4)$, \( \beta' = \max{\left\{2\beta, \gamma + 1/\gamma\right\}}\) and the set \(\bbt\) are assumed to be fixed.
\begin{df}\label{def:3.3.3}
	\textbf{A.} For any set \( t\subset\mathbb{R}^2\), sequence \( \bs{t} = \left(t(l);1\leq l\leq n\right)\) will be called a \emph{net} or a \emph{net of length $n$ for $t$} if \(t(1)\subset\ldots\subset t(n)=t\) and \(
	\left(t(l), t(l)\setminus t(l-1)\right)\in\tau_{\beta'\gamma'}
	\), for \(2\leq l\leq n\); with \(\mathcal{N}_n^t\) we will denote the space of such nets.
	
	\textbf{B.} If \(t_1\subset\mathbb{R}^2\), we will say that \(\bs{t}^\ast = \left( t^\ast(l);1\leq l\leq n^\ast \right) \) is a \emph{station} or a \emph{station of length $n^\ast$ for $t_1$}, if $\emptyset=t^\ast(1)\subset\ldots\subset t^\ast(n^\ast)=t_1$ and for every set $t$ satisfying \( (t,t_1)\in\tau_{\beta'/2}\) the sequence \( \left( (t\setminus t_1)\cup t^\ast(l); 1\leq l\leq n^\ast \right) \) is a net (of length $n^\ast$ for $t$); with \(\mathcal{S}_{n^*}^{t_1}\) we will denote the space of such stations.
\end{df}

Let us emphasize that the families of nets \(\mathcal{N}_n^t\) and stations \(\mathcal{S}_n^{t_1}\) are defined for the fixed numbers \(\beta'\) and \(\gamma'\).
In the present part A of this section we will prove Theorem \ref{thm:fund} assuming that the following theorem is valid.

\begin{thm}
\label{thm:withast}
	For the fixed numbers \(\beta'\), \(\gamma'\), there exists a number \(n^*\in\mathbb{N}\), independent of \(\bbt\) (satisfying \(\mathcal{T}\ni\bbt\subset[-1,1]\times[-1,1]\), \(d(\bbt)\leq 1\)), such that the following condition is valid:
	\begin{itemize}
		\item[{\normalfont(\textasteriskcentered)}] for \((t,t_1)\in\tau_{\beta'/2}\) the set \(\mathcal{S}_{n^*}^{t_1}\) is not empty.
	\end{itemize}
\end{thm}

The number \(n^*\) is considered fixed throughout part A of this section.

For any \(t\in\mathcal{T}\) and $n \in \N$ we denote 
\begin{namelist}{lll}
\item{$\bullet$} \( P_xt=\left\{x\in\R\colon \exists y\colon (x,y)\in t\right\} \), that is, the projection of \(t\) onto the \(x\)-axis.

\item{$\bullet$} A sequence of numbers \( \boldsymbol{\chi} = \left(\chi(l);1\leq l\leq n\right) \) will be called a \textit{skeleton} or a \textit{skeleton of length $n$} if \( 0\leq\chi(l)-\chi(l-1)\leq\gamma \), for \( 2\leq l\leq n \). If \(t\in\mathcal{T}\), such skeleton of length $n$ will be called a \textit{skeleton for \(t\)} whenever \( P_xt = \left[\chi(1),\chi(n)\right] \).
\end{namelist}

The symbol \(\mathcal{R}_n\), or \(\mathcal{R}_n^t\), will denote the space of all skeletons of length \(n\), or skeletons of length \(n\) for \(t\) respectively.

For \( t\in\mathcal{F}_{\beta'/2}\), \(\bs{t}=\left(t(l);1\leq l\leq n\right)\in\mathcal{N}_n^t\),
\( \boldsymbol{\chi} = \left(\chi(l);1\leq l\leq n\right) \in\mathcal{R}_n\); \(i,j\in\mathbb{Z}\), we will write 
\begin{namelist}{lll}
\item{$\bullet$} \(\bs{\chi}\prec_{ij}^t\bs{t}\) when $t\cap \left(x\leq\chi(l)+i\gamma\right) \subset t(l)\subset \left(x\leq\chi(l)+j\gamma\right)\text{, for any }1\leq l\leq n\text.$
\end{namelist}
Let us emphasize that in this definition \(\bs{\chi}\) is not necessarily a skeleton for \(t\), although it will often be the case in arguments of the paper.

We will make use of the following elementary fact.
\begin{re}
\label{rem:elem_tcapU}
	If \(t\in\mathcal{T}\) and \(U\) is an open set for which \(t\cap U\neq\emptyset\), then \(\operatorname{int}{(t)}\cap U\neq \emptyset\).
\end{re}
\medskip
\noindent
\textbf{Proof.}
	For any point \(\bs{x}\in t\) we have \( \bs{x}\in\clo\left[\bigcup_{\varepsilon\in(0,1)}\left((1-\varepsilon)\bs{x} +\varepsilon\operatorname{int}{(t)}\right)\right] \subset \clo{\left(\operatorname{int}{(t)}\right)}\) by convexity of \(t\). \qed

In the following the assumption that \(\gamma\in\left(0,1/4\right)\) is essential.
We also fix \(n_0\geq 1/\gamma+1\), so, for any \(t\in\mathcal{F}_\beta\), there exists a skeleton \(\bs{\chi}\in\mathcal{R}^t_{n_0}\), as \(d(t)\leq d(\bbt)\leq 1\).

\begin{lem}
\label{lem:withalpha}
	For \(n_1=(n_0-1)n^*+1\) the following condition is valid
	\begin{itemize}
		\item[(\(\alpha\))] for \(t\in\mathcal{F}_{\beta}\), there exist a net \(\bs{t}^+ \in\mathcal{N}_{n_1}^t\) and a skeleton \(\bs{\chi}^+\in\mathcal{R}_{n_1}^t\) satisfying \( \bs{\chi}^+ \prec_{3,5}^t\bs{t}^+ \).
	\end{itemize}
\end{lem}

The number \(n_1\) will be fixed throughout part A of this section.

\medskip
\noindent
\textbf{Proof of Lemma \ref{lem:withalpha}}
	Let us denote \( \varrho = \gamma/2 +1/{2\gamma}\) and a skeleton \( \bs{\chi} = \left(\chi(l); 1\leq l\leq n^0\right) \in\mathcal{R}_{n^0}^t \).
	For the ball \( B_l:= B\left((\chi(l)+4\gamma-\varrho, 0), \varrho\right)\) and the square \( S=[-1,1]\times[-1,1]\), in view of \(\gamma \in(0,1/4)\), we have \( S\cap \left(x\leq \chi(l)+3\gamma\right)\subset S\cap B_l \subset S\cap \left(x\leq\chi(l) +4\gamma\right)\). Thus, for \( t(l):= t\cap B_l \), \( l \in \left\{ 1,\ldots , n^0\right\}\), we obtain \( t\cap \left(x\leq\chi(l)+3\gamma\right) \subset t(l)\subset \left(x\leq\chi(l)+4\gamma\right) \). Consequently, we have shown that \(\bs{\chi}\prec_{3,4}^t\bs{t}\), for \(\bs{t}=\left(t(l);1\leq l\leq n^0\right)\).
	According to Remark \ref{rem:elem_tcapU}, for \( U=\left(x<\chi(l)+3\gamma\right)\), the sets \( t(l) \) have nonempty interiors, and, by the definition of \( B_l\), are convex and satisfy \( t(l) \in \mathcal{F}_{\max\{\beta,\varrho\}} =  \mathcal{F}_{\beta'/2}\).
	
	For \( 2 \leq l \leq n^0\), we have \( (t(l),t(l)\sm t(l-1))\in\tau_{\beta'/2}\) and, for the number \(n^\ast\) defined in Theorem \ref{thm:withast}, there exists a station \( \bs{t}^\ast(l)=\left(t^*(l,k); 1\leq k\leq n^*\right)\) of length \(n^*\) for \( t(l)\setminus t(l-1)\). Consider the concatenations:
	\begin{eqnarray*}
		\bs{t}^+ &=&
		\left(t(1)\cup t^*(2,k); 1\leq k\leq n^*\right) \circ \ldots\circ
		\left(t(n^0-1)\cup t^*(n^0,k); 1\leq k\leq n^*\right)\circ \left(t(n^0)\right) \text,\\
		\bs{\chi}^+ &=&
		\underbrace{\left(\chi(1),\ldots,\chi(1)\right)}_{n^*} \circ \ldots\circ
		\underbrace{\left(\chi(n^0-1),\ldots,\chi(n^0-1)\right)}_{n^*}\circ\left((\chi(n^0))\right)
	\end{eqnarray*}
	(where \(\left(t(n^0)\right)\) and \(\left(\chi(n^0)\right)\) are one-element sequences).
	Then we have \(\bs{\chi}^+\in\mathcal{R}_{n_1}^t\).
	As
	\( t(l-1)\cup t^*(l,k)\subset t(l)\subset t\cap(x\leq \chi(l)+4\gamma)\subset t\cap (x\leq\chi(l-1)+5\gamma) \), we can see that \( \bs{\chi}^+\prec_{3,5}^t\bs{t}^+\). By definition of the station \(\bs{t}^*(l)\), elements of the sequence \(\bs{t}^+\) lie in \(\mathcal{F}_{\beta'}\) and \(\bs{t}^+\in\mathcal{N}_{n_1}^t\). \qed

\medskip

Now we will consider a mirror version of the condition (\(\alpha\)) in the presented lemma. 
\begin{namelist}{lll}
\item{$\bullet$} For \(t\in\mathcal{F}_\beta\), any sequence \(\bs{t} = \left(t(k); 1\leq k\leq n\right)\) will be called an \textit{anti-net} or an \textit{anti-net of length $n$ for \(t\)} if \(t=t(1)\supset\ldots\supset t(n)\) and \( \left(t(l),t(l)\setminus t(l+1)\right) \in \tau_{\beta'\gamma'} \), for \(1\leq l\leq n-1\). The symbol \(\mathcal{A}_n^t\) will denote the space of such anti-nets.
\end{namelist}

For \(\bs{t}=\left(t(l);1\leq l\leq n\right)\in\mathcal{A}_n^t\) and \( \boldsymbol{\chi} = \left(\chi(l);1\leq l\leq n\right) \in\mathcal{R}^t_n\) we will write 
\begin{namelist}{lll}
\item{$\bullet$} \(\bs{t}\succ_{-3,-5}^t\bs{\chi}\) when $ t\cap \left(x\geq\chi(l)-3\gamma\right) \subset t(l)\subset \left(x\geq\chi(l)-5\gamma\right)\text{, for any }1\leq l\leq n\text.$
\end{namelist}
We stress that \(\bs{\chi}\) is a skeleton for \(t\).

\begin{cor}
\label{cor:withbeta}
	For \(n_1=(n_0-1)n^*+1\) we also have:
	\begin{itemize}
		\item[(\(\beta\))] for any \(t\in\mathcal{F}_\beta\) there exists an anti-net \(\bs{t}^-\in\mathcal{A}_{n_1}^t\) and a skeleton \(\bs{\chi}^-\in\mathcal{R}_{n_1}^t\) satisfying \(\bs{t}^-\succ_{-3,-5}^t \bs{\chi}^-\).
	\end{itemize}
\end{cor}
\medskip
\noindent
\textbf{Proof.}
	In this proof for set \(t\subset\mathbb{R}^2\), sequence of sets \( \bs{t} = \left(t(l); 1\leq l\leq n\right)\) in \(\mathbb{R}^2\) and sequence of numbers \( \bs{\chi} = \left(\psi(l); 1\leq l\leq n\right) \) we will denote
	\begin{eqnarray*}
		I_x t &=& \left\{ (-x,y)\colon (x,y\in t) \right\}\text,\\
		I_x \bs{t} &=& \left(I_xt(n+1-l); 1\leq l\leq n\right)\text,\\
		I_-\bs{\chi} &=& \left( -\chi(n + 1 - l ); 1\leq l\leq n\right)\text.
	\end{eqnarray*}
	
	For the fixed set \(\bbt\) and arbitrary numbers \(\varrho\geq 1\), \(\sigma>0\) we set
	\begin{eqnarray*}
	\mathcal{F}_\varrho^x &=& \left\{ t\in\mathcal{T}\colon t\subset I_x\bbt
		\text{, $t$ is $\varrho$-convex at every }
		\bs{x}\in\partial t\cap\operatorname{int}(I_x\bbt)\right\}\text,\\
	\tau_\varrho^x &=& \left\{(t,t_1)\in\tau_0\colon t,t\sm t_1\in\mathcal{F}_\varrho^x\right\}\text,\\
	\tau_{\varrho\sigma}^x &=& \left\{(t,t_1)\in\tau_\varrho^x\colon d(t_1)\leq\sigma\right\}\text,
	\end{eqnarray*}
and for the fixed $\beta', \gamma'$ and any sets $t,t_1$ and number $n \in \N$,

\begin{itemize}
\item[] \( \left(t(l); 1\leq l\leq n\right)\in\mathcal{N}^t_{n,x}\), if \(t(1)\subset\ldots\subset t(n)=t\), \\
\(\left(t(l),t(l)\sm t(l-1)\right)\in\tau^x_{\beta'\gamma'}\), for \(2\leq l\leq n\);
\item[] \(\left(t^*(l); 1\leq l\leq n\right)\in\mathcal{S}^{t_1}_{n,x}\), if for any set \(\bar{t}\) satisfying \((\bar{t},t_1)\in\tau_{\beta'/2}^x\), we have \\
$\left((\bar{t}\sm t_1)\cup t^*(l); 1\leq l\leq n\right) \in \mathcal{N}^{\bar{t}}_{n,x}$
\end{itemize}
	
	For the previously fixed number \(n_1=(n_0-1)n^*+1\), we thus have the equivalences:
	\begin{eqnarray*}
		\bs{\chi}\in\mathcal{R}_{n_1}^{I_xt} &\iff& I_{-}\bs{\chi}\in\mathcal{R}^t_{n_1}\text, \\
		\bs{t}^+\in\mathcal{N}^{I_xt}_{n_1,x} &\iff& I_x\bs{t}^+\in\mathcal{A}^t_{n_1}\text,
	\end{eqnarray*}
	and, for \(\bs{\chi}\in\mathcal{R}_{n_1}^{I_xt}\), \(t\in\mathcal{N}_{n_1,x}^{I_xt}\),
	\[
	\bs{\chi}\prec_{3,5}^{I_xt} \bs{t} \iff I_x\bs{t}^+\succ_{-3,-5}^tI_-\bs{\chi} \text.
	\]
	
Additionally, Theorem \ref{thm:withast} applied with \(I_x\bbt\) in place of \(\bbt\) yields that $\text{for }(t,t_1)\in\tau_{\beta'/2}^x\text{ the set }\mathcal{S}^{t_1}_{n^*,x}\text{ is not empty.}$
	Therefore, the statement of Lemma \ref{lem:withalpha} holds for \(I_x\bbt\), in place of \(\bbt\).
	It means that, for \(n_1=(n_0-1)n^*+1\) we have
\begin{namelist}{lll}
\item{$(\alpha_x)$} for any set \(\bar{t}\in\mathcal{F}_\beta^x\) there exists a skeleton \(\bs{\chi}_{\bar{t}}\in\mathcal{R}_{n_1}^{\bar{t}}\) and a net \(\bs{t}_{\bar{t}}\in\mathcal{N}_{n_1,x}^{\bar{t}}\) satisfying
	\(\bs{\chi}_{\bar{t}}\prec_{3,5}^{\bar{t}}\bs{t}_{\bar{t}}\).
\end{namelist}
	
For any fixed set \(t\in\mathcal{F}_\beta\), with the assumed notation, it is enough to take \(\bs{t}^-=I_x\bs{t}_{I_xt}\) and \(\bs{\chi}^-=I_-\bs{\chi}_{I_xt}\). \qed

\medskip

The proof of Theorem \ref{thm:fund} that will follow from Theorem~\ref{thm:withast}, Lemma~\ref{lem:withalpha}, and Corollary \ref{cor:withbeta} is still quite complicated. Although the lemma below may look technical, it reveals an outline of the reasoning behind our proof of the theorem.
Let us recall that for sequences of sets \(\bs{r} = \left(r(l); 1\leq l\leq n\right)\), \(\bs{s} = \left(s(l); 1\leq l\leq n\right)\) and \(\bs{p} = \left(p(l); 1\leq l\leq m\right)\), the relations \(\bs{s}|_1=\bs{r}\), \(\bs{s}|_{-1}=\bs{r}\) imply
\[
r(l) = \begin{cases} s(l), & \text{for } l\in2\mathbb{N}-1 \\ \emptyset, &\text{for }l\in2\mathbb{N} \end{cases} \text,\quad
r(l) = \begin{cases} s(l), & \text{for } l\in2\mathbb{N} \\ \emptyset, &\text{for }l\in2\mathbb{N}-1 \end{cases} \text,
\]
respectively. The relation \(\bs{s}\ssucc\bs{p}\), in turn, implies that \(\bs{p}\) coincides with 
subsequence \(\left(s(l(k)); 1\leq k\leq m\right)\) of \(\bs{s}\), and \(s(l)=\emptyset\), for \(l\in\{1,\ldots,n\}\sm\{l(1),\ldots,l(m)\}\). Moreover, \(\bs{J}\bs{s}=s(1)\), \(\bs{S}\bs{s}=s(n)\) are the first and the last elements of the sequence \(\bs{s}\).
\begin{lem}\label{lem:osmy_IstPopDusz}
	If a set \(t\in\mathcal{F}_{\beta}\), a number \(n\in\mathbb{N}\), a net
	\(\bs{t}^+\in\mathcal{N}_n^t\), an anti-net \(\bs{t}^-\in\mathcal{A}_n^t\), and a skeleton
	\(\bs{\chi}\in\mathcal{R}_n^t\) satisfy \(\bs{t}^-\succ_{-3,-5}^t\bs{\chi}\) and
	\(\bs{\chi}\prec_{-1,6}^t\bs{t}^+\), then there exists a population of souls
	\(\left\langle\bs{t}',\bs{t}_1'\right\rangle\in\left(\tau_{\beta'\gamma'}\right)^{2n}\) for which
	\begin{itemize}
		\item[\textup{1)}] \(\bs{t}'_1|_1 \ssucc \Delta_1\bs{t}^+\),\quad  \(\bs{t}'_1|_{-1}\ssucc\Delta_{-1}\bs{t}^-\),
		\item[\textup{2)}] \(\bigcup\bs{t}'=t\),\quad \(d_x(\bs{t}')\leq 11\gamma\),
		\item[\textup{3)}] \(\bs{J}\bs{t}'=\bs{J}\bs{t}^+\),\quad \(\bs{S}\bs{t}'=\bs{S}\bs{t}^-\).
	\end{itemize}
\end{lem}

In our reasoning we will make use of the following elementary fact.
\begin{re}\label{rem:element_zbiry}
	For any sets \(A\), \(B\), \(C\) and \(D\) we have the implication
	\[
{A\setminus B\subset C\subset D} \implies {A\setminus B = (A\cap C)\setminus (B\cap D)}\text.
	\]
\end{re}

\medskip

\noindent
\textbf{Proof.}
	For any \(\bs{x}\in A\setminus B\) we have \(\bs{x}\in A\), \(\bs{x}\in C\) and \(\bs{x}\notin B\), thus \(\bs{x}\notin B\cap D\). On the other hand, for \(\bs{x}\in (A\cap C)\setminus (B\cap D)\) we have \(\bs{x}\in A\), as well as \(\bs{x}\in C\subset D\) and \(\bs{x}\notin B\cap D\), hence \(\bs{x}\in A\), \(\bs{x}\notin B\). \qed

\medskip

\noindent
\textbf{Proof of \textup{Lemma \ref{lem:osmy_IstPopDusz}}}
	Let us denote
	\[
	\bs{t}^\otimes=\left(t^\otimes(l); 1\leq l\leq n\right) \coloneqq \left(t^+(l)\cap t^-(l); 1\leq l \leq n\right)\text,
	\]
	for the sequences of sets \(\left(t^+(l); 1\leq l \leq n\right) = \bs{t}^+\) and \(\left( t^-(l); 1\leq l \leq n\right) = \bs{t}^-\).
	
	{\bfseries Step I\;\;} The sequence \( \bs{t}^\otimes\) satisfies the following conditions:
	\begin{itemize}
	\item[1')] \(\Delta_1\bs{t}^\otimes =\Delta_1\bs{t}^+\),\quad  \(\Delta_{-1}\bs{t}^\otimes =\Delta_{-1}\bs{t}^-\),
	\item[2')] \(\bigcup\bs{t}^\otimes=t\),\quad \(d_x(\bs{t}^\otimes)\leq 11\gamma\),
	\item[3')] \(\bs{J}\bs{t}^\otimes=\bs{J}\bs{t}^+\),\quad \(\bs{S}\bs{t}^\otimes=\bs{S}\bs{t}^-\),
	\item[4)] \(t^\otimes(l)\in\mathcal{F}_{\beta'}\), for \(1\leq l\leq n\),\; and\; \(t^\otimes(l)\cap t^\otimes(l-1)\in\mathcal{F}_{\beta'}\), for \(2\leq l\leq n\),
	\item[5)] for \(1\leq l<\bar{l}\leq n\) the sets \(t^\otimes(l)\setminus t^\otimes(l+1)\) and \(t^\otimes(\bar{l})\setminus t^\otimes(\bar{l}-1)\) are separated.
	\end{itemize}
	
	{\bfseries Proof of Step I\;\;} The condition 3') results from equalities \(t^+(n) = t\) and \(t^-(1)  = t\), when \(\bs{t}^+\) is a net and \(\bs{t}^-\) is an anti-net for $t$.
	
	Let us show 2'). Directly from the definition of relations \(\prec_{-1,6}^t\) and \(\succ_{-3,-5}^t\), for \(\left(\chi(l); 1\leq l\leq n\right):= \bs{\chi}\), we have 
	\[
	t^+(l)\subset\left(x\leq\chi(l)+6\gamma\right)\text{ and } t^-(l)\subset\left(x\geq\chi(l)-5\gamma\right) \text, \]
	for \(1\leq l\leq n\), thus \(d_x(\bs{t}^\oplus)\leq 11\gamma\). Similarly, we obtain the inclusion
	\[ t^+(l)\cap t^-(l) \supset t\cap\left(\chi(l)-3\gamma\leq x\leq \chi(l)-\gamma\right) \text,\]
	which yields \(\bigcup\bs{t}^\otimes \supset t\cap \left(\chi(1)-3\gamma\leq x\leq \chi(n)-\gamma\right) = t\cap \left(x\leq\chi(n)-\gamma\right)\). By 3') we also have \(t^\otimes(n)=t^-(n)\supset t\cap(x\geq \chi(n)-3\gamma)\), thus finally \(\bigcup\bs{t}^\otimes  = t\).
	
	To obtain 4), we note that each set \(s=t^\otimes(l)\cap t^\otimes(l-1)\), for \(l\in\{2,\ldots,n\}\), satisfies \(s=t^+(l-1)\cap t^-(l)\), which implies that it is closed, convex and \(\beta'/2\)-convex at every point \(\bs{x} \in (\partial s)\cap\operatorname{int}{(\bbt)}\), by analogous properties of the sets \(t^+(l-1),t^-(l)\). To show that \(\operatorname{int}{(s)}\neq \emptyset\), let us consider two cases:
	\begin{itemize}
		\item [a)] \( t\cap \left(x<\chi(l-1)-\gamma\right)\neq \emptyset\),
		\item [b)] \( t\subset \left(x\geq \chi(l-1)-\gamma\right)\).
	\end{itemize}
	
	In the case a), by the fact that \(\chi(l)-\chi(l-1)\leq \gamma\), we have the open set \( U = \left(\chi(l)-3\gamma<x<\chi(l-1)-\gamma\right)\) satisfying \( U\cap t\neq\emptyset \) and obvously \( U\cap t\subset s\). It is thus enough to note that \(\operatorname{int}{(U\cap t)}\neq \emptyset\), in accordance with Remark \ref{rem:elem_tcapU}.
	
	In the case b) we have \(t\subset t\cap\left(x\geq\chi(l)-3\gamma\right)\subset t^-(l)\),
	so \(t^\otimes(l)=t^+(l-1)\) and we use the fact that \(\operatorname{int}({t^+(l-1)})\neq\emptyset\).
	
	We have thus obtained that \(t^\otimes(l)\cap t^\otimes(l-1)\in\mathcal{F}_{\beta'}\), for \(2\leq l\leq n\). Similarly, it can be shown that \(t^\otimes(l)\in\mathcal{F}_{\beta'}\), for \(1\leq l\leq n\).
	
	Now we show 1'). For \(2\leq l\leq n\), let us denote \(A=t^+(l)\), \(B=t^+(l-1)\), \(C=t^-(l)\), \(D=t^-(l-1)\). Then we have \(C\subset D \) and \(A\setminus B\subset t\cap\left(x>\chi(l-1)-\gamma\right)\subset t\cap\left(x>\chi(l)-2\gamma\right)\), \(C\supset t\cap\left(x>\chi(l)-3\gamma\right)\) which implies \(A\setminus B\subset C\). As a result, \(t^+(l)\setminus t^+(l-1)=A\setminus B = (A\cap C)\setminus (B\cap D) = t^\otimes(l)\setminus t^\otimes(l-1)\), by Remark \ref{rem:element_zbiry}. We have thus shown that \(\Delta_1\bs{t}^\otimes =\Delta_1\bs{t}^+\).
	
	The equality \(\Delta_{-1}\bs{t}^\otimes =\Delta_{-1}\bs{t}^-\) can be obtained by a similar argument.
	
	The condition 5), in view of already proven 1'), can be deduced from the relations
	\( t^+(\bar{l})\setminus t^+(\bar{l}-1) \subset \left(x\geq\chi(\bar{l}-1)-\gamma\right)\),
	\( t^-(l)\setminus t^-(l+1) \subset \left(x\leq\chi(l+1)-3\gamma\right)\),
	which implies that \( t^+(\bar{l})\setminus t^+(\bar{l}-1) \) and
	\( t^-(l)\setminus t^-(l+1) \) are separated for \(l<\bar{l}\).
	
	{\bfseries Step II\;\;} For the sequence of sets \(\bs{t}'=\left(t'(k); 1\leq k\leq 2n\right)\) defined by the relation \(t'(2l-1)=t'(2l)=t^\otimes(l)\), for \(1\leq l\leq n\),
	there is a population of souls \(\left\langle\bs{t}',\bs{t}'_1\right\rangle \in \left(\tau_{\beta'\gamma'}\right)^{2n} \) for which conditions 1), 2), 3) are satisfied.
	
	{\bfseries Proof of Step II\;\;} Directly from the definition of \(\bs{t}'\) it follows that there is a sequence \(\bs{t}'_1\) defined by 
	\begin{equation}\label{eq:2}
		\bs{t}_1'|_1 = \Delta_1\bs{t}'\text{,}\qquad \bs{t}_1'|_{-1} = \Delta_{-1}\bs{t}' \text.
	\end{equation}
	By the condition 4) combined with \eqref{eq:2}, there exists a sequence of souls \(\langle\bs{t}',\bs{t}'_1\rangle\) from the set \((\tau_{\beta'})^{2n}\). Moreover, the condition 1') implies that \(\langle\bs{t}',\bs{t}'_1\rangle\in\left(\tau_{\beta'\gamma'}\right)^{2n}\). It also implies that the condition 1) is satisfied. Similarly, 2) follows form 2'), and 3) follows from 3').
	
	By its definition, the sequence \(\langle\bs{t}',\bs{t}'_1\rangle\) is regular, so it is enough to show that it is consistent as well. With $\left(r(k); 1\leq k\leq 2n\right):= \bs{t}'_1|_{-1}$, $\left(s(k); 1\leq k\leq 2n\right):= \bs{t}'_1|_1$, and numbers \(k\), \(\bar{k}\) satisfying \(r(k)\neq\emptyset\), \(s(\bar{k})\neq\emptyset\), \(k<\bar{k}\), we have \(k=2l\), \(\bar{k}=2\bar{l}-1\), for some \(l,\bar{l}\in\left\{1,\ldots,n\right\}\), \(l<\bar{l}\). Then it follows that \(r(k)=t^\otimes(l)\setminus t^\otimes(l+1)\), \(s(\bar{k})=t^\otimes(\bar{l})\setminus t^\otimes(\bar{l}-1)\), thus the sets \(r(k)\), \(s(\bar{k})\) are separated, in accordance with the condition 5). Therefore the sequence \(\langle\bs{t}',\bs{t}'_1\rangle\) is consistent. \qed
	
	\medskip
	
Now we proceed to the proof of Theorem \ref{thm:fund}, which relies on a combination of Theorem \ref{thm:withast}, Lemma \ref{lem:withalpha}, and Corollary \ref{cor:withbeta}, with the use of Lemma \ref{lem:osmy_IstPopDusz}.
We will start with some auxiliary arguments concerning sequences of sets and numbers.
	
For any numbers \(m,\hat{m}\in\mathbb{N}\), \(\hat{m}\geq m\), any sequences \(\bs{a}=(a(k); 1\leq k\leq m)\), \(\hat{\bs{a}}=(\hat{a}(l); 1\leq l\leq \hat{m})\), and systems \((\bs{a}^1,\ldots,\bs{a}^i)\), \((\hat{\bs{a}}^1,\ldots,\hat{\bs{a}}^i)\) of sequences \(\bs{a}^h=\left(a^h(k); 1\leq k\leq m\right)\), \(\hat{\bs{a}}^h=\left(\hat{a}^h(k); 1\leq k\leq \hat{m}\right)\), for \(h\in\{1,\ldots ,i\}\), we will write
\begin{itemize}
	\item \((\hat{\bs{a}}^1,\ldots,\hat{\bs{a}}^i) \vartriangleright_{\hat{m}} (\bs{a}^1,\ldots,\bs{a}^i)\), if there exists a non-decreasing surjective map  \[k\colon\{1,\ldots,\hat{m}\}\to\{1,\ldots,m\}\]satisfying \(\hat{a}^h(l)=a^h(k(l))\), for \(h\in\{1,\ldots,i\}\), \(l\in\{1,\ldots,\hat{m}\}\); the system \((\hat{\bs{a}}^1,\ldots,\hat{\bs{a}}^i)\) will be called a \emph{stretch} of the system \((\bs{a}^1,\ldots,\bs{a}^i)\),
	\item \(\hat{\bs{a}} \vartriangleright_{\hat{m}} \bs{a}\) if \(\hat{a}(l)=a(k(l))\), for \(l\in\{1,\ldots,\hat{m}\}\) and some non-decreasing surjective map \(k\colon\{1,\ldots,\hat{m}\}\to\{1,\ldots,m\}\); the sequence \(\hat{\bs{a}}_1\) will be called a \emph{stretch} of the sequence \(\bs{a}\). Therefore, \(\hat{\bs{a}}\) is obtained from \(\bs{a}\), by possibly repeating some of its elements multiple times.
\end{itemize}

For sequences of sets \(\bs{s}\) and \(\hat{\bs{s}}\) satisfying \(\hat{\bs{s}} \vartriangleright_{\hat{m}}\bs{s}\), and a soul \((t,t_1)\) we thus have
\begin{equation}
\label{eq:complex3}
	\begin{gathered}
		\bs{J}\hat{\bs{s}} = \bs{J}\bs{s}\text,\quad \bs{S}\hat{\bs{s}} = \bs{S}\bs{s}\text,\quad
		\Delta_1\hat{\bs{s}} \ssucc \Delta_1\bs{s}\text,\quad \Delta_{-1}\hat{\bs{s}} \ssucc \Delta_{-1}\bs{s}\text, \\
		\bs{s}\in\mathcal{N}_m^{t\setminus t_1} \implies \hat{\bs{s}}\in\mathcal{N}_{\hat{m}}^{t\setminus t_1}\text,\quad
		\bs{s}\in\mathcal{A}_m^{t} \implies \hat{\bs{s}}\in\mathcal{A}_{\hat{m}}^{t}\text,\\
		\bs{s}\in\mathcal{S}_m^{t_1} \implies \hat{\bs{s}}\in\mathcal{S}_{\hat{m}}^{t_1}\text,
	\end{gathered}
\end{equation}	
Whenever \((\hat{\bs{s}},\hat{\bs{\chi}}) \vartriangleright_{\hat{m}} (\bs{s},\bs{\chi})\) for a skeleton \(\bs{\chi}\), then, \(\hat{\bs{\chi}}\) is also a skeleton and
\begin{equation}\label{eq:complex3.1}
\begin{gathered}
	\bs{\chi}\prec_{ij}^{t\setminus t_1} \bs{s} \implies \hat{\bs{\chi}}\prec_{ij}^{t\setminus t_1} \hat{\bs{s}}\text,\quad\text{for }i,j\in\mathbb{Z}\text,\\
	\bs{s}\succ_{-3,-5}^t\bs{\chi} \implies \hat{\bs{s}}\succ_{-3,-5}^t\hat{\bs{\chi}} \text.
\end{gathered}
\end{equation}
For any sequence of numbers \(\bs{\psi}=(\psi(k); 1\leq k\leq m)\) let us denote
\[
	\|\bs{\psi}\|_\infty=\max{\{\lvert\psi(k)\rvert\colon 1\leq k\leq m\}}
	\text.
\]
Then, for any skeletons \(\bs{\chi}\), \(\bs{\psi}\), and any sequence of sets \(\bs{s}\), of the same length \(m\in\mathbb{N}\), we have
\begin{equation}\label{eq:complex3.2}
	\left({\bs{\chi}\prec_{3,5}^{t\setminus t_1} \bs{s}\quad\text{and}\quad
	\|\bs{\chi}-\bs{\psi}\|_\infty \leq \gamma}\right) \implies \bs{\psi}\prec_{2,6}^{t\setminus t_1} \bs{s}
	\text.
\end{equation}

\begin{namelist}{lll}
	\item{$\bullet$} We will now fix a soul \((t,t_1)\in\tau_{\beta\gamma}\), and for the numbers \(n_1\), \(n^*\) (fixed in Lemma and Theorem \ref{thm:withast}, respectively), we also fix a station \(\bs{t}^*_{t_1}\in\mathcal{S}_{n^*}^{t_1}\), a net \(\bs{t}^+_{t\setminus t_1}\in\mathcal{N}_{n_1}^{t\setminus t_1}\), an anti-net \(\bs{t}^-_t\in\mathcal{A}_{n_1}^t\), and skeletons \(\bs{\chi}^+_{t\setminus t_1}\in\mathcal{R}_{n_1}^{t\setminus t_1}\), \(\bs{\chi}^-_{t}\in\mathcal{R}_{n_1}^{t}\). It will be assumed that \(\bs{t}_t^-\succ_{-3,-5}^t\bs{\chi}_t^-\) and \(\bs{\chi}_{t\setminus t_1}^+\prec_{3,5}^{t\setminus t_1}\bs{t}_{t\setminus t_1}^+\).
\end{namelist}

\begin{lem}\label{lem:existsstretch10}
	There exist stretches
	\[
	\left(\hat{\bs{t}}^-(t,t_1), \hat{\bs{\chi}}^-(t,t_1)\right)
	\vartriangleright_{2{n_1}} \left({\bs{t}}^-_t, {\bs{\chi}}^-_t\right),
	\]
	\[
	\hat{\bs{t}}^+(t,t_1)
	\vartriangleright_{2{n_1}} {\bs{t}}^+_{t \setminus t_1}
	\]
	satisfying \( \hat{\bs{\chi}}^-(t,t_1) \prec_{2,6}^{t\setminus t_1} \hat{\bs{t}}^+(t,t_1)\). 
\end{lem}

Sequences that are subject of the above lemma will from now be considered fixed. For reader's convenience we prove the following elementary fact.
\begin{lem}\label{lem:anotherelemforseq11}
	For skeletons \(\bs{\chi}=(\chi(k); 1\leq k\leq m)\), \(\bs{\psi}=(\psi(k); 1\leq k\leq n)\), satisfying \(\lvert \chi(1) - \psi(1)\rvert ,\lvert \chi(m) - \psi(n)\rvert\leq \gamma \) there exist respective stretches \(\hat{\bs{\chi}}\vartriangleright_{m+n}\bs{\chi}\), \(\hat{\bs{\psi}}\vartriangleright_{m+n}\bs{\psi}\) satisfying \( \|\hat{\bs{\chi}}-\hat{\bs{\psi}}\|_\infty \leq \gamma\).
\end{lem}
\medskip
\noindent
\textbf{Proof.}
	The lemma is valid for \(m+n=2\), i.e. \(m=n=1\). Indeed, for \(\bs{\chi}=\left(\chi(1)\right)\) and \(\bs{\psi}=\left(\psi(1)\right)\) we can take \(\hat{\bs{\chi}}=\left(\chi(1),\chi(1)\right)\), \(\hat{\bs{\psi}}=\left(\psi(1),\psi(1)\right)\). Let us assume that the statement of the Lemma is true for a fixed value of \(m+n=: p\geq 2\). We will show that it is also true when \(m+n=p+1\). Without loss of generality we can assume that \(\chi(m)\geq \psi(n)\). If \(m=1\), then relation \(\|\hat{\bs{\chi}}-\hat{\bs{\psi}}\|_\infty \leq \gamma\) can be obtained for stretches
	\[
	\hat{\bs{\chi}}=\underbrace{\left(\chi(1),\ldots,\chi(1)\right)}_{n+1}
	\quad
	\text{and}\quad\hat{\bs{\psi}}=\left(\psi(1),\ldots,\psi(n),\psi(n)\right)
	\text.
	\]
	
	Let \(m\geq 2\). Then, we have \(\lvert \chi(m-1)-\psi(n)\rvert \leq\gamma \) and, as assumed, there exists stretches \(\bs{\chi}^1\vartriangleright_{m+n-1}\left(\chi(1),\ldots,\chi(m-1)\right)\) and \(\bs{\psi}^1\vartriangleright_{m+n-1}\bs{\psi}\) satisfying \(\|\bs{\chi}^1-\bs{\psi}^1\|_\infty \leq \gamma\). Thus, it is enough to consider concatenations \(\hat{\bs{\chi}}=\bs{\chi}^1\circ\left(\chi(m)\right)\) and \(\hat{\bs{\psi}}=\bs{\psi}^1\circ\left(\psi(m)\right)\), in which single element is appended. \qed

\medskip
\noindent
\textbf{Proof of \textup{Lemma \ref{lem:existsstretch10}}}
	By assumption $d(t_1) \leq \gamma$, we have $|\chi_t^-(1) - \chi_{t \setminus t_1}^+(1)| \leq \gamma$, $|\chi_t^-(n_1) - \chi_{t \setminus t_1}^+(n_1)| \leq \gamma$ for $(\chi_t^-(l); 1 \leq l \leq n_1) = \bs{\chi}_t^-$, $(\chi_{t \setminus t_1}^+(l); 1 \leq l \leq n_1) = \bs{\chi}_{t \setminus t_1}^+$. By Lemma \ref{lem:anotherelemforseq11}, for \(\bs{\chi}=\bs{\chi}^+_{t\setminus t_1}\), \(\bs{\psi}=\bs{\chi}^-_t\), \(m=n=n_1\), we can find stretches
	\[
	\left(\hat{\bs{t}}^-(t,t_1), \hat{\bs{\chi}}^-(t,t_1)\right)
	\vartriangleright_{2{n_1}} \left({\bs{t}}^-_t, {\bs{\chi}}^-_t\right),
	\]
	\[
	(\hat{\bs{t}}^+(t,t_1), \hat{\bs{\chi}}^+(t,t_1))
	\vartriangleright_{2{n_1}} ({\bs{t}}^+_{t \setminus t_1}, {\bs{\chi}}^+_{t \setminus t_1})
	\]
	satisfying \( \left\|\hat{\bs{\chi}}^+(t,t_1)- \hat{\bs{\chi}}^-(t,t_1) \right\|_\infty\leq \gamma\). Then, \(\hat{\bs{\chi}}^+(t,t_1) \prec_{3,5}^{t\setminus t_1} \hat{\bs{t}}^+(t,t_1)\), cf. \eqref{eq:complex3.1}, and \(\hat{\bs{\chi}}^-(t,t_1) \prec_{2,6}^{t\setminus t_1} \hat{\bs{t}}^+(t,t_1)\), cf. \eqref{eq:complex3.2}. \qed

Recall that for sequences of sets \(\bs{r} = \left(r(l); 1\leq l\leq m\right)\), \(\bs{s} = \left(s(l); 1\leq l\leq m\right)\) we write \(\bs{r} \perp \bs{s}\) if for any \(l\in\{1,\ldots,m\} \) we have \(r(l)=\emptyset\) or \(s(l)=\emptyset\).

An important idea that will be used to construct offspring of souls is contained in the following lemma, for \(\bs{t}^{+}(t,t_1)\), \(\bs{t}^-(t,t_1)\), \(\bs{\chi}^-(t,t_1)\) appearing in Lemma \ref{lem:existsstretch10}.

\begin{lem}
\label{lem:existstretches12}
	There exist stretches
	\[
	\left(
	\doublehat{\bs{t}}^{+}(t,t_1), \doublehat{\bs{t}}^{-}(t,t_1), \doublehat{\bs{\chi}}^{-}(t,t_1)\right)
	\vartriangleright_{2{n_1}+n^*} \left(\hat{t}^{+}(t,t_1), \hat{\bs{t}}^{-}(t,t_1), \hat{\bs{\chi}}^-(t,t_1)\right)
	\]
	and \( \hat{\bs{t}}^*(t,t_1) \vartriangleright_{2{n_1}+n^*} \bs{t}^*_{t_1}\) satisfying
	\begin{itemize}
		\item [\emph{1)}] \(\Delta_1\doublehat{\bs{t}}^{+}(t,t_1) \perp \Delta_1\hat{\bs{t}}^*(t,t_1)\),
		\item [\emph{2)}] \(\doublehat{\bs{\chi}}^{-}(t,t_1) \prec_{-1,1}^{t_1} \hat{\bs{t}}^*(t,t_1)\).
	\end{itemize}
\end{lem}

The sequences \(\doublehat{\bs{t}}^{+}(t,t_1)\), \(\doublehat{\bs{t}}^{-}(t,t_1)\), \(\doublehat{\bs{\chi}}^{-}(t,t_1)\) and \(\hat{\bs{t}}^*(t,t_1)\) will be considered fixed.

\medskip
\noindent
\textbf{Proof of \textup{Lemma \ref{lem:existstretches12}}}\quad
	Since the soul \((t,t_1)\) is fixed, we will write for short
	\begin{eqnarray*}
		\left(t^*(l); 1\leq l\leq n^*\right)&\coloneqq& \bs{t}^*_{t_1}\text, \\
		\left(\hat{t}^+(l); 1\leq l\leq 2{n_1}\right)&\coloneqq& \hat{\bs{t}}^+(t,t_1) \text, \\
		\left(\hat{t}^-(l); 1\leq l\leq 2{n_1}\right)&\coloneqq& \hat{\bs{t}}^-(t,t_1)\text,\\
		\left(\hat{\chi}^-(l); 1\leq l\leq 2{n_1}\right)&\coloneqq& \hat{\bs{\chi}}^-(t,t_1)\text.
	\end{eqnarray*}
For
	\[ l^0 \coloneqq \min{\left\{l\in\{1,\ldots,2{n_1}\}\colon t_1\subset \left(x\leq\hat{\chi}^-(l)+\gamma\right)\right\}}
	\text,
	\]
	in view of the inequality \(d(t_1)\leq\gamma\), we have
	\begin{equation}\label{eq:e4}
		t_1\subset \left(\hat{\chi}^-(l^0)-\gamma\leq x\leq\hat{\chi}^-(l^0)+\gamma\right)\text.
	\end{equation}
	Let us set
	\begin{align*}
		\hat{\bs{t}}^*(t,t_1) &= {\underbrace{(\emptyset,\ldots,\emptyset)}_{l^0}}
			\circ{\bs{t}^*_{t_1}} \circ{\underbrace{(t_1,\ldots,t_1)}_{2{n_1}-l^0}}
			\text, \\
		\doublehat{\bs{t}}^+(t,t_1) &= {\left(\hat{t}^+(1),\ldots,\hat{t}^+(l^0)\right)} \circ
			{\underbrace{\left(\hat{t}^+(l^0),\ldots,\hat{t}^+(l^0)\right)}_{n^*}} \circ
			\left(\hat{t}^+(l^0+1),\ldots,\hat{t}^+(2{n_1})\right)
			\text, \\
		\doublehat{\bs{t}}^-(t,t_1) &= {\left(\hat{t}^-(1),\ldots,\hat{t}^-(l^0)\right)} \circ
			{\underbrace{\left(\hat{t}^-(l^0),\ldots,\hat{t}^-(l^0)\right)}_{n^*}} \circ
			\left(\hat{t}^-(l^0+1),\ldots,\hat{t}^-(2{n_1})\right)
			\text, \\
		\doublehat{\bs{\chi}}^-(t,t_1) &= {\left(\hat{\chi}^-(1),\ldots,\hat{\chi}^-(l^0)\right)} \circ
			{\underbrace{\left(\hat{\chi}^-(l^0),\ldots,\hat{\chi}^-(l^0)\right)}_{n^*}} \circ
			\left(\hat{\chi}^-(l^0+1),\ldots,\hat{\chi}^-(2{n_1})\right)
			\text.
	\end{align*}
	 In order to show that condition 2) is satisfied let us denote \(\hat{n}=2{n_1}+n^*\) and
\begin{eqnarray*}
\left(\hat{t}^*(k); 1\leq k\leq\hat{n}\right) &=& \hat{\bs{t}}^*(t,t_1)\text,\\
\left(\doublehat{\chi}^-(k); 1\leq k\leq\hat{n}\right) &=& \doublehat{\bs{\chi}}^-(t,t_1)\text.
\end{eqnarray*}
	 By the relation~\eqref{eq:e4}, for \(k\in\{1,\ldots,\hat{n}\}\), we obtain the implications
	 \begin{multline*}
	 \hat{t}^*(k) \neq\emptyset \implies k>l^0 \implies \doublehat{\chi}^-(k)+\gamma\geq\hat{\chi}^-(l^0)+\gamma \\ 
	 \implies
	 t_1\cap\left(x\leq\doublehat{\chi}^-(k)+\gamma\right)\supset t_1\supset \hat{t}^*(k)
	 \text,
	 \end{multline*}
	 thus \( \hat{t}^*(k) \subset t_1\cap\left(x\leq\doublehat{\chi}^-(k)+\gamma\right) \), for any \(k\in\{1,\ldots,\hat{n}\}\).
	 We also obtain the implications
	 \[
	 t_1\cap\left(x\leq\doublehat{\chi}^-(k)-\gamma\right) \neq \emptyset\implies
	 k>l^0+n^* \implies \hat{t}^*(k)=t_1 \text,
	 \]
	which leads to \(t_1\cap\left(x\leq\doublehat{\chi}^-(k)-\gamma\right)\subset \hat{t}^*(k)\), and proves statement 2).
	
	Relation 1) is a direct consequence of definitions of the sequences \(\hat{\bs{t}}^*(t,t_1)\) and \(\doublehat{\bs{t}}^+(t,t_1)\). \qed

As a means of preparation for the proof of Theorem \ref{thm:fund}, we will further demonstrate certain general properties of nets and stations for the fixed soul \((t,t_1)\).
The number \(\hat{n}=2{n_1}+n^*\) will also be fixed.

\begin{lem}
\label{lem:l13}
	If \(\bs{r}\in\mathcal{S}^{t_1}_{\hat{n}}\), \(\bs{s}\in\mathcal{N}^{t\setminus t_1}_{\hat{n}}\), \(\bs{\chi}\in\mathcal{R}_{\hat{n}}\), for which \(\bs{\chi}\prec_{-1,1}^{t_1}\bs{r}\), \(\bs{\chi}\prec_{2,6}^{t\setminus t_1}\bs{s}\), and \(\Delta_1 \bs{r} \perp \Delta_1 \bs{s}\), then we have the following relations:
	\begin{itemize}
		\item[\emph{1)}] \(\bs{r}\oplus\bs{s}\in\mathcal{N}^t_{\hat{n}}\),
		\item[\emph{2)}] \( \bs{\chi} \prec_{-1,6}^t {(\bs{r}\oplus\bs{s})}\),
		\item[\emph{3)}] \( {t_1\otimes\Delta_1(\bs{r}\oplus\bs{s})} \ll
							{(t\setminus t_1)\otimes \Delta_1((\bs{r}\oplus\bs{s}))}\)
		\item[\emph{4)}] \( {\left(t_1\otimes\Delta_1(\bs{r}\oplus\bs{s})\right)\perp\left((t\setminus t_1)\otimes\Delta_1(\bs{r}\oplus\bs{s})\right)}\).
	\end{itemize}
\end{lem}

The task of showing that \((\bs{r}\oplus\bs{s})\in\left(\mathcal{F}_{\beta'}\right)^{\hat{n}}\), that is property 1) in Lemma \ref{lem:l13} turns out to be somewhat tedious. In the proof we will make use of the following observation.

\begin{lem}\label{lem:l14}
	If for some sets \(s_0\subset s_1\subset s_2\) we have \(s_0,s_2\in\mathcal{F}_{\beta'}\) and the sets \(s_1\setminus s_0\), \(s_2\setminus s_1\) are separated, then \(s_1\in \mathcal{F}_{\beta'}\).
\end{lem}
\medskip
\noindent
\textbf{Proof.}
	First let us note that \(\operatorname{int}{(s_1)}\neq\emptyset\), since \(s_0\in\mathcal{F}_{\beta'}\).
	
	Let \(\bs{x}\in(\partial s_1)\cap\operatorname{int}{(\bbt)}\). Denoting \(\Delta_1=s_1\setminus s_0\), \(\Delta_2=s_2\setminus s_1\), for some neighborhood \(V\) of \(\bs{x}\) we have the alternatives: 1) \(V\cap\Delta_1=\emptyset\) or 2) \(V\cap\Delta_2=\emptyset\). If statement 1) is true, then \(V\cap s_1 = V\cap s_0\). If statement 2) is true, then \(V\cap s_1=V\cap s_2\). In any case, the set \(s_1\) is \(\beta'\)-convex at point \(\bs{x}\). Similarly one can show the implication \(\bs{x}\in\partial s_1\implies \bs{x}\in s_1\), which is equivalent to \(s_1\) being a closed set.
	
	To obtain convexity of \(s_1\), let us consider arbitrary \(\bs{x}\in s_0\), \(\bs{y}\in\Delta_1\). There is a point \(\bs{z}\in\left[\bs{x},\bs{y}\right]\cap\partial\Delta_1\), satisfying \(\left[\bs{z},\bs{y}\right]\setminus \bs{z}\subset\Delta_1\). Then, \(\bs{z}\in s_2\), \(\bs{z}\notin\Delta_2\), which implies \(\bs{z}\in s_0\), \(\left[\bs{x},\bs{z}\right]\subset s_0\), and, as a result, we have \(\left[\bs{x},\bs{y}\right]\subset s_1\). Next, let us consider \(\bs{y},\bar{\bs{y}}\in\Delta_1\), such that $[\bs{y}, \bar{\bs{y}}] \not\subset \Delta_1$, then $[\bs{y}, \bar{\bs{y}}] \subset s_2$ and $[\bs{y}, \bar{\bs{y}}] \not\subset \Delta_1 \cup \Delta_2$. Then there is a point \(\bs{x}\in\left[\bs{y},\bar{\bs{y}}\right]\cap s_0\). As previously, it follows that \(\left[\bs{x},\bs{y}\right],\left[\bs{x},\bar{\bs{y}}\right]\subset s_1\). This finally shows that the set \(s_1\) is convex. \qed

\medskip
\noindent
\textbf{Proof of \textup{Lemma \ref{lem:l13}}}
For sequences \(\left(r(l); 1\leq l\leq \hat{n}\right)=\bs{r}\), \(\left(s(l); 1\leq l\leq \hat{n}\right)=\bs{s}\), and \(\left(\chi(l); 1\leq l\leq \hat{n}\right)=\bs{\chi}\) we have
\[
t\cap\left(x\leq\chi(l)-\gamma\right)\subset t\cap\left(r(l)\cup s(l)\right)\subset t\cap \left(x\leq\chi(l)+6\gamma\right)
\text,
\]
that is relation 2) holds. Moreover, for \(1\leq k<l\leq\hat{n}\), with \(r(0)\coloneqq r(1)\), we have
\(r(k)\setminus r(k-1)\subset\left(x\leq\chi(k)+\gamma\right)\) and \(s(l)\sm s(l-1)\subset \left(x>\chi(l-1)+2\gamma\right)\subset \left(x>\chi(k)+2\gamma\right)\), therefore the sets \(r(k)\sm r(k-1)\) and \(s(l)\sm s(l-1)\) are separated. That is, relation 3) holds.

Relation 4) results directly from definitions of sequences \(\bs{r}\) and \(\bs{s}\) and the fact that \(\Delta_1\bs{r}\perp\Delta_1\bs{s}\).

We also have
\( d\left(\Delta_1(\bs{r}\oplus\bs{s})\right) = \max{ \left\{  d\left(\Delta_1\bs{r}\right), d\left(\Delta_1\bs{s}\right)\right\} } \leq \gamma\), by which 1) will follow from the condition:
\begin{equation}\label{eq:e5}
	r(l)\cup s(l)\in\mathcal{F}_{\beta'}\text{ for }1\leq l\leq \hat{n} \text.
\end{equation}
Notice that for any fixed \(l\in\{1,\ldots,\hat{n}\}\) we have \(r(l)\cup (t\sm t_1)\in\mathcal{F}_{\beta'}\) and \(s(l)\in \mathcal{F}_{\beta'}\), by definition of a station and a net respectively. Denoting \(s_0=s(l)\), \(s_1=r(l)\cup s(l)\), \(s_2=r(l)\cup (t\sm t_1)\), \(\Delta_1=s_1\sm s_0=r(l)\), \(\Delta_2=s_2\sm s_1 = (t\sm t_1)\sm s(l)\), we have in turn \(\Delta_1\subset\left(x\leq \chi(l)+\gamma\right)\), \(\Delta_2\subset\left(x\geq\chi(l)+3\gamma\right)\). Therefore \(\Delta_1\) and \(\Delta_2\) are separated and \(s_1\in\mathcal{F}_{\beta'}\), by Lemma \ref{lem:l14}. That is, condition \eqref{eq:e5} is satisfied. \qed

\medskip

Our argumentation concerning the arbitrarily fixed soul \((t,t_1)\), in view of Lemma \ref{lem:l13}, lead to the following corollary, for $\hat{n} = 2n_1 + n*$.
\begin{cor}\label{cor:c15}
	{\bfseries\textup{A.}} Let us set \(\bs{t}^{\#}(t,t_1) = \hat{\bs{t}}^{*}(t,t_1) \oplus \doublehat{\bs{t}}^+(t,t_1)\). Then we have:
	\begin{itemize}
		\item[\textup{1)}] \(\bs{t}^{\#}(t,t_1) \in\mathcal{N}^t_{\hat{n}}\),
		\item[\textup{2)}] \(\doublehat{\bs{\chi}}^{-}(t,t_1) \prec_{-1,6}^t \bs{t}^{\#}(t,t_1) \),
		\item[\textup{3)}] \( \left(t_1\otimes\Delta_1\bs{t}^{\#}(t,t_1) \right) \perp \left( (t\sm t_1)\otimes \bs{t}^{\#}(t,t_1) \right) \),\\
		\(t_1\otimes \Delta_1\bs{t}^{\#}(t,t_1) \ll (t\sm t_1)\otimes \Delta_1\bs{t}^{\#}(t,t_1)\),
		\item[\textup{4)}] \(\bs{J}\bs{t}^{\#}(t,t_1)=\bs{J}\bs{t}^+_{t\sm t_1}\),
		\item[\textup{5)}] \(t_1\otimes\Delta_1\bs{t}^{\#}(t,t_1)\ssucc\Delta_1\bs{t}^*_{t_1}\),
		\item[\textup{6)}] \((t\sm t_1)\otimes\Delta_1\bs{t}^{\#}(t,t_1)\ssucc\Delta_1\bs{t}^+_{t\sm t_1}\).
	\end{itemize}
	
	{\bfseries \textup{B.}} For sequence \(\doublehat{\bs{t}}^-(t,t_1)\) we have
	\begin{itemize}
		\item[\textup{7)}] \(\doublehat{\bs{t}}^-(t,t_1) \in \mathcal{A}^t_{\hat{n}}\),
		\item[\textup{8)}] \(\doublehat{\bs{t}}^-(t,t_1) \succ_{-3,-5}^t \doublehat{\bs{\chi}}^-(t,t_1)\),
		\item[\textup{9)}] \(\bs{S}\doublehat{\bs{t}}^-(t,t_1) = \bs{S}\bs{t}^-_t\),
		\item[\textup{10)}] \(\Delta_{-1}\doublehat{\bs{t}}^-(t,t_1) \ssucc \Delta_{-1}\bs{t}^-_t\).
	\end{itemize}
\end{cor}

\medskip
\noindent
\textbf{Proof.}
	Part B.\ of the above statement is an immediate consequence of the relation \(\left( \doublehat{\bs{t}}^-(t,t_1), \doublehat{\bs{\chi}}^-(t,t_1) \right) \vartriangleright_{\hat{n}} \left(\bs{t}^-_t,\bs{\chi}^-_t\right) \), which in turn can be obtained from Lemmas \ref{lem:existsstretch10} and \ref{lem:existstretches12}, cf. properties of stretches \eqref{eq:complex3}, \eqref{eq:complex3.1}, since, for arbitrary system of sequences, a stretch of a stretch is itself a stretch of the original system.
	
	In order to prove part {A.} of the corollary, let us denote \(\bs{r} = \hat{\bs{t}}^*(t,t_1)\), \(\bs{s}=\doublehat{\bs{t}}^+(t,t_1)\), \(\bs{\chi} = \doublehat{\bs{\chi}}^-(t,t_1)\). By statements 1) and 2) of Lemma \ref{lem:existstretches12}, we have \(\Delta_1\bs{r}\perp\Delta_1\bs{s}\), \(\bs{\chi}\prec_{-1,1}^{t_1}\bs{r}\). In view of the relation
	\(
		\left(\doublehat{\bs{t}}^-(t,t_1), \doublehat{\bs{\chi}}^-(t,t_1)\right)
	\vartriangleright_{\hat{n}}
	\left(\bs{t}^-_t,\bs{\chi}^-_t\right)
	\)
	(following Lemma \ref{lem:existstretches12} and Lemma \ref{lem:existsstretch10}), we also have \(\bs{\chi}\prec_{2,6}^{t\sm t_1}\bs{s}\). Since \(\bs{t}^*_{t_1}\in\mathcal{S}_{n^*}^{t_1}\), \(\bs{t}^+_{t\sm t_1}\in\mathcal{N}_{n}^{t\sm t_1}\), it follows that \(\bs{r}\in\mathcal{S}_{\hat{n}}^{t_1}\) and \(\bs{s}\in\mathcal{N}_{\hat{n}}^{t\sm t_1}\). Thus, Lemma \ref{lem:l13} can be used.
	
	The conditions 1), 2) and 3) with 4) in Lemma \ref{lem:l13} imply that, respectively, conditions 1), 2) and 3) of part {A.}\ are satisfied.
	
	Conditions 5) and 6) will follow from the definitions of the sequence \(\bs{t}^{\#}(t,t_1)\) and the relations \(\hat{\bs{t}}^*(t,t_1) \vartriangleright_{\hat{n}}\bs{t}^*_{t_1}\), \(\doublehat{\bs{t}}^+(t,t_1) \vartriangleright_{\hat{n}}\bs{t}^+_{t\sm t_1}\). Similarly, we can obtain condition 4), considering that \(\bs{J}\hat{\bs{t}}^*=\emptyset\), by the definition of a station. \qed

\medskip
\noindent
\textbf{Proof of \textup{Theorem \ref{thm:fund}}}
	Based on Theorem \ref{thm:withast}, Lemma \ref{lem:withalpha} and Corollary \ref{cor:withbeta}, for our numbers \(n_1,n^*\in\mathbb{N}\) (given in Lemma \ref{lem:withalpha} and Theorem \ref{thm:withast}, respectively, and independent of \(\bbt\)),
	and an arbitrary soul \((t,t_1)\in\tau_{\beta,\gamma}\) we have, in particular, defined sequences
	\begin{equation}\label{eq:e6dep}
		\begin{aligned}
			&\bs{t}_{t_1}^*\in\mathcal{S}^{t_1}_{n^*}
			\text,\quad
			\bs{t}_{t\sm t_1}^+\in\mathcal{N}^{t\sm t_1}_{n_1}\text,\quad
			\bs{t}^-_t\in\mathcal{A}^t_{n_1}
			\text,\\
			&\text{where }
			\bs{t}_{t_1}^*\text{, }
			\bs{t}_{t\sm t_1}^+\text{, and }
			\bs{t}^-_t\text{ depend on }t_1\text{, }t\sm t_1\text{, and }t\text{ respectively.}
		\end{aligned}
	\end{equation}
	
	Next, we have obtained sequences \(\bs{t}^{\#}(t,t_1)\), \(\doublehat{\bs{t}}^-(t,t_1)\), and \(\doublehat{\bs{\chi}}^-(t,t_1)\), described in Corollary \ref{cor:c15}. According to conditions 1), 2), 7), and 8) therein, we can apply Lemma \ref{lem:osmy_IstPopDusz} for \(\bs{t}^+=\bs{t}^{\#}(t,t_1)\), \(\bs{t}^-=\doublehat{\bs{t}}^-(t,t_1)\), \(\bs{\chi}^-=\doublehat{\bs{\chi}}^-(t,t_1)\), and \(\hat{n}=2n_1+n^*\) in place of \(n\), to obtain a population of souls, which we denote by \(\left\langle \bs{t}'(t,t_1), \bs{t}'_1(t,t_1) \right\rangle\).
	
	The proof of Theorem \ref{thm:fund} reduces to the following sequences of observations.
	
	{\bfseries Step I\;\;} Values obtained by the function \(\tau_{\beta,\gamma}\ni(t,t_1) \mapsto \left\langle \bs{t}'(t,t_1), \bs{t}'_1(t,t_1) \right\rangle\) are populations in \(\left(\tau_{\beta'\gamma'}\right)^{m'}\) for \(m'=2\hat{n}\). This follows directly from Lemma~\ref{lem:osmy_IstPopDusz}.
	
	{\bfseries Step II\;\;} Function \(\tau_{\beta,\gamma}\ni(t,t_1) \mapsto  \bs{t}'(t,t_1)\) is a filling, satisfying \(d_x(\bs{t}'(t,t_1))\leq 11\gamma\), for \((t,t_1)\in\tau_{\beta\gamma}\). This is implied by statements 2) and 3) in Lemma~\ref{lem:osmy_IstPopDusz}, together with 4) and 9) in Corollary~\ref{cor:c15}.

	{\bfseries Step III\;\;} The function \(\tau_{\beta,\gamma}\ni(t,t_1) \mapsto \bs{t}'_1(t,t_1)|_1 \) is a dust. This follows from the first relation in condition 1) of Lemma~\ref{lem:osmy_IstPopDusz} and statements 3), 5) and 6) in Corollary~\ref{cor:c15}, combined with the dependence reminded of in \eqref{eq:e6dep}.
	
	{\bfseries Step IV\;\;} The function \(\tau_{\beta,\gamma}\ni(t,t_1) \mapsto \bs{t}'_1(t,t_1)|_{-1} \) is an anti-dust. This is implied by the second relation in condition 1) of Lemma~\ref{lem:osmy_IstPopDusz} and statement 10) of Corollary~\ref{cor:c15}, again, combined with \eqref{eq:e6dep}. \qed

\subsubsection*{B. Construction of stations with the use of nets}

In order to obtain offspring of souls stipulated in Theorem \ref{thm:fund}, it is left for us to prove  Theorem \ref{thm:withast}. We will begin with another reduction of the problem.

Let us recall that the numbers \(\beta\geq 1\); \(\gamma,\gamma'\in(0,1/4)\), and thus \(\beta'=\max{\{2\beta,\gamma + 1/\gamma\}}\) are fixed, as well as the set \(\bbt\) (satisfying \(\mathcal{T}\ni \bbt\subset [-1,1]\times[-1,1]\), \(d(\bbt)\leq 1\)).

For any soul \((t,t_1)\in\tau_{\beta'/2}\) and number \(n\in\mathbb{N}\) we denote
\begin{namelist}{lll}
	\item{$\bullet$} \(
\mathcal{N}^{(t,t_1)}_n = \left\{\left(t(l); 1\leq l\leq n\right)\in\mathcal{N}_n^t\colon t(1)=t\sm t_1\right\}
\text,
\)
\end{namelist}
for the family of nets \(\mathcal{N}^t_n\) given in Definition \ref{def:3.3.3}.

\begin{thm}\label{thm:astast16}
	Assume that for some \(n^*\) independent of \(\bbt\) the following condition is satisfied.
	\begin{itemize}
		\item[{\normalfont(\textasteriskcentered\textasteriskcentered)}] For any soul \( (t,t_1) \in\tau_{\beta'/2}\) there exists a net \(\bs{t}\in\mathcal{N}_{n^*}^{(t,t_1)}\).
	\end{itemize}
	Then, condition \textup{(\textasteriskcentered)} in \emph{Theorem \ref{thm:withast}} is satisfied as well.
\end{thm}

The above theorem is a straightforward consequence of the following.

\begin{thm}\label{thm:t17}
	For \(n\in\mathbb{N}\), \((t,t_1)\in\tau_{\beta'/2}\), and net $\left(t(l); 1\leq l\leq n\right)\in\mathcal{N}_n^{(t,t_1)}$, sequence \(\left(t(l)\sm t(1); 1\leq l\leq n\right)\) is a station in \(\mathcal{S}_n^{t_1}\).
\end{thm}

Our proof of Theorem \ref{thm:t17} will be based on an interesting fact, expressed in the theorem below. A sequence of sets \(\bs{\Delta}=\left(r_0,r_1,r_2,s_0,s_1,s_2\right)\) will be called a \emph{pack} if 
\(r_0,r_1,r_2,s_0,s_2\in\mathcal{F}_{\beta'}\) and
\begin{equation*}
	\begin{gathered}
		r_0\subset r_1\subset r_2\text,\qquad s_0\subset s_1\subset s_2\text, \\
		r_1\sm r_0=s_1\sm s_0\text,\qquad r_2\sm r_0=s_2\sm s_0\text.
	\end{gathered}
\end{equation*}
Pack \(\bs{\Delta}\) and the differences \(\Delta_1=r_1\sm r_0\) and \(\Delta_2=r_2\sm r_0\) will be fixed throughout the proof of the following.
\begin{thm}\label{thm:t18}
	For pack \(\bs{\Delta}\) we have \(s_1\in\mathcal{F}_{\beta'}\).
\end{thm}

A proof of the above theorem will follow from the next three lemmas. For any set \(\Delta\subset\mathbb{R}^2\), we will use the following two components of its boundary:
\begin{namelist}{lll}
	\item{$\bullet$} \quad \(
\partial^+\Delta=(\partial\Delta)\cap\Delta\) and \(\partial^-\Delta=(\partial\Delta)\sm\Delta
\).
\end{namelist}

\begin{lem}\label{lem:l19}
	We have the following.
	\begin{itemize}
		\item[\textbf{A.}] If \(p\subset q\subset\mathbb{R}^2\), \(p=\clo{(p)}\),  \(\Delta=q\sm p\), then
		\[
		q =\operatorname{clos}{(q)} \iff \partial^-\Delta = p\cap\operatorname{clos}{(\Delta)}
		\text.
		\]
		\item[\textbf{B.}] For pack \(\bs{\Delta}\), we have \(\partial^-\Delta_1\subset\partial^-\Delta_2\).
		\item[\textbf{C.}] For pack \(\bs{\Delta}\), the set \(s_1\) is convex and closed.
	\end{itemize}
\end{lem}

\medskip
\noindent
\textbf{Proof.}
	To prove A. it is enough to notice that \(\partial^-\Delta=\left(p\cap\clo(\Delta)\right)\cup\partial^-q\).
	
	Statement B. follows from A., as \(\partial^-\Delta_1=r_0\cap\clo(\Delta_1)\subset r_0\cap\clo(\Delta_2)=\partial^-(\Delta_2)\).
	
	For showing the convexity of \(s_1\), let us take arbitrary \(\bs{x},\bs{y}\in s_1\). If \(\bs{x},\bs{y}\in s_0\), then \([\bs{x},\bs{y}]\subset s_1\), by assumed convexity of \(s_0\).
	
	In the case when \(\bs{x},\bs{y}\in\Delta_1\), by convexity of \(r_1\), we have \([\bs{x},\bs{y}]\subset r_1\), thus \([\bs{x},\bs{y}]\cap(\Delta_2\sm \Delta_1)=\emptyset\). While, by convexity of \(s_2\), it follows that \([\bs{x},\bs{y}]\subset s_2 = s_1\cup(\Delta_2\sm\Delta_1)\), thus \([\bs{x,\bs{y}}]\subset s_1\).
	
	If \(\bs{x}\in s_0\), \(\bs{y}\in\Delta_1\), then, as \(s_2\) is convex, we have \([\bs{x},\bs{y}]\subset s_2\), and, as \(s_0\) is closed and convex, there exists a point \(\bs{z}\in[\bs{x},\bs{y}]\) satisfying \([\bs{x},\bs{z}]\subset s_0\), \([\bs{z},\bs{y}]\sm\{\bs{z}\}\subset \Delta_2\).
	Thus, \(\bs{z}\in\partial^-\Delta_2\), and then \(\bs{z}\in r_0\), in accordance with statement A., with \(p=r_0\), \(q=r_2\).
	By convexity of \(r_1\), we have \([\bs{z},\bs{y}]\sm\{\bs{z}\}\subset \Delta_1\), and finally \([\bs{x},\bs{y}]\subset s_1\).
	
	To see that \(s_1\) is closed, note that A. used with \(p=s_0\), \(q=s_2\) yields \(\partial^-\Delta_2\subset s_0\), thus \(\partial^-\Delta_1\subset s_0\), by statement B.
	Lastly, with \(p=s_0\), \(q=s_1\), A. implies that \(s_1\) is closed. \qed

\begin{lem}\label{lem:l20}
	If  \(\bs{x}\in(\partial s_0)\cap\operatorname{int}(\bbt)\), \(\varepsilon>0\), then
	\[
	\bs{x}\in\operatorname{int}{\left(\operatorname{conv}\left(B(\bs{x},\varepsilon)\sm s_0\right)\right)}
	\text{.}
	\]
\end{lem}
\medskip
\noindent
\textbf{Proof.}
	There is a ball \(B\in\mathcal{B}_{\beta'}\), satisfying \(\bs{x}\in\partial B\), \(s_0\cap B(\bs{x},\varepsilon_1)\subset B\), for some \(\eps_1\in(0,\eps)\). Then, we have
	\begin{multline*}
	\bs{x}\in\operatorname{int}{\left(\operatorname{conv}\left(B(\bs{x},\varepsilon_1)\sm B\right)\right)} \\
	\subset \operatorname{int}{
		\left[\operatorname{conv}{\left(B(\bs{x},\varepsilon_1)\sm (s_0\cap B(\bs{x},\varepsilon_1))\right)}\right]
		}
	\subset \operatorname{int}{\left(\operatorname{conv}\left(B(\bs{x},\varepsilon)\sm s_0\right)\right)}
	\text.
	\end{multline*} \qed

\begin{lem}\label{lem:l21}
	For the pack \(\bs{\Delta}\) and any point  \(\bs{x}\in({\partial}{s_1})\cap\operatorname{int}{(\bbt)}\), the set \(s_1\) is \(\beta'\)-convex at \(\bs x\).
\end{lem}
\medskip
\noindent
\textbf{Proof.}
	In the case when \(\bs{x}\notin\clo{(\Delta_1)}\), the sets $s_1$, \(s_0\) are identical in some neighborhood of \(\bs{x}\), thus, we can use the assumption that \(s_0\in\mathcal{F}_{\beta'}\).
	
	If \(\bs{x}\in\partial^+\Delta_1\), then \(\bs{x}\notin r_0\), \(\bs{x}\notin s_0\), and some neighborhood \(V\) of \(\bs{x}\) satisfies \(V\cap r_0=\emptyset\) and \(V\cap s_0=\emptyset\).
	Thus, we have \(s_1\cap V=\Delta_1\cap V=r_1\cap V\), and we can use the assumption that \(r_1\in\mathcal{F}_{\beta'}\).
	
	Assume that \(\bs{x}\in\partial^-\Delta_1\). By statement B. of Lemma \ref{lem:l19}, we have \(\bs{x}\in\partial^-\Delta_2\) and proceed as follows.
	
	If \(\bs{x}\in\partial s_2\), then \(s_2\) is \(\beta'\)-convex at \(\bs{x}\), as well as is the subset \(s_1\subset s_2\).
	
	If \(\bs{x}\in\operatorname{int}{(s_2)}\), then some neighborhood of \(\bs{x}\) is a subset of \(s_0\cup\Delta_2\). By Lemma \ref{lem:l20}, since \(s_0\in\mathcal{F}_{\beta'}\), there is also a neighbourhood \(V\) of \(\bs{x}\) contained in \(\operatorname{conv}{(\Delta_2)}\). As a result, \(V\subset \Delta_2\cup s_0\) and \(V\subset \Delta_2\cup r_0\), as both \(s_2\), and \(r_2\) are convex. Finally, \(V\cap s_0=V\cap r_0\), therefore \(V\cap s_1=V\cap (s_0\cup\Delta_1) = V\cap (r_0\cup\Delta_1)=V\cap r_1\), and we can use the assumption that \(r_1\in\mathcal{F}_{\beta'}\). \qed

\medskip
\noindent
\textbf{Proof of \textup{Theorem \ref{thm:t18}}}
	The theorem follows from statement C. of Lemma \ref{lem:l19}, Lemma \ref{lem:l21}, and the observation that \(\operatorname{int}{(s_1)}\supset\operatorname{int}{(s_0)}\neq\emptyset\). \qed

\medskip
\noindent
\textbf{Proof of \textup{Theorem \ref{thm:t17}}}
	For an arbitrary soul \((\bar{t},t_1)\in\tau_{\beta'}\), with the same disturbance \(t_1\) and any \(l\in\{1,\ldots,n\}\), it is enough to show that \(s_1\coloneqq (\bar{t}\sm t_1)\cup\left(t(l)\sm t(1)\right)\in\mathcal{F}_{\beta'}\). Let us denote \(r_0=t\sm t_1\), \(r_1=t(l)\), \(r_2=t\), \(s_0=\bar{t}\sm t_1\), \(s_2=\bar{t}\). This way we obtain a pack \((r_0,r_1,r_2,s_0,s_1,s_2)\), and we can use Theorem~\ref{thm:t18}. \qed

\medskip
\noindent
\textbf{Proof of \textup{Theorem \ref{thm:astast16}}}
	The theorem follows by a straightforward application of Theorem \ref{thm:t17}. \qed

\subsubsection*{C. On concepts of strict convexity}

In the construction of nets demanded in condition (\textasteriskcentered\textasteriskcentered) of Theorem \ref{thm:astast16}, we will use an equivalent definition of the class \(\mathcal{F}_\varrho\) of \(\varrho\)-convex sets, for the fixed set \(\mathbb{T}\)

We will fix the number \(\varrho\geq 1\) and use the following notation.
\begin{itemize}
	\item[] We set \(\mathcal{F} = \left\{t\in\mathcal{T}\colon t\subset\bbt\right\}\).
	\item[] For a convex, closed set \(s\subset\mathbb{R}^2\), and \(\bs{x}\in\mathbb{R}^2\sm s\), \(P_s\bs{x}\) will be the minimal distance projection: \(P_s\bs{x}=\bs{y}\), if \(\bs{y}\in s\) and \(\|\bs{x}-\bs{y}\| = \min{\left\{\|\bs{x}-\bs{z}\|\colon \bs{z}\in s\right\}}\);
	\item[] we also write \(B(s,\bs{x})=B\), if \(B\in \mathcal{B}_\varrho\) and \(P_B\bs{x} = P_s\bs{x}\).
	\item[] If additionally \(s\in\mathcal{F}\) and \(\bs{x}\in\bbt\sm s\), we define \(E(s,\bs{x})\) to be the connected component of \(\bbt\sm B(s,\bs{x})\) containing \(\bs{x}\) as a member.
	\item[] We set \(C(s,\bs{x})=\bbt\sm E(s,\bs{x})\), if \(E(s,\bs{x})\) is defined.
	\item[] Lastly, we denote
	\begin{multline*}
	\mathcal{C}_\varrho = \Big\{\bbt\sm E\colon E\text{ is a connected component of }\bbt\sm B\text,\\\text{ for some }B\in\mathcal{B}_\varrho
	\text{, } \operatorname{int}(\bbt\cap B)\neq\emptyset
	\text{, }\bbt\sm B\neq\emptyset\Big\}
	\text.
	\end{multline*}
\end{itemize}

\begin{thm}\label{thm:t24}
	For any set \(t\in\mathcal{F}\), the following conditions are equivalent:
	\begin{itemize}
		\item[\(\alpha\)\textup)]	\(t\in\mathcal{F}_\varrho\),
		\item[\(\beta\)\textup)]	\(t=\bigcap_{C\in\mathcal{D}\cup\{\bbt\}}C\), for some family \(\mathcal{D}\subset\mathcal{C}_\varrho\),
		\item[\(\gamma\)\textup)]	\(t\subset C(t,\bs{x})\), for every \(\bs{x}\in\bbt\sm t\).
	\end{itemize}
\end{thm}

Our proof of the theorem will follow from a sequence of lemmas concerning the conditions \(\alpha\)), \(\beta\)), \(\gamma\)).

\begin{lem}\label{lem:l25}
	For any \(t\in\mathcal{F}\), condition \(\gamma\)) implies \(\beta\)), and condition \(\beta\)) implies \(\alpha\)).
\end{lem}
\medskip
\noindent
\textbf{Proof.}
	To show that \(\gamma)\implies\beta)\), it is enough to take \(\mathcal{D}=\left\{C(t,\bs{x})\colon\bs{x}\in\bbt\sm t\right\}\).
	
	To prove that \(\beta)\implies\alpha)\), consider \(t\neq\bbt\). We define a function \(f\colon\mathcal{C}_\varrho\to\mathcal{B}_\varrho\), continuous with respect to the Hausdorff metric, by the relation: \(f(\bbt\sm E)=B\), when \(E\) is a connected component of \(\bbt\sm B\). Then, for \(\bs{x}_0\in\partial t\cap\operatorname{int}{(\bbt)}\) we can find sets \(C_1,C_2,\ldots\in\mathcal{D}\) satisfying \(\bs{x}_0\in\partial\bigcap_{i\in\mathbb{N}}C_i\).
	By a subsequence argument, we can assume that \(d\left(f(C_i), B_0\right)\to 0\), when \(i\to\infty\), for some ball \(B_0\in\mathcal{B}_\varrho\) satisfying \(\bs{x}_0\in\partial B_0\). This, in turn, implies the \(\varrho\)-convexity of \(t\) at \(\bs{x}_0\). \qed

\medskip

The proof of the implication \(\alpha)\implies\gamma)\) is somewhat more laborious.

\begin{lem}\label{lem:l26}
	If \(t\in\mathcal{F}_\varrho\), \(\bs{x},\bs{y}\in t\), \(\bs{x}\in\operatorname{int}{(\bbt)}\), then we have \([\bs{x},\bs{y}]\sm\{\bs{x},\bs{y}\}\subset\operatorname{int}{(t)}\).
\end{lem}

\medskip

\noindent
\textbf{Proof.}
	In view of convexity of \(\bbt\), the assumption \(\bs{x}\in\operatorname{int}{(\bbt)}\) implies \(I\subset\operatorname{int}{(\bbt)}\) with \(I=[\bs{x},\bs{y}]\sm\{\bs{x},\bs{y}\}\).
	For any \(\bs{z}\in I\) and any ball \(B\in\mathcal{B}_\varrho\) satisfying \(\partial B\ni \bs{z}\), we have \(V\cap I\not\subset B\), for any neighborhood \(V\) of \(\bs{z}\).
	Thus, supposition of \(\bs{z}\in\partial t\) contradicts the assumption that \(t\in\mathcal{F}_\varrho\). \qed

\medskip

\noindent
\begin{lem}\label{lem:l27}
	If a set \(t\in\mathcal{T}\), for some ball \(B\in\mathcal{B}_\varrho\) and point \(\bs{x}\in\partial B\cap\partial t\), satisfies the conjunction
	\[
	\exists{\eps_1>0}\colon B(\bs{x},\eps_1)\cap B\subset t
	\quad\text{and}\quad
	\forall{\eps>0}\colon B(\bs{x},\eps)\cap t\not\subset B
	\text,
	\]
	then \(t\) cannot be \(\varrho\)-convex at \(\bs{x}\).
\end{lem}

\medskip
\noindent
\textbf{Proof.}
	Let \(\eps_1\) satisfy \(B(\bs{x},\eps_1)\cap B\subset t\). For any ball \(B_1\in\mathcal{B}_\varrho\) and any number \(\eps\in(0,\eps_1)\), we have the following two implications. If \(B_1\neq B\), \(\bs{x}\in\partial B_1\), then \(B(\bs{x},\eps)\cap B\not\subset B_1\), so \(B(\bs{x},\eps)\cap t\not\subset B_1\). If \(B_1=B\), then \(B(\bs{x},\eps)\cap t\not\subset B_1\) by assumptions. Since \(\eps\in(0,\eps_1)\) has been arbitrarily chosen, the set \(t\) is not \(\varrho\)-convex at \(\bs{x}\). \qed

\medskip

For \(\bs{x},\bs{y}\in\bbt\), \(\bs{x}\neq\bs{y}\), we call \(l\) a \emph{bow with endpoints \(\bs{x},\bs{y}\)}, if \(l\) is an arc of some circle \(\partial B\), \(B\in\mathcal{B}_\varrho\), with endpoints \(\bs{x},\bs{y}\), and \(d(l)=\|\bs{x}-\bs{y}\|\). Thus, we always have \(d(l)\leq 1\), as \(d(\bbt)\leq 1\).

\begin{lem}\label{lem:l28}
	For any \(t\in\mathcal{F}_\varrho\) and any bow \(l\) with endpoints \(\bs{x},\bs{y}\), if \(\bs{x},\bs{y}\in t\) and \(\bs{x}\in\operatorname{int}{(\bbt)}\), \(l\subset\bbt\), then \(l\subset t\).
\end{lem}
\medskip
\noindent
\textbf{Proof.}
	Changing coordinate system if necessary, we have
	\[
	l=\left(\lvert x\rvert\leq\alpha,\ y=\sqrt{\varrho^2-x^2}-\sqrt{\varrho^2-\alpha^2}\right)
	\text,
	\]
	for some \(\alpha\in(0,1)\).
	Let us denote
	\[
	l_\lambda=\left(\lvert x\rvert\leq\lambda,\ y=\sqrt{\varrho^2-x^2}-\sqrt{\varrho^2-\lambda^2}\right)
	\text,
	\]
	for \(\lambda\in[0,\alpha]\).
	Assume that \(l\not\subset t\). Since \(l=l_\alpha\), we have \(l_\lambda\not\subset t\) for \(\lambda\) close to \(\alpha\).
	On the other hand, by Lemma \ref{lem:l26}, we have \(\bs{0}\in\operatorname{int}{(t)}\), thus \(l_\lambda\subset t\) for \(\lambda\) close to \(0\).
	Let \(\lambda_0\) be the highest number in \((0,\alpha)\) for which \(l_{\lambda_0}\subset t\). Then, the set \(k\coloneqq l_{\lambda_0}\cap\partial t\) is nonempty and closed.
	Let us consider point \(\bs{x}_0=(x_0,y_0)\in k\) of the smallest possible value of \(x_0\).
	In view of Lemma \ref{lem:l26}, the line segment
	\(\left(\lvert x\rvert<\alpha, y=0\right)\), is contained in \(\operatorname{int}{(t)}\), therefore the endpoints of \(l_{\lambda_0}\) are not elements of \(k\).
	We have then \(\bs{x}_0\in\operatorname{int}{(\mathbb{T})}\), by \(\lambda_0\in(0,\alpha)\), moreover, for the ball \(B_0\) satisfying \(l_{\lambda_0}\subset\partial B_0\), we have
	\[
	\forall{\eps>0}\colon B(\bs{x}_0,\eps)\cap t\not\subset B_0
	\quad\text{and}\quad
	\exists{\eps>0}\colon  B(\bs{x}_0,\eps)\cap B_0\subset t
	\text.
	\]
	According to Lemma \ref{lem:l27} we have obtained a contradiction with the assumption that \(t\in\mathcal{F}_\varrho\), which completes the proof. \qed

\begin{lem}\label{lem:l29}
	Let \(t\in\mathcal{F}_\varrho\), \(\bs{x}\in\bbt\sm t\). We then have
	\(\partial^-E(t,\bs{x})\cap\operatorname{int}{(t)}=\emptyset\).
\end{lem}
\medskip
\noindent
\textbf{Proof.}
	Suppose that there exists a point \(\bs{z}\in\partial^-E(t,\bs{x})\cap\operatorname{int}{(t)}\).
	Changing coordinate system if necessary, we can assume that \(\bs{x}=(0,\chi)\), for some \(\chi>0\), \(P_t\bs{x}=\bs{0}\), \(B(t,\bs{x})=B((0,-\varrho),\varrho)\), \(\bs{z}=(\xi,\sqrt{\varrho^2-\xi^2}-\varrho)\), for some \(\xi\in(0,1)\).
	
	For \(\eps>0\) let us denote line segment \(I_\eps\coloneqq\left(x=0,0<y\leq\eps\right)\).
	The arc
	\[l=\left(0\leq x\leq\xi,y=\sqrt{\varrho^2-x^2}-\varrho\right)\text,\] with endpoints \(\bs{0}\) and \(\bs{z}\), satisfies \(l\subset \partial^-E(t,\bs{x})\) with \(E(t,\bs{x})\supset I_\eps\cup(\bs{z}+I_\eps)\), for some \(\eps>0\).
	Then, by connectedness of \(E(t,\bs{x})\), we have \(l+I_{\eps_1}\subset E(t,\bs{x})\), for some \(\eps_1\in(0,\eps)\).
	
	Now we can see that in the strip \(\clo(l+I_{\eps_1})\) we can find bows with endpoints \(\bs{0}\) and \(\bs{v}\), which we may denote \(l_{\bs{v}}\), satisfying \(\bs{v}\in\bs{z}+I_{\eps_1}\).
	Adjusting \(\eps_1\) to be suitably small, we can assure that \(\bs{z}+I_{\eps_1}\subset\operatorname{int}(t)\).
	By Lemma \ref{lem:l28}, we have \(l_{\bs{v}}\subset t\), for \(\bs{v}\in\bs{z}+I_{\eps_1}\).
	Lastly, let us note that for those arcs we also have \(l_{\bs{v}}\not\subset(y\leq 0)\), which contradicts the fact that \(P_t(0,\chi)=\bs{0}\), with some \(\chi>0\). \qed

\medskip
\noindent
\textbf{Proof of \textup{Theorem \ref{thm:t24}}}
	In view of Lemma \ref{lem:l25}, we will show that \(\alpha)\implies\gamma)\).
	For \(t\in\mathcal{F}_\varrho\) let us suppose that there exists \(\bs{z}\in E(t,\bs{x})\cap t\).
	Then, \(B(\bs{z},\eps)\cap B(t,\bs{x})=\emptyset\), for some \(\eps>0\). Remark \ref{rem:elem_tcapU} yields existence of a ball \(B(\bs{z}_1,\eps_1)\subset B(\bs{z},\eps)\cap\operatorname{int}{(t)}\).
	As a result, \(B(\bs{z}_1,\eps_1)\) and \(\{\bs{z}\}\) are contained in the same connected component of \(\bbt\sm B(t,\bs{x})\), and we have \(B(\bs{z}_1,\eps_1)\subset E(t,\bs{x})\cap t\).
		
	In a suitable coordinate system we have \(\bs{x}=(0,\chi)\), \(P_t\bs{x}=\bs{0}\), \(t\subset(y\leq 0)\), thus \(\bs{z}_1\in(y<0)\), \(B(t,\bs{x})=B((0,-\varrho),\varrho)\) and \(\bs{z}_1\notin B((0,-\varrho),\varrho)\).
	Therefore, there exists a point \(\bs{y}\in\partial B((0,-\varrho),\varrho)\cap \left([\bs{0},\bs{z}_1]\sm \{\bs{0},\bs{z}_1\}\right)\).
	Then, \(\bs{y}\in\operatorname{int}{(t)}\), as \(\bs{z}_1\in\operatorname{int}(t)\),
	and \(\bs{y}\in\partial^-E(t,\bs{x})\), as \(\bs{z_1}\in E(t,\bs{x})\), \(\bs{0}\in\partial^- E(t,\bs{x})\), with \(E(t,\bs{x})\) being a connected component of \(\bbt\sm B((0,-\varrho),\varrho)\). We have obtained a contradiction with Lemma \ref{lem:l29}, which ends the proof. \qed

\subsubsection*{D. Constructions of nets}

Our proof of the condition (\textasteriskcentered\textasteriskcentered) in Theorem \ref{thm:astast16} will be based on the following two theorems.

\begin{thm}\label{thm:t22}
	There are numbers \(\delta >0\) and \(j\in\mathbb{N}\), independent on $\mathbb{T}$, such that for any soul \((t,t_1)\in\tau_{\beta'/2}\) satisfying \(d(t,t\sm t_1)\leq \delta\),
	the family \(\mathcal{N}_j^{(t,t_1)}\) is not empty.
\end{thm}

\begin{thm}\label{thm:t23}
	For any \(\delta>0\) there is a number \(k\in\mathbb{N}\), independent on $\mathbb{T}$, such that for any soul \((t,t_1)\in\tau_{\beta'/2}\), there are sets \(t(i)\in\mathcal{F}_{\beta'/2}\), \(1\leq i\leq k\), satisfying \(t\sm t_1=t(1)\subset\ldots\subset t(k)=t\) and \(d(t(i),t(i-1))\leq\delta\), for \(2\leq i\leq k\).
\end{thm}

Using these theorems it is easy to see the validity of the following.
\begin{thm}\label{thm:t33}
	Condition \textup{(\textasteriskcentered\textasteriskcentered)} of \textup{Theorem \ref{thm:astast16}} is satisfied.
\end{thm}
\medskip
\noindent
\textbf{Proof.}
	Let us fix numbers \(\delta>0\), \(j\in\mathbb{N}\) as in Theorems \ref{thm:t22}. For the number \(\delta\) we choose \(k\in\mathbb{N}\) as in Theorem \ref{thm:t23}. Then we set \(n^*=(k-1)j\). In particular numbers \(\delta\), \(j\), \(k\), and \(n^*\) do not depend on \(\bbt\). With this notation, for any soul \((t,t_1)\in\tau_{\beta'/2}\), we can find sets \(t(1),\ldots,t(k)\in\mathcal{F}_{\beta'/2}\), satisfying \(t\sm t_1=t(1)\subset\ldots\subset t(k)=t\), \(d(t(i),t(i+1))\leq \delta\), for \(2\leq i\leq k\).
	Next, for \(i\in\{2,\ldots,k\}\) we can find nets \(\bs{t}(i)\in\mathcal{N}_j^{(t(i),t(i)\sm t(i-1))}\). Their concatenation \(\bs{t}=\bs{t}(2)\circ\ldots\circ\bs{t}(k)\) is a member of \(\mathcal{N}_{n^*}^{(t,t_1)}\). \qed

\medskip

The proof of Theorem \ref{thm:fund}, given as a consequence of Theorem \ref{thm:withast} in part A. of this section, is thus completed, as we have the following.

\medskip

\noindent
\textbf{Proof of \textup{Theorem \ref{thm:withast}}}
	According to Theorem \ref{thm:t33}, condition (\textasteriskcentered\textasteriskcentered) is satisfied. By Theorem \ref{thm:astast16}, so is the condition (\textasteriskcentered). \qed

\medskip

Our proofs of Theorems \ref{thm:t22} and \ref{thm:t23} are based on the global characterization of the class \(\mathcal{F}_\varrho\) given in Theorem \ref{thm:t24}.

Proceeding to the proof of Theorem \ref{thm:t22}, we fix number \( \bar{\varrho} > \varrho \) and introduce some additional notation.

For any \(\bs{z}\in\mathbb{R}^2\), \(s\subset\mathbb{R}^2\), \(s\neq\emptyset\), we set \(d(\bs{z},s) = \inf{\left\{\|\bs{z}-\bs{y}\|\colon \bs{y}\in s\right\}}\).

\begin{namelist}{lll}
	\item{$\bullet$} With \(\delta(\varrho,\bar{\varrho})\) we denote the smallest positive number for which conditions \(\|\bs{z}\|\geq\gamma'/3-\delta(\varrho,\bar{\varrho})\), \(d(\bs{z},B((0,-\varrho),\varrho))\leq\delta(\varrho,\bar{\varrho})\), for \(\bs{z}\in\mathbb{R}^2\), imply \(\bs{z}\in B((0,-\bar{\varrho}),\bar{\varrho})\).
\end{namelist}


For \(s\in\mathcal{T}\) and \(\bs{x}\in\bbt\sm s\), analogically to the sets \(B(s,\bs{x})\), \(E(s,\bs{x})\), \(C(s,\bs{x})\), we define:
\begin{itemize}
	\item[] \(\bar{B}(s,\bs{x})=B\), if \(B\in\mathcal{B}_{\bar{\varrho}}\), \(P_B\bs{x}=P_s\bs{x}\),
	\item[] \(\bar{E}(s,\bs{x})\), to be the connected component of \(\bbt\sm\bar{B}(s,\bs{x})\) containing \(\bs{x}\),
	\item[] \(\bar{C}(s,\bs{x})=\bbt\sm\bar{E}(s,\bs{x})\), thus
	\begin{equation}\label{eq:e9}
	\bar{C}(s,\bs{x})\supset C(s,\bs{x})\text,\quad\text{for }\bs{x}\in\bbt\sm s\text.
	\end{equation}
	\item[] We obviously put
	\begin{multline*}
		\mathcal{C}_{\bar{\varrho}} = \Big\{\bbt\sm E\colon E\neq\emptyset\text{ is a connected component of }\bbt\sm B \\
	\text{for some ball }B\in\mathcal{B}_{\bar{\varrho}}
	\text{, }\operatorname{int}(\bbt\cap B)\neq\emptyset\text{, }\bbt\sm B\neq \emptyset\Big\} \text.
	\end{multline*}
\end{itemize}

Let us note that Theorem \ref{thm:t24} is also valid when \(\varrho\) is replaced with \(\bar{\varrho}\), therefore, for any set \(t\in\mathcal{T}\), the following conditions are equivalent:
\begin{itemize}
	\item[\(\bar{\alpha}\)\textup)]	\(t\in\mathcal{F}_{\bar{\varrho}}\),
	\item[\(\bar{\beta}\)\textup)]	\(t=\bigcap_{C\in\mathcal{D}\cup\{\bbt\}}C\), for some family \(\mathcal{D}\in\mathcal{C}_{\bar{\varrho}}\).
\end{itemize}

We will use the following observation.

\begin{lem}\label{lem:l30}
	For \(t\in\mathcal{F}_\varrho\), \(\bs{x}\in\bbt\sm t\), \(\bs{z}\in E(t,\bs{x})\), we have \(d(\bs{z},B(t,\bs{x})) \leq d(\bs{z}, C(t,\bs{x}))\).
\end{lem}
\medskip
\noindent
\textbf{Proof.}
	It is enough to notice that \(\|\bs{z}-\bs{v}\|\geq d(\bs{z}, B(t,\bs{x}))\), for \(\bs{v}\in C(t,\bs{x})\cap F\), for some connected component \(F\) of \(\bbt\sm B(t,\bs{x})\) that is disjoint with \(E(t,\bs{x})\). For such points \(\bs{v}\) we have \([\bs{z},\bs{v}]\cap B(t,\bs{x})\neq\emptyset\), which completes the proof. \qed

\begin{lem}\label{lem:l31}
	For sets \(r\subset t\), \(\Delta\subset t\sm r\), satisfying \(r,t\in\mathcal{F}_\varrho\), \(d(r,t)\leq \delta(\varrho,\bar{\varrho})\), \(d(\Delta)\leq\gamma'/3\), there exists a set \(\bar{r}\in\mathcal{F}_{\bar{\varrho}}\) satisfying \(r\cup\Delta\subset\bar{r}\subset t\) and \(d(\bar{r}\sm r)\leq\gamma'\).
\end{lem}
\medskip
\noindent
\textbf{Proof.}
	Denote \(\sigma=\gamma'/3\).
	In the case when \(\Delta=\emptyset\), we set \(\bar{r}=r\). If \(\Delta\neq\emptyset\) and \(d(t\sm r,\Delta)\leq\sigma\), we use the relation \(\mathcal{F}_\varrho\subset\mathcal{F}_{\bar{\varrho}}\).
	We take \(\bar{r}=t\), and the bound \(d(\bar{r}\sm r)\leq 3\sigma=\gamma'\) follows from the triangle inequality. Otherwise, we set
	\[
	\bar{r}=t\cap \bigcap_{\substack{\bs{x}\in t\sm r \\ d(\bs{x},\Delta)>\sigma}} \bar{C}(r,\bs{x}) \text.
	\]
	
	Since \(t\in\mathcal{F}_\varrho\subset\mathcal{F}_{\bar{\varrho}}\), and \(\bar{\beta}\textup{)}\implies\bar{\alpha}\textup{)}\), we have $\bar{r}\in\mathcal{F}_{\bar{\varrho}}$ and \(d(\bs{x},\Delta)\leq\sigma\), for \(\bs{x}\in\bar{r}\sm r\), thus \(d(\bar{r}\sm r)\leq 3\sigma=\gamma'\), by triangle inequality. In accordance with \eqref{eq:e9}, we also have \(r\subset\bar{r}\).
	
	The proof will be completed by showing that \(\Delta\subset\bar{r}\). For any \(\bs{x}\in t\sm r\) satisfying \(d(\bs{x},\Delta)>\sigma\) and its projection \(\bs{y}=P_r\bs{x}\), we have \(d(\bs{y},\Delta)>\sigma-\delta(\varrho,\bar{\varrho})\).
	Moreover, for \(\bs{z}\in\Delta\) we have
	\begin{equation*}
		\begin{aligned}
			d(\bs{z},B(r,\bs{x})) &\leq d(\bs{z},C(r,\bs{x}))\quad&\text{(cf. Lemma \ref{lem:l30})}\\
			&\leq d(\bs{z},r)\quad &\text{(as \(r\subset C(r,\bs{x})\))}\\
			&\leq \delta(\varrho,\bar{\varrho})
			\quad&\text{(since \(\bs{z}\in\Delta\subset t\))}\text.
		\end{aligned}
	\end{equation*}
	In a suitable coordinate system we then have \(\bs{y}=\bs{0}\), \(B(r,\bs{x})=B((0,-\varrho),\varrho)\), \(\bar{B}(r,\bs{x})=\bar{B}((0,-\bar{\varrho}),\bar{\varrho})\), and
	\begin{gather*}
		\|\bs{z}\|=\|\bs{z}-\bs{y}\|\geq d(\bs{y},\Delta)>\sigma-\delta(\varrho,\bar{\varrho})\text, \\
		d(\bs{z},B((0,-\varrho),\varrho))\leq \delta(\varrho,\bar{\varrho})
		\text.
	\end{gather*}
	Therefore, \(\bs{z}\in\bar{B}(r,\bs{x})\), by the definition of \(\delta(\varrho,\bar{\varrho})\). Moreover, \(\bs{z}\in\bar{C}(r,\bs{x})\), by the fact that \(\bs{z}\in\Delta\subset\bbt\). The inclusion \(\Delta\subset\bar{r}\) is thus shown. \qed

\medskip

In the forthcoming proof we use the assumption that \(\bbt\subset[-1,1]\times[-1,1]\).

\medskip
\noindent
\textbf{Proof of \textup{Theorem \ref{thm:t22}}}
	We can find a number \(j\in\mathbb{N}\) such that there exists a covering \(\bigcup_{2\leq h\leq j}\Delta_h\) of the square \([-1,1]\times[-1,1]\) satisfying \(d(\Delta_h)\leq\gamma/3\), for \(2\leq h\leq j\). For this \(j\) we further fix numbers \(\varrho(h)\) such that \(\beta'/2=\varrho(1)<\ldots<\varrho(j)=\beta'\), and define \(\delta = \min{\left\{\delta(\varrho(h-1),\varrho(h))\colon 2\leq h\leq j\right\}}\).
	Clearly, \(j\) and \(\delta\) do not depend on \(\bbt\).
	
	Let \((t,t_1)\in\tau_{\beta'/2}\) be a soul. We set \(t(1)=t\sm t_1\). Having defined \(t(h-1)\), for \(2\leq h\leq j\), satisfying \(t\sm t_1\subset t(h-1)\subset t\), \(t(h-1)\in \mathcal{F}_{\varrho(h-1)}\), let us denote \(r=t(h-1)\), \(\varrho=\varrho(h-1)\), \(\bar{\varrho}=\varrho(h)\), and \(\Delta=\Delta_h\cap(t\sm t(h-1))\). Then, \(\delta(\varrho,\bar{\varrho})\geq\delta\geq d(t,r)\), \(d(\Delta)\leq\gamma'/3\), and we can apply Lemma \ref{lem:l31}. For thus obtained \(\bar{r}\), we set \(t(h)=\bar{r}\), which also assures that \(t\sm t_1\subset t(h)\subset t\) and \(t(h)\in \mathcal{F}_{\varrho(h)}\).
	The resulting sequence of sets \(t(1)\subset\ldots\subset t(j)\) satisfies
	\(t(1)=t\sm t_1\), \(t(h)\in \mathcal{F}_{\varrho(h)}\subset\mathcal{F}_{\beta'}\), \(d(t(h)\sm t(h-1))\leq \gamma'\), for \(2\leq h\leq j\) and \(t(j)\supset(t\sm t_1)\cup \bigcup_{2\leq h\leq j}(\Delta_h\cap t_1)=t\). This ends the proof. \qed

\medskip

Now we proceed to the proof of Theorem \ref{thm:t23}. To that end we fix number \(\varrho=\beta'/2\) and for \(\delta>0, \bs{x},\bs{y}\in\bbt, \bs{x} \neq \bs{y}\), we define
\begin{itemize}
	\item[$\bullet$]
	\(
	\eps(\delta) = \inf{
		\left\{
		\|\bs{\bar{z}}-\bs{z}\|\colon \bs{\bar{z}}\in B((0,-\varrho),\varrho)\text{, }\bs{z}\in B(\bs{0},1)\sm B((0,\delta-\varrho),\varrho)
		\right\}
	}
	\),
	\item[]
	\( B(\{\bs{y}\},\bs{x})  = B\)\quad if \(B\in\mathcal{B}_\varrho\),\quad
	\(P_B\bs{x}=\bs{y}\),
	\item[]
	$C(\{\bs{y}\}, \bs{x}) = \bbt \setminus E$ if $E$ is connected component of the set $\bbt \setminus B(\{\bs{y}\},\bs{x})$ and $\bs{x} \in E$.
\end{itemize}

In the forthcoming proof we use the assumption that
\(d(\bbt)\leq 1\), which somewhat simplifies the argument.

\begin{lem}\label{lem:l32}
	For \(r\in\mathcal{F}_\varrho\), \(\delta>0\), there exists set \(s\in\mathcal{F}_\varrho\) satisfying \(r\subset s\), \(d(r,s)\leq\delta\), and the implication
	\(
	\bs{x}\in\bbt\sm r\text{, }d(\bs{x},r)\leq\eps(\delta)\implies \bs{x}\in s
	\).
\end{lem}
\medskip
\noindent
\textbf{Proof.}
	If \(d(\bbt,r)\leq \delta\), it is enough to take \(s=\bbt\). Assume then that \(d(\bbt,r)> \delta\).
	For any \(\bs{x}\in\bbt\) satisfying \(d(\bs{x},r)>\delta\), we can find \(\bs{y}\in[P_r\bs{x},\bs{x}]\) for which \(\|\bs{y}-P_r\bs{x}\|=\delta\) (thus \(\bs{y}\neq\bs{x}\) and \(C(\{\bs{y}\},\bs{x})\supset C(r,\bs{x})\)).
	
	For \(\bs{z}\in\bbt\) such that \(d(\bs{z},r)\leq \eps(\delta)\) we show that \(\bs{z}\in C(\{\bs{y}\},\bs{x})\).
	
	If \(\bs{z}\notin E(r,\bs{x})\), then \(\bs{z}\in C(r,\bs{x})\subset C(\{\bs{y}\},\bs{x})\).
	
	However, if \(\bs{z}\in E(r,\bs{x}), d(\bs{z}, r) \leq \eps(\delta)\), we can use Lemma \ref{lem:l30} to obtain inequality \(d(\bs{z},B(r,\bs{x}))\leq\eps(\delta)\).
	With suitably adjusted coordinate system we have  \(P_r\bs{x}=\bs{0}\), \(\bs{y}=(0,\delta)\), \(B(r,\bs{x})=B((0,-\varrho),\varrho)\), \(B(\{\bs{y}\},\bs{x})=B((0,\delta-\varrho),\varrho)\).
	Suppose that \(\bs{z}\in\bbt\sm C(\{\bs{y}\},\bs{x})\).
	As \(d(\bbt)\leq 1\), we have both \(\bs{z}\in B(\bs{0},1)\sm B((0,\delta-\varrho),\varrho)\) and \(d(\bs{z},B((0,-\varrho),\varrho))\leq\eps(\delta)\). Contradiction with the definition of \(\eps(\delta)\) shows that \(\bs{z}\in\bbt\cap B(\{\bs{y}\},\bs{x})\subset C(\{\bs{y}\},\bs{x})\).
	
	Lastly, setting \(s=\bigcap_{\substack{\bs{x}\in \bbt \\ d(\bs{x},r)>\delta}}C(\{\bs{y}\},\bs{x})\), in view of the implication \(\beta\text)\implies\alpha\text)\) in Theorem \ref{thm:t24}, completes the proof. \qed

\medskip

\noindent
\textbf{Proof of \textup{Theorem \ref{thm:t23}}}
	Let us fix \(k\geq 1/{\eps(\delta)}+1\). For \((t,t_1)\in\tau_{\beta'/2}\), let us denote \(t(1)=t\sm t_1\), and having defined \(t(i)\), for \(i\in\{1,\ldots,k-1\}\), we apply Lemma \ref{lem:l32} with \(r=t(i)\) to obtain \(s\) described therein. We set \(t(i+1)=s\cap t\). Then, \(d(t(i),t(i+1))\leq\delta\).
	
	Next, let us consider \(\bs{x}\in t\), for which we choose a point \(\bs{x}_1\in t\sm t_1\), together with points \(\bs{x}_i\in[\bs{x}_1,\bs{x}]\) satisfying \(\|\bs{x}_i-\bs{x}_{i-1}\|<\eps(\delta)\), for \(2\leq i\leq k\), \(\bs{x}_k=\bs{x}\). By induction, \(\bs{x}_i\in t(i)\) for \(1\leq i\leq k\), therefore \(\bs{x}\in t(k)=t\). Finally, \(t(k)=t\), by the arbitrariness of \(\bs{x}\in t\), and \(k\) clearly does not depend on \(\bbt\). \qed

\medskip

Finally, the proof of Theorem \ref{thm:fund} is completed.

\section{Existence of convex Peano curve filling convex compact set $\mathbb{T}$}

Recall that $\T$ is a family of convex, compact subsets of $\R^2$ with non-empty interiors. We can now demonstrate that the announced Fundamental Construction for any $\TT \in \T$ is feasible. A sequence of sets $\ttt = (t(K); 1 \leq K \leq M)$ forms a population if $t([K,L]) \in \T$ for $1 \leq K \leq L \leq N$ (for unions $t([K,L]) = \bigcup_{K \leq J \leq L} t(J)$). We also denote (as previously) $\bigcup \ttt = t([1,M])$, $d (\ttt) = \text{max} \{ d(t(K)): 1 \leq K \leq M \}$, $I \ttt = (t(M+1 - K); 1 \leq K \leq M)$ and $I^0 \ttt = \ttt, I^n \ttt = I(I^{n-1} \ttt)$ for $n \in \N$. For a sequence of sequences of sets $(\ttt(K); 1 \leq K \leq M)$ we denote anti-ordering $\mathcal{A}(\ttt(K); 1 \leq K \leq M) = I^0 \ttt(1) \circ I^1 \ttt(2) \circ \dots \circ I^{M-1} \ttt(M)$. (cf. (\ref{eq3.2.1})).

The existence of the Fundamental Construction for $\TT$ is stated by the following theorem:

\begingroup
\setcounter{thm}{1}
\renewcommand\thethm{\Alph{thm}}
\begin{thm}
\label{thmB}
For $\TT \in \T$, there exist numbers $M(j), m(j) \in \N$ and populations of sets $\ttt(j) = (t(j; K); 1 \leq K \leq M(j)) \in \T^{M(j)}$, for $j \in \N$, satisfying the following conditions:
\begin{enumerate}[1)]
\item[\textup{1)}] $\bigcup \ttt(1) = \TT$;
\item[\textup{2)}] for $j \in \N$ we have $M(j+1) = M(j)m(j)$ and $t(j;K) = t(j+1; [(K-1)m(j) +1, Km(j)])$ for $K \in \{ 1, \dots, M(j) \}$;
\item[\textup{3)}] $ d(\ttt(j)) \rightarrow 0$ as $j \rightarrow \infty$.
\end{enumerate}
\end{thm}
\endgroup
\setcounter{thm}{0}

\textbf{Proof.} Without loss of generality, we assume $d(\TT) \leq 1$, $\TT \subset [-1,1] \times [-1,1]$. We will base our argument on Theorem \ref{thm3.3.1} and Corollary \ref{cor3.3.2}. We fix numbers $\gamma(j) \in (0,\frac{1}{4})$ for $j \in \N$ such that $\gamma(j) \rightarrow 0$ as $j \rightarrow \infty$ and $\beta(1) \geq 1$, $\beta(j+1) = \text{max} \{ 2 \beta(j), \gamma(j) + \frac{1}{\gamma(j)} \}$ for $j \in \N$. Then we fix numbers $m'(j) \in 2 \N$ for $j \in \N$ satisfying the following conditions:
\begin{enumerate}[I)]
\item For odd $j$, there exists offspring of souls
\begin{equation}
\label{eq3.4.1}
\tau_{\beta(j) \gamma(j)} \ni (t,t_1) \rightarrow \mathopen{<}\ttt'(j;t,t_1), \ttt'_1(j;t,t_1)\mathclose{>} \in (\tau_{\beta(j+1), \gamma(j+1)})^{m'(j)}
\end{equation}
satisfying $d_x(\ttt'(j; t,t_1)) \leq 11 \gamma(j)$, in notation \eqref{eq:3.3.1}, according to Theorem \ref{thm3.3.1}, for the fixed \(\bbt\).

\item For even $j$, there exists offspring of souls (\ref{eq3.4.1}) satisfying $d_y(\ttt'(j;t,t_1)) \leq 11 \gamma(j)$, in notation \eqref{eq:3.3.3}, according to Corollary \ref{cor3.3.2}.
\end{enumerate}

We define numbers $M(j)$ and sequences of souls $\mathopen{<} \ttt(j), \ttt_1(j)\mathclose{>} \in (\tau_{\beta(j) \gamma(j)})^{M(j)}$ for $j \geq 1$ according to the following induction. For $j=1$, we set $M(1) = 2$ and $\ttt(1) = (\TT, \TT), \ttt_1(1)= (\emptyset, \emptyset)$. Assume now that a sequence of souls
\begin{equation*}
\mathopen{<}\ttt(j), \ttt_1(j)\mathclose{>} = ((t(j;K), t_1(j;K)); 1 \leq K \leq M(j)) \in (\tau_{\beta(j)\gamma(j)})^{M(j)}
\end{equation*}
is defined. For $1 \leq K \leq M(j)$ denote by $\mathopen{<} \ttt'(j;K), \ttt'_1(j;K) \mathclose{>}$ the offspring of a parent $(t(j;K), t_1(j;K))$, according to I) for $j \in 2\N -1$ and according to II) for $j \in 2\N$. Then, using anti-ordering of a sequence of sequences of sets, we define $\mathopen{<}\ttt(j+1), \ttt_1(j+1)\mathclose{>} = \mathopen{<}\mathcal{A}(\ttt'(j;K); 1 \leq K \leq M(j)), \mathcal{A}(\ttt'_1(j;K); 1 \leq K \leq M(j))\mathclose{>}$. Thus also we have \\ $\mathopen{<}\ttt(j+1), \ttt_1(j+1)\mathclose{>} \in (\tau_{\beta(j+1) \gamma(j+1)})^{M(j+1)}$.

Notice that the proposed two-term sequence $\mathopen{<}\ttt(1), \ttt_1(1)\mathclose{>} = ((\TT, \emptyset), (\TT, \emptyset))$ is a population of souls in $\tau_{\beta(1) \gamma(1)}$, according to Definition \ref{df3.1.2}. According to Theorem \ref{thm3.2.3}, the sequences $\mathopen{<}\ttt(j), \ttt_1(j)\mathclose{>}$ are also populations of souls, recursively for $j = 2,3,\dots$. The obtained bases $\ttt(j)$ are therefore populations of sets, according to Theorem \ref{thm3.1.6}. Condition 1) is fulfilled, as well as condition 2) (according to (\ref{eq3.2.3})).

According to notations \eqref{eq:3.3.1} and \eqref{eq:3.3.3}, we also have
\(d_x(\ttt(1)) \geq d_x(\ttt(2)) \geq \ldots\) and \(d_y(\ttt(1))\geq d_y(\ttt(2))\geq \ldots\), together with the inequalities \(d_x(\ttt(j+1))\leq 11\gamma(j)\), for \(j\in 2N-1\), and \(d_y(\ttt(j+1))\leq 11\gamma(j)\), for \(j\in 2N\),
thus condition 3) is satisfied as well. \qed

\medskip
\noindent
\textbf{Proof of Theorem \ref{thm0}.} The case when int$(\TT) = \emptyset$ is trivial. When int$(\TT) \neq \emptyset$, we apply Theorems \ref{thm2.0.1} and \ref{thmB}. \qed

\medskip

The relationships between convex Peano curves will be investigated in our next publication. Here, we will confine ourselves to the following remark concerning affine mappings.

Let's consider any convex Peano curve $f$ filling a given compact convex set $\TT \subset \R^2$, and denote the function $g: \R \rightarrow \TT$ as follows:
\begin{equation*}
g(u) = \left\{
\begin{array}{ccc}
f(0) & \text{for} & u < 0 \\
f(u) & \text{for} & u \in [0,1] \\
f(1) & \text{for} & u > 1.
\end{array}
\right.
\end{equation*}

Next, for a convex set $K$ in a real or complex linear topological space $\XX$ and a continuous functional $\Lambda$ on $\XX$, assume there exist numbers $a, b \in \Lambda(K)$ such that $a \neq b$. Without loss of generality, one can assume that $\operatorname{Re}a < \operatorname{Re}b$. Then the function $h: \XX \rightarrow \TT$, defined as $h(\x) = g\left(\frac{\operatorname{Re}(\Lambda(\x) - a)}{\operatorname{Re}(b-a)}\right)$, is a continuous surjection, where the images $h(L)$ are convex for all convex subsets $L \subset K$. Thus, we have:

\begin{cor}
\label{cor3.4.2}
If for a convex set $K$ in a real or complex linear topological space $\XX$, there exists a continuous linear functional on $\XX$ that is non-constant on $K$, then, for a compact convex set $\TT \subset 
\R^2$, there exists a continuous surjection $h: \XX \rightarrow \TT$, where the images of convex subsets $L$ of $K$ are convex, and $h(K) = \TT$.
\end{cor}

It is therefore evident that convex mappings can differ drastically from affine mappings.

\bigskip

\begin{center}{\centering \textsc{Acknowledgments}}\end{center}

Sincere gratitude is due to Professor Jalanta Misiewicz, whose assistance has fundamentally influenced the shape of this work. I would also like to thank Professor Andrzej Komisarski for valuable discussions, as well as Professors Stanisław Goldstein and Władysław Wilczyński for their valuable information on the literature. Additional thanks go to Małgorzata Kozioł, Jakub Olejnik, and Damian Prusinowski for their technical and linguistic assistance, and for their patience.


\newpage
\newpage


\bigskip

\textbf{Adam Paszkiewicz}, Faculty of Mathematics and Computer Science, University of Lodz, Poland,
email: \texttt{adam.paszkiewicz@wmii.uni.lodz.pl}
\end{document}